\documentclass[11pt, reqno]{amsart}

\usepackage{hyperref} 
\usepackage{amssymb, mathrsfs,bbold}
\usepackage{amsmath, amsthm}
\usepackage{enumerate}
\usepackage{dsfont}
\usepackage{color}
\usepackage[a4paper]{geometry}
\usepackage[utf8]{inputenc}   
\usepackage{stmaryrd}
\usepackage{titlesec, wasysym}
\usepackage{appendix}
\usepackage{todonotes, fancyhdr}
\usepackage{chngcntr}
\usepackage{cancel}

\usepackage{tikz-cd}
\usetikzlibrary{calc}

\usepackage{setspace}
\renewcommand{\baselinestretch}{0.99}

\usepackage[all]{xy}
\numberwithin{subsection}{section}
\numberwithin{subsubsection}{subsection}
\numberwithin{equation}{section} 

\newenvironment{Dem}[1][\unskip]{%
    \begin{list}{\hspace{1.15cm}\textbf{\textit{Proof #1 --}}}{   
        \setlength{\topsep}{0pt}%
        \setlength{\leftmargin}{0pt}%
        \setlength{\rightmargin}{0pt}%
        \setlength{\listparindent}{0pt}%
        \setlength{\itemindent}{0pt}%
        \setlength{\parsep}{0pt}%
        \addtolength{\leftmargin}{0pt} 
        \addtolength{\rightmargin}{0pt}%
    } \item }{\hfill $\RHD$\end{list}\smallskip}

\newenvironment{Dem*}[1][\unskip]{%
    \begin{list}{\hspace{0cm}{\sf \textbf{{\small Proof #1 --}}}}{   %
        \setlength{\topsep}{0pt}%
        \setlength{\leftmargin}{0pt}%
        \setlength{\rightmargin}{0pt}%
        \setlength{\listparindent}{0pt}%
        \setlength{\itemindent}{0pt}%
        \setlength{\parsep}{0pt}%
        \addtolength{\leftmargin}{20pt}%
        \addtolength{\rightmargin}{0pt}%
    } \item }{\hfill $\rhd$\end{list}\smallskip}

\renewcommand\thesection       {\arabic{section}}
\renewcommand\thesubsection    {\thesection{\boldmath $.$}\arabic{subsection}}
\renewcommand\thesubsubsection    {\thesection{\boldmath $.$}\arabic{subsection}{\boldmath $.$}\arabic{subsubsection}} 

\titleformat{\section}[block] 
{\filcenter\normalfont\sffamily\bfseries\large}  
{{\hspace{-0.87cm}}\thesection \hspace{0.2em} --\vspace{0cm}}{0.5em}{} 

\titleformat{\subsection}[runin]
{\filcenter\normalfont\sffamily\bfseries}  
{{\hspace{0cm}}\thesubsection \hspace{0.15em} -- \vspace{0.1cm}}{.2em}{}   
\titlespacing{\subsection}{-0pc}{1.5ex plus .1ex minus .2ex}{0pc}   

\titleformat{\subsubsection}[runin]
{\filcenter\normalfont\sffamily\bfseries}   
{\filright\sffamily{\hspace{0cm}}\thesubsubsection\hspace{0.2em} --}{.5em}{}\titlespacing{\subsection}{-0pc}{1.5ex plus .1ex minus .2ex}{0pc}

\newtheoremstyle{mystyle}
{3pt}               
{3pt}               
{\it }                      
{}                      
{\bfseries}      
{}                      
{0.5em}                 
{\hspace{0cm}\textit{#2 --} {\hspace{-0.02cm}}\textit{#1}}    

\theoremstyle{mystyle}

\newtheorem{thm}{Theorem.}   
\newtheorem*{thm*}{Theorem}

\newtheorem{lem}[thm]{{Lemma}. }
\newtheorem{prop}[thm]{{Proposition.}}
\newtheorem{defn}[thm]{{Definition.}}
\newtheorem*{rem*}{Remark.}

\newtheoremstyle{mystyle3}
{3pt}               
{3pt}               
{\it }                      
{}                      
{\bfseries}      
{}                      
{0.5em}                 
{\hspace{0cm}{\textbf{\textit{#2}} --} {\hspace{-0.02cm}}{\textbf{\textit{#1}}}}   

\theoremstyle{mystyle3}

\newtheoremstyle{mystyle2}
{3pt}               
{3pt}               
{\it }                      
{}                      
{\sffamily}    
{}                      
{0.5em}                 
{\llap{#2 }{\it #1{\hspace{0.2cm}--}}}

\theoremstyle{mystyle2}

\newtheorem*{definition*}{Definition}
\newtheorem*{theorem*}{Theorem}
\newtheorem*{Remark*}{Remark}
\newtheorem*{lem*} {Lemma}
\newtheorem*{defn*} {Definition}
\newtheorem*{prop*} {Proposition}
\newtheorem*{cor*} {Corollary}

\newcommand{\ssk}{\smallskip}

\renewcommand{\epsilon}{\varepsilon}
\newcommand{\eps}{\epsilon}

\newcommand\bbE{\mathbb{E}}

\newcommand\bbN{\textbf{\textsf{N}}}
\newcommand\bbR{\textbf{\textsf{R}}}

\newcommand{\mcD}{\mathcal{D}}
\newcommand{\mcE}{\mathcal{E}}
\newcommand{\mcF}{\mathcal{F}}

\newcommand{\mcL}{\mathcal{L}}

\newcommand{\mcN}{\mathcal{N}}

\newcommand{\mcP}{\mathcal{P}}

\newcommand\mcW{\mathcal W}

\newcommand{\bfZ}{\mathbf{Z}}

\newcommand{\bfu}{{\bf u}}
\newcommand{\bfv}{{\bf v}}

\newcommand{\bbL}{\mathbb{L}}
\newcommand{\bbP}{\mathbb{P}}

\newcommand{\scrC}{\ensuremath{\mathscr{C}}}

\makeatletter
\newcommand*{\defeq}{\mathrel{\rlap{%
                     \raisebox{0.3ex}{$\m@th\cdot$}}%
                     \raisebox{-0.3ex}{$\m@th\cdot$}}%
                     =}
\makeatother

\makeatletter
\newcommand*{\eqdef}{=\mathrel{\rlap{%
                     \raisebox{0.3ex}{$\m@th\cdot$}}%
                     \raisebox{-0.3ex}{$\m@th\cdot$}}%
                     }
\makeatother

\newcommand{\dd}{\text{d}}
\newcommand{\dab}{\mcE^\beta_T}
\newcommand{\norme}[1]{\Vert #1 \Vert}

\newcommand{\dw}{d_{L_T^\infty\mathcal{W}_{p,C^\alpha}}}

\newcommand{\scrL}{\mathscr{L}}
\definecolor{greeen}{rgb}{0.0,0.56,0.15}

\begin{document}

\begin{center}
{\LARGE\sffamily{\textbf{Mean field singular stochastic PDEs}   \vspace{0.5cm}}}
\end{center}

\begin{center}
{\sf I. BAILLEUL and N. MOENCH}
\end{center}

\vspace{1cm}

\begin{center}
\begin{minipage}{0.8\textwidth}
\renewcommand\baselinestretch{0.7} \scriptsize \textbf{\textsf{\noindent Abstract.}} We study some systems of interacting fields whose evolution is given by some singular stochastic partial differential equations of mean field type. We provide a robust setting for their study and prove a well-posedness result and a propagation of chaos result. 
\end{minipage}
\end{center}

\vspace{0.6cm}

{\sf 
\begin{center}
\begin{minipage}[t]{11cm}
\baselineskip =0.35cm
{\scriptsize 

\center{\textbf{Contents}}

\vspace{0.1cm}

\textbf{1. Introduction\dotfill 
\pageref{SectionIntro}}

\textbf{2. Additive noise\dotfill
\pageref{SectionAdditiveNoise}}

\textbf{3. Basics on paracontrolled calculus\dotfill 
\pageref{SectionBasics}}

\textbf{4. Mean field equations with a regular interaction in the diffusivity\dotfill 
\pageref{SectionRegularInteraction}}

\textbf{5. Mean field equations with pointwise interaction in the diffusivity\dotfill 
\pageref{SectionGeneral}}

\textbf{6. Propagation of chaos\dotfill 
\pageref{SectionChaos}}

\textbf{A. Enhancing some random noises\dotfill 
\pageref{SectionEnhancingNoises}}

\textbf{B. Proof of Theorem \ref{thm_nonexplosionA}\dotfill 
\pageref{section_appendix_nonexplosion}}

}\end{minipage}
\end{center}
}   \vspace{1cm}

\section{Introduction}
\label{SectionIntro}

Let $(\xi^i)_{i\geq 1}$ stand for a sequence of independent, identically distributed, random spacetime distributions on the $2$-dimensional torus $\mathbf{\mathsf{T}}^2$ defined on a common probability space $(\Omega, \mcF, \mathbb{P})$. Assume in this introduction that the $\xi^i$ are almost surely some continuous functions of time with values in the space of $(\alpha-2)$-H\"older regular distributions over $\mathbf{\mathsf{T}}^2$, for some $2/3<\alpha<1$. The archetype of such a noise is given by the time-independent space white noise. We study a system of interacting fields whose evolution is given by the following system of {\it singular} stochastic partial differential equations (SPDEs)
\begin{equation} \label{EqMainSystem}
(\partial_t-\Delta)u^i_t = f(u^i_t,\mu_t^n)\,\xi_t^i + g(u^i_t,\mu_t^n), \qquad (1\leq i\leq n),
\end{equation}
where 
\[
\mu_t^n \defeq \frac{1}{n}\sum_{i=1}^n\delta_{u^i_t}
\]
is the running time empirical measure of the system, a probability measure on some appropriate function space. Some (possibly random) initial conditions in that function space are given. We write here $u_t^i$ for $u^i(t)$. 

The first works on mean field type systems of stochastic partial differential equations are due to Chiang, Kallianpur \& Sundar \cite{CKS} and Kallianpur \& Xiong \cite{KX}. They consider the case of a constant diffusivity $f=1$ and a spacetime white noise $\xi$ and relate the corresponding equation to some models of large assemblages of interacting neurons. The dynamics of interacting polymer systems can also be modelled by this special kind of equations \cite{EShen, Criens}. All these and subsequent works consider the case of a spacetime white noise $\xi$. This allows to use the tools of stochastic/It\^o calculus to make sense of the equations considered and study them. However the relevance of spacetime white noise as a model of real life noise is questionable in a number of situations. The present work is a first step to building a robust approach of this kind of systems that applies to a large class of noises that goes beyond the semimartingale class of noises of stochastic calculus. We take profit for that purpose of the tools that were developed for the study of the so-called singular stochastic partial differential equations.

\ssk

Recall the rule of thumb: {\sl One can make sense of the product of two distributions with some given H\"older regularities if and only if the sum of their regularity exponents is positive.} The term {\sl singular} in the expression {\sl singular stochastic partial differential equation} refers to the fact that, given the low regularity of the noise, the regularizing effect of the heat resolvent does not give a priori a sufficient regularity to the $u^i$ to make sense of the products $f(u^i,\mu_t^n)\,\xi^i$ in \eqref{EqMainSystem} -- neither as a space product at any fixed time or as a product of parabolic distributions/functions. The diffusivity term $f(u^i,\mu_t^n)$ is expected to have at best parabolic regularity $\alpha$ while the product $f(u^i,\mu_t^n)\,\xi^i$ is well defined if and only if $\alpha+(\alpha-2)>0$. This condition does not hold in our case where $\alpha<1$. 

The settings of regularity structures and paracontrolled calculus have been developed in the last ten years to deal precisely with this kind of problem and one can indeed use either of them to make sense of Equation \eqref{EqMainSystem} as an equation of the form 
\begin{equation} \label{EqReformulationSystem}
(\partial_t-\Delta) {\sf u}^{(n)} = {\sf f}({\sf u}^{(n)})\,\xi^{(n)} + {\sf g}({\sf u}^{(n)}),
\end{equation}
for some $n$-dimensional unknown ${\sf u}^{(n)} = (u^1,\dots,u^n)$ and some multi-dimensional noise 
\[
\xi^{(n)} = \big(\xi^1,\dots,\xi^n\big),
\]
and to identify some conditions on the functions ${\sf f}$ and ${\sf g}$ under which \eqref{EqReformulationSystem} has a unique solution over a given time interval. This way of proceeding does not take profit from the specific structure of the mean field type equation \eqref{EqMainSystem}. It is in particular unclear how to prove a propagation of chaos result for the interacting field system from this point of view. The necessity of a point of view tailor-made to mean field-type dynamics gets even clearer if one looks at what should most naturally be the limit dynamics of a fixed field in the system \eqref{EqMainSystem} when $n$ tends to $\infty$, say the field with label $i=1$. Based on symmetry/exchangeability considerations, this field is expected to be a solution of the equation
\begin{equation} \label{EqMainEq}
(\partial_t-\Delta)u_t = f(u_t, \mcL(u_t)) \xi_t + g(u_t , \mcL(u_t)),
\end{equation}
where $\mcL(u_t)$ stands for the law of the random variable $u_t=u(t)$ and $\xi$ stands for a random distribution with the same law as the $\xi^i$. The dynamics \eqref{EqMainEq} comes with an initial condition $u_0$ in some function space. This type of equations appears as the large scale picture of all the above mentioned works. One has for instance in the study of the large $N$ linear sigma model of quantum field theory the equation
\begin{equation} \label{EqLinearSigma}
(\partial_t-\Delta) \Phi = -\bbE[\Phi^2] \, \Phi + \xi,
\end{equation}
where $\xi$ stands here for a spacetime white noise. This equation was studied in depth by Shen, Smith, Zhu \& Zhu in \cite{SZZ22-3D, SSZZ22}. The additive character of this equation makes it possible to use some tools that are not available in the more general setting of the dynamics \eqref{EqMainEq} with multiplicative noise. Stochastic fluid models, or simplified climate models, also provide some motivation for studying mean field type equations of the form \eqref{EqMainEq}. See for instance the recent works of Ren, Tang \& Wang \cite{RTW} and Crisan, Holm \& Korn \cite{CHK} and the references therein.

\ssk

Our first aim in this work is to develop a setting within which one can make sense of the system \eqref{EqMainSystem} and the equation \eqref{EqMainEq} in a unified way, for a large class of spacetime noises $\xi$, for some mean field-dependent non-constant diffusivity $f$.

\ssk

Denote by $z$ and $z'$ some generic spacetime points. The choice of some particular class of functions $f$ and $g$ in \eqref{EqMainSystem} and \eqref{EqMainEq} is guided by the physics of the phenomenon modeled by the system \eqref{EqMainSystem}. The phenomenology of the linear sigma model leads for instance to \eqref{EqLinearSigma}, with its quadratic, but signed, diffusivity $-\bbE[\Phi^2]$. To make things concrete we consider in this introduction the case where $f(u,\mu)$ and $g(u,\mu)$ depend linearly on their measure argument and are of the form
\begin{equation} \label{EqModelDiffusivity} \begin{split}
f(u,\mu)(z) = \iint F\big(u(z), v(z')\big) k(z,z') \, dz' \mu(dv) = \bbE\bigg[\int F\big(u(z),V(z')\big)k(z,z') \, dz'\bigg]
\end{split} \end{equation}
for $u$ a real-valued function on $\mathbf{\mathsf{T}}^2$, for a random function $V$ defined on $(\Omega,\mcF,\bbP)$, with law $\mu$, and a real-valued function $F$ on $\bbR^2$. Think of the kernel $k$ as a parameter that captures the range of the interaction between the different fields in the system, with extreme cases $k(z,z')=1$ and $k(z,z')=\delta_z(z')$, and intermediate cases represented by some bounded $C^2$ kernels, with bounded derivatives, for instance. The physics behind the two extreme cases is very different and we will technically deal with them in a different way. For $k\in C^2$ we talk of a \emph{regular interaction}. For $k(z,z')=\delta_z(z')$ we talk of a \emph{pointwise interaction}. Our main result reads informally as follows.

\ssk

\begin{thm} \label{ThmMainSketch}
One can design a setting where Equation \eqref{EqMainEq} makes sense.
\begin{itemize}
	\item[{\it (a)}] Under some Lipscthiz type regularity assumptions on $f$ and $g$ there exists a deterministic positive time $T$ such that the system \eqref{EqMainSystem} and Equation \eqref{EqMainEq} have some unique solutions on the time interval $[0,T]$.    \vspace{0.1cm}
	
	\item[{\it (b)}] The law of any fixed finite tuple of fields in the field system \eqref{EqMainSystem} converges to a tuple of independent, identically distributed, solutions of \eqref{EqMainEq} as $n$ tends to $\infty$, on the time interval $[0,T]$. 
\end{itemize}
\end{thm}

\ssk

So there is propagation of chaos for the system \eqref{EqMainSystem}, with mean field dynamics given by the mean field type equation \eqref{EqMainEq}. See Theorem \ref{SolMcKeanSingular} and Theorem \ref{ThmPropagationChaos} for some proper versions of Theorem \ref{ThmMainSketch}. We note here that Theorem \ref{ThmMainSketch} applies to some nonlinear functions $f,g$ of their measure argument. We note as well that one can prove some much stronger results in the case of an additive noise ($f=1$). The equation is in particular well-posed globally in time. See Theorem \ref{ThmAffine} in Section \ref{SectionAdditiveNoise}.

\ssk

{\it \S1. Paracontrolled approach.} While Equation \eqref{EqMainEq} and the system \eqref{EqMainSystem} share the common feature of being singular, in the sense that they involve some ill-defined products, the mean field interaction in \eqref{EqMainEq} causes a different kind of problem. A close situation was studied by Bailleul, Catellier \& Delarue in their analysis of mean field type random rough differential equations \cite{BailleulCatellierDelarue}. They used therein the language of controlled paths. We design in the present work an approach similar to \cite{BailleulCatellierDelarue,roughchaos} for the study of \eqref{EqMainEq}, using the language of paracontrolled calculus to build our setting. 

The original form of paracontrolled calculus was introduced by Gubinelli, Imkeller \& Perkowski in \cite{GIP}; one can find a nice short account of the basics of paracontrolled calculus in Gubinelli \& Perkowski's lecture notes \cite{GPEnsaios}. The paracontrolled structure involves an operator $\varolessthan$ that will be introduced below. Denote by $\omega\in\Omega$ a generic chance element and write $X(\omega)$ for $-(\partial_t-\Delta)^{-1}(\xi(\omega))$, and $\overline{X}$ for an independent copy of the random variable $X$, both defined on the probability space $(\Omega\times\Omega,\mcF\otimes\mcF,\bbP^{\otimes 2})$. We use a notion of paracontrolled field that is tailor made to capture not only the paracontrolled structure of $u$ needed to make sense of its product with $\xi$ but also of the structure needed to describe the mean field specific spacetime function 
$$
(t,x)\mapsto f\big(u_t,\mcL(u_t)\big)(x).
$$ 
This comes under the form of a definition saying that a random field $\phi(\omega)$ is $\omega$-paracontrolled by a reference field $X(\omega)$ of parabolic H\"older regularity $\alpha$ if one has almost surely
\begin{equation} \label{EqDefnOmegaPC}
\phi(\omega) \approx (\delta_z \phi)(\omega) \varolessthan X(\omega) + \overline{\bbE}\big[(\delta_\mu \phi)(\omega,\cdot) \varolessthan \overline{X}(\cdot)\big]
\end{equation}
up to a remainder with better regularity, for some random functions $(\delta_z \phi)(\omega)$ and $(\delta_\mu \phi)(\omega,\varpi)$ on $\mathbf{\mathsf{T}}^2$ that depend on $\omega$ and an additional independent chance element $\varpi$ that is averaged out in the $\overline{\bbE}$ expectation. A precise definition, conveying in particular the meaning of the notations $\delta_z\phi, \delta_\mu \phi$, is given in Section \ref{SubsectionParacStructure}. This definition will play a key role in our construction of a robust solution theory for \eqref{EqMainEq}. {\color{black} One can think of them as the classical and some generalized Gubinelli derivatives of $\phi$ with respect to $X$ and $\overline{X}$.}

\ssk

{\it \S2. Enhanced noise and renormalization.} Setting up a framework for the study of a given singular stochastic partial differential equation driven by a realization $\xi(\omega)$ of a random noise usually requires that we enhance the noise with the additional datum of some quantities that do not make sense analytically $\omega$-wise. In the archetypal example of the $2$-dimensional parabolic Anderson model equation we are given a space white noise $\xi$ that is almost surely of space H\"older regularity $-1-\eta$ for all $\eta>0$. The dynamics is
\begin{equation} \label{EqPAMIntro}
(\partial_t-\Delta)v = v\xi.
\end{equation}
Enhancing the noise consists in building a random variable $\xi^{(2)}$ that plays the role of the $\omega$-wise ill-defined product of $\xi(\omega)$ and $(\partial_t-\Delta)^{-1}(\xi(\omega))$. This random variable is given by the $\bbL^2(\Omega,\mathbb{P})$ limit of the renormalized regularized quantity
\[
\xi^\epsilon \, (\partial_t-\Delta)^{-1}(\xi^\epsilon) - C^\epsilon,
\]
where $\xi^\epsilon$ stands for an arbitrary smooth regularization of $\xi$ by convolution that converges to $\xi$ in the space of distributions with H\"older regularity $-1-\eta$, and $C^\epsilon$ is an explicit constant that diverges to $+\infty$ as a multiple of $\vert\hspace{-0.05cm}\log\epsilon\vert$. This limit random variable is suggestively denoted by $(\xi X)(\omega)$. The fact that the naive approximation $\xi^\epsilon(\omega) (\partial_t-\Delta)^{-1}(\xi^\epsilon(\omega))$ is not converging leads to the interpretation of the solution $v$ to \eqref{EqPAMIntro} as a limit (in probability) of solutions $v^\epsilon$ to the renormalized equation
\[
(\partial_t-\Delta)v^\epsilon = v^\epsilon\xi^\epsilon - C^\epsilon v^\epsilon,
\]
rather than as a limit of solutions to the parabolic Anderson model equation driven by the regularized noise $\xi^\epsilon$. We talk in this setting of the pair of random variables $(\xi, \xi^{(2)})$ as an {\sl enhanced noise}. 

A richer enhancement of the noise $\xi$ is needed in the analysis of the mean field equation \eqref{EqMainEq}. Not only do we need to add the random variable $(\xi X)(\omega)$ to our notion of enriched noise, but the description \eqref{EqDefnOmegaPC} of an $\omega$-controlled field should make it plain that we also need to add a doubly random variable that plays the role of the analytically ill-defined product of $\xi(\omega)$ and $(\partial_t-\Delta)^{-1}(\overline{\xi}(\varpi))$, where $(\omega,\varpi)\in\Omega^2$ and we work with the product probability $\mathbb{P}^{\otimes 2}$ on $(\Omega^2,\mcF^{\otimes 2})$. Luckily, the independence of $\xi$ and $\overline{\xi}$ allows to define a doubly random variable $\big(\xi \overline{X})(\omega,\varpi)$ as the $L^2(\Omega^2,\mathbb{P}^{\otimes 2})$ limit of the regularized quantity

\[
\xi^\epsilon(\partial_t-\Delta)^{-1}(\overline{\xi}^\epsilon)
\]
{\it without any renormalization}. This will lead us to the interpretation of a solution to Equation \eqref{EqMainEq} as the limit in probability as $\epsilon>0$ goes to $0$ of the solution $u^\epsilon$ to the renormalized equation
\[
(\partial_t-\Delta)u^\epsilon = f(u^\epsilon,\mcL(u^\epsilon_t))\,\xi^\epsilon_t - C^\epsilon (f'f)\big(u^\epsilon,\mcL(u^\epsilon_t)\big) + g(u^\epsilon,\mcL(u^\epsilon_t)),
\]
where $f'$ stand for the derivative of $f$ with respect to its first argument and $C^\epsilon$ is a constant diverging as $\epsilon>0$ goes to $0$. See Theorem \ref{SolMcKeanSingular} for a proper statement.

\bigskip

\textbf{\textit{Organization of this work.}} We treat the elementary case of a system \eqref{EqMainSystem} and Equation \eqref{EqMainEq} with additive noise ($f=1$) in Section \ref{SectionAdditiveNoise}. Some very robust results can be obtained in this simple setting, leading in particular to a simple proof of propagation of chaos for the corresponding system of interacting fields for an essentially arbitrary random noise $\xi$ of the form $(\partial_t-\Delta)Z$, for $Z\in C_TC^\alpha$. This setting covers the case where $\xi$ is a spacetime white noise. No tools from paracontrolled calculus are needed to deal with this case. 

We use the language of paracontrolled calculus to study some more general equations or systems. We recall what we need from this domain in Section \ref{SectionBasics} and study Equation \eqref{EqMainEq} in the simple setting of a diffusivity of the form \eqref{EqModelDiffusivity} with a $C^2$ kernel $k$ in Section \ref{SectionRegularInteraction}. {\color{black} The classical tools of paracontrolled calculus are sufficient to make sense of \eqref{EqMainEq} and prove a local in time well-posedness result.}

The setting of Section \ref{SectionBasics} is not sufficient to deal with \eqref{EqMainEq} {\color{black} when the interaction kernel is the Dirac $\delta_x(y)$ kernel. The interacting fields have in that case a purely pointwise interaction.} We introduce a notion of mean field enhancement of the noise in Section \ref{SubsectionEnhancedNoise}, and describe an associated notion of paracontrolled structure in Section \ref{SubsectionParacStructure}. This structure specific to the mean fiel setting allows to give a fine description of some function of the form $f(u,\mcL(u))$. The well-posed character of \eqref{EqMainEq} in an ad hoc space of paracontrolled functions is the object of Section \ref{SubsectionSolvingEquation}. The quantitative regularity result that we obtain for the solution $u$ of \eqref{EqMainEq} as a function of the enhanced noise entails in Section \ref{SectionChaos} a propagation of chaos result for the system \eqref{EqMainSystem}.

We gather in some appendices some side results. Appendix \ref{SectionEnhancingNoises} provides a pedestrian approach of the probabilistic construction of the enhancement of a relatively large class of noises. The mean field nature of \eqref{EqMainEq} requires that we deal with some non-exploding dynamics. The global well-posedness problem for classical singular stochastic PDEs is a tricky problem that was solved recently in some generality in the regime that we consider. Our proof of Theorem \ref{ThmMainSketch} involves a variation of some a priori estimates of Shen, Zhu \& Zhu \cite{ShenZhuZhu} whose proof is detailed in Appendix \ref{section_appendix_nonexplosion}.

\ssk

For a general overview of the domain of propagation of chaos for mean field type dynamics we recommend the review articles \cite{ChaintronDiez1, ChaintronDiez2} by Chaintron \& Diez.

\medskip

\textbf{\textit{Notations.}} We gather here a number of notations that we will use frequently. We write $a \defeq b$ to define $a$ as being equal to $b$.
\begin{itemize}	
	\item[--] {\it For $\gamma\in\bbR\backslash\bbN$ the Besov-H\"older spaces $C^\gamma$ over $\mathbf{\mathsf{T}}^2$ and their norms $\Vert\cdot\Vert_\gamma$ are defined from some Littlewood-Paley projectors $\Delta_i : \mcD'(\mathbf{\mathsf{T}}^2)\rightarrow C^\infty(\mathbf{\mathsf{T}}^2)$ setting 
	\[
	\Vert f\Vert_\gamma \defeq \sup_{i\geq -1} 2^{i\gamma} \Vert \Delta_i(f)\Vert_\infty
	\]
	and defining $C^\gamma$ as the closure of $C^0$ under this norm. Set
	\[
	\Delta_{<j} \defeq \sum_{i\leq j-1}\Delta_i,
	\]
	and define the paraproduct $f\varolessthan g$ of any two distributions $f,g$ on $\mathbf{\mathsf{T}}^2$ as
	\begin{equation} \label{EqParaproduit}
	f \varolessthan g \defeq  \sum_{i\geq 1} \Delta_{<i-1}(f)\,\Delta_i(g).
	\end{equation}
	The resonant operator $f\odot g$ is defined formally as
	\begin{equation} \label{EqResonance}
	f\odot g \defeq \sum_{\vert i-j\vert\leq 1} \Delta_j(f)\,\Delta_i(g).
	\end{equation}
	For any Banach space $E$, regularity exponent $r>0$ and time horizon $0<T<\infty$ we set 
	\[
	C_T^r E \defeq C^r([0,T],E)
	\] 
	and write $L^\infty_TE$ for $L^\infty([0,T];E)$. We will also need the parabolic H\"older space $\scrC^\alpha_T$ on $[0,T]\times\mathbf{\mathsf{T}}^2$, for $\alpha>0$, which is isometric to $C^{\alpha/2}_TL^\infty(\mathbf{\mathsf{T}}^2) \cap C_T C^\alpha(\mathbf{\mathsf{T}}^2)$ equipped with its natural norm. For $u \in \scrC^\alpha_T$ we denote by $u_t$ its value at time $t$ -- this is a function on $\mathbf{\mathsf{T}}^2$.   \vspace{0.1cm}
	
	\item[--] We will denote by $(P_t)_{t \geq 0}$ the semigroup generated by the Laplace-Beltrami operator $\Delta$ on an ad hoc function space. Recall the elementary estimate
	\begin{equation} \label{EqHeatEstimate}
	\norme{P_t(u)}_{C^{\gamma+\delta}}\lesssim_T \min(t,1)^{-\delta/2}\norme{u}_{C^\gamma},
	\end{equation}
	for any $\gamma\in\bbR$ and $\delta\geq 0$. We write 
	\[
	\mathscr{L} \defeq \partial_t-\Delta
	\]
	for the heat operator and define the inverse operator 
	\[
	\mathscr{L}^{-1}(\zeta)(t) \defeq \int_0^t P_{t-s}(\zeta_s) \, ds.
	\]
	We have for any $0<T<\infty$ the continuity estimate
\begin{equation} \label{EqBoundIntegration}
\Vert \scrL^{-1}(\zeta)\Vert_{C_TC^\alpha} \lesssim_T \Vert \zeta\Vert_{C_T C^{\alpha-2}}.
\end{equation}	}
	
	\item[--] {\it We denote by $\bbL^p(\Omega \,; E)$ the space of $E$-valued random variables in $\bbL^p(\Omega,\mcF,\bbP)$. We denote by $\mathcal{L}(Z)$ the law of a random variable $Z$.}   \vspace{0.1cm}
	
	\item[--] {\it For an integrability exponent $1\leq p<\infty$ we denote by $\mcP_p(E)$ the set of probability measures on the metric space $(E,d)$ that have a moment of order $p$. We also denote by $\mcW_{p,E}$ the $p$-Wasserstein metric on $\mcP_p(E)$ 
	
	\[
	\mcW_{p,E}(P,Q) = \inf_\gamma\left\{\left(\int_E d(x,y)^p\,\gamma(dx,dy)\right)^{1/p}\,;\,\gamma\in\mcP_0(E\times E) \textrm{ has marginals }P \textrm{ and }Q\right\}.
	\]
	For $\alpha>0$ we define a distance on $L^\infty_T\mathcal{P}_p(C^\alpha)$ setting
	\[
	\dw(\mu,\mu') \defeq \sup_{t\in[0,  T]} \mathcal{W}_{p,C^\alpha}(\mu_t,\mu'_t).
	\]   }
	
	\item[--] {\it For a measure $\mu$ on a metric space $E$ and $\phi\in C_b(E)$ we write $\mu(\phi)$ for the integral of $\phi$ with respect to $\mu$.}   \vspace{0.1cm}

	\item[--] {\it We fix throughout this work some regularity exponents}
	\begin{equation} \label{EqConditionAlphaBeta}
	\frac{2}{3}<\beta<\alpha<1.
	\end{equation}	
\end{itemize}

\medskip

\section{Additive noise}
\label{SectionAdditiveNoise}

Fix $0<T_0<\infty$ and $1\leq p<\infty$. For $Z$ a function on $[0,T_0]\times \mathbf{\mathsf{T}}^2$ we set
\[
\zeta = \mathscr{L}(Z).
\] 
Following Coghi, Deuschel, Friz \& Maurelli \cite{CDFM} we begin our work by studying the case of a mean field type equation with additive noise 
\begin{equation} \label{EqMcKV0}
(\partial_t-\Delta)u = \zeta + g(u, \mcL(u_t))
\end{equation}
and random initial condition $u_0$, assuming that the random variable $(Z,u_0)$ is an element of $\bbL^p\big(\Omega, C_{T_0}C^r \times C^r\big)$ with $Z$ null at time $0$, {\color{black} for some regularity exponent $r>0$}. This setting includes the case where $\zeta$ is a spacetime white noise. We rewrite \eqref{EqMcKV0} in integral form
\begin{equation} \label{EqMcKV0Integral}
u_t = P_t(u_0) + Z_t + \int_0^t P_{t-s} \big(g(u_s,\mcL(u_s))\big) \, ds  \qquad  (0\leq t\leq T_0).
\end{equation}
No singular product is involved in the study of this equation and we are able to solve it with some classical elementary tools. We prove in Section \ref{SubsectionAdditive} that Equation \eqref{EqMcKV0Integral} is well-posed if $g$ is Lipschitz continuous in the sense of \eqref{EqAssumptionG} below. The law of the solution to \eqref{EqMcKV0Integral} turns out to be a Lipschitz continuous function of the law of $(Z,u_0)$ in the Wasserstein $p$-space. This strong result leads in Section \ref{SubsectionAdditiveChaos} to a propagation of chaos result for an associated system of interacting fields.

\medskip

\subsection{Additive mean field equation$\boldmath{.}$   \hspace{0.15cm}}
\label{SubsectionAdditive} 

For $\mu\in\mathcal{P}_p(C_{T_0}C^r)$ and $t\in[0,T_0]$, we write $\mu_t$ for the image measure of $\mu$ in $C^r$ by the $t$-time coordinate map $u\in C_{T_0}C^r\mapsto u_t \in C^r$. We make the following regularity assumption on $g$.

\medskip

\noindent \textbf{\textsf{Assumption on $g$ --}} {\it There exists a constant $L$ such that for every $v_1, v_2\in C^r$ and $\nu_1,\nu_2 \in \mcP_p(C^r)$ we have}
\begin{equation} \label{EqAssumptionG} \begin{split} 
\big\Vert g(v_1, \nu_1) - g(v_2,\nu_2)\big\Vert _{C^{r-2}}  \leq L\big( \norme{v_1-v_2} _{C^r} + \mathcal{W}_{p,C^r}(\nu_1,\nu_2) \big).
\end{split} \end{equation}

\medskip

{\color{black} For instance, for a globally Lipschitz function $G : C^r\times C^r\rightarrow C^{r-2}$ the functions
\[
g(v,\mu) = \int G(v,w)\mu(dw)\quad \textrm{or} \quad G\Big(v,\int \hspace{-0.1cm}w\mu(dw)\Big)
\]
satisfy the Lipschitz bound \eqref{EqAssumptionG}. The particular case $g(v,\mu)=c\int\hspace{-0.1cm}w\mu(dw)$, for some constant $c$, is already interesting.} The Schauder estimate \eqref{EqBoundIntegration} ensures that the map
\[
\Phi(u) = P_t(u_0) + Z_t + \int_0^t P_{t-s} \big(g(u_s,\mu_s)\big) ds
\]
is a well defined map from $C_{T_0} C^r$ into itself. 

\ssk

\begin{prop}
Suppose Assumption \eqref{EqAssumptionG} holds. For any finite $T_0>0$, for any $\mu\in \mathcal{P}_p(C_{T_0} C^r), u_0\in C^r$ and $Z\in C_{T_0}C^r$ the map $\Phi$ has a unique solution $u^\mu(Z,u_0)\in C_{T_0}C^r$.
\end{prop}

\ssk

\begin{Dem}
For $u,v\in C_{T_0}C^r$, using Assumption \eqref{EqAssumptionG} and \eqref{EqBoundIntegration}, we have
\[
\norme{\Phi(u)_t-\Phi(v)_t}_{C^r}  \leq \int_0^t \big\Vert P_{t-s}\big(g(u_s,\mu_s)\big) - P_{t-s}\big(g(v_s,\mu_s)\big)\big\Vert_{C^r} ds \leq \int_0^t  L \norme{u_s-v_s}_{C^r} \, ds.
\]
Denote by $\Delta_k(0,t)$ the simplex $\{0\leq s_1 \leq \dots\leq s_k\leq t\} $ and write $ds$ for $ds_1\dots ds_k$. An iteration of the previous bound gives
\[
\norme{\Phi^{\circ k}(u)_t - \Phi^{\circ k}(v)_t}_{C^r} \leq L^k \int_{\Delta_k(0,t)} \norme{u_{s_k}-v_{s_k}}_{C^r} ds \leq\frac{(LT_0)^k}{k!} \norme{u-v}_{C_{T_0}C^r}.
\]
The map $\Phi^{ \circ k}$ is thus contracting for $k$ large enough, so it has a unique fixed point.
\end{Dem}

\ssk

We now work with $(Z,u_0)$ random, an element of $\bbL^p\big(\Omega, C_{T_0}C^r \times C^r\big)$.

\ssk

\begin{prop}
For every $\mu\in\mathcal{P}_p(C_{T_0}C^r)$ the law of $u^\mu(Z,u_0)$ belongs to $\mcP_p(C_{T_0}C^r)$.
\end{prop}

\ssk

\begin{Dem}
Write $\delta_{\mathbf{0}}$ for Dirac distribution on the null function $\mathbf{0}$. We have the estimate
\begin{align*}
\norme{u_t^\mu}_{C^\alpha} &\leq C\Big( \norme{u_0}_{C^r} + \norme{Z_t}_{C^\alpha} + \int_0^t  \norme{ g(u_s^\mu,\mu_s)}_{C^r} ds \Big)  \\
&\leq C\left( \norme{u_0}_{C^r}+\norme{Z_t}_{C^r}+\int_0^t \norme{g(0,\delta_{\mathbf{0}})}+L\Big(\norme{u_s}_{C^r}+\mcW_{p,C^r}(\mu_s,\delta_{\mathbf{0}}) \Big)\dd s \right)  \\
&\leq  C\Big( \norme{u_0}_{C^r}+\norme{Z_t}_{C^r} + T_0 \norme{g(0,\delta_{\mathbf{0}})}_{C^\alpha} + L T_0 \mcW_{p,C_{T_0}C^r}(\mu,\delta_{\mathbf{0}}) \Big) + CL \int_0^t\norme{u_s}_{C^r} \dd s,
\end{align*}
for some positive constant $C$. We get the inequality
\[
\norme{u_t}_{C^r}\leq  C\Big( \norme{u_0}_{C^r}+\norme{Z_t}_{C^r} + T_0 \norme{g(0,\delta_{\mathbf{0}})}_{C^r} + T_0\mcW_{p,C_TC^r}(\mu,\delta_{\mathbf{0}}) \Big)e^{CLt}
\]
from Gronwall lemma, from which the conclusion follows since $Z\in \bbL^p(\Omega,C_{T_0}C^r)$.
\end{Dem}

\ssk

Set
\[
\Psi : \left\{ \begin{array}{ccc}  \mathcal{P}_p(C_{T_0}C^r) \times \bbL^p(\Omega, C_{T_0}C^r\times C^r) & \rightarrow & \mathcal{P}_p(C_{T_0} C^r) \\ \big(\mu, (Z,u_0)\big) & \mapsto & \mcL\big(u^\mu(Z,u_0)\big) \end{array} \right. 
\]
We define a \textbf{\textit{solution to \eqref{EqMcKV0Integral}}} as a fixed point of the map 
\[
\Psi\big(\cdot,(Z,u_0)\big) : \mathcal{P}_p(C_{T_0} C^r)\rightarrow \mathcal{P}_p(C_{T_0} C^r).
\]

\ssk

\begin{thm} \label{ThmAffine}
Suppose the Lipschitz regularity assumption \eqref{EqAssumptionG} on $g$ holds. Then \eqref{EqMcKV0Integral} has a unique solution denoted by $u(Z,u_0)$. We have the Lipschitz estimate
\begin{equation} \label{EqLispchitzEstimate}
\mathcal{W}_{p,C_{T_0}C^r}\big(\mcL(u(Z,u_0)) , \mcL(u(Z'  , u_0' ))\big) \lesssim \mathcal{W}_{p,C_{T_0}C^r\times C^r}\big(\mcL(Z,u_0) , \mcL(Z'  , u_0' )\big),
\end{equation}
for an implicit multiplicative constant that depends only on $g,p$ and $T_0$.
\end{thm}

\ssk

\begin{Dem}
Fix $(Z,u_0)$ and use the shorthand notation $\Psi_{Z,u_0}(\cdot)$ for $\Psi\big(\cdot,(Z,u_0)\big)$. For $\mu,\nu\in \mathcal{P}_p(C_TC^r)$ write $u^\mu$ and $u^{\nu}$ for $u^\mu(Z, u_0)$ and $u^{\nu}(Z, u_0)$, respectively. One has 
\[
u^\mu_t - u_t^{\nu} = \int_0^t \Big(P_{t-s}\big(g(u^\mu_s,\mu_s)\big) - P_{t-s}\big(g(u^{\nu}_s,\nu_s)\big)\Big)ds,
\]
and 
\[
\big\Vert u^\mu_t-u^{\nu}_t \big\Vert_{C^r}^p \leq C\int_0^t \Big(\big\Vert u^\mu_s - u^{\nu}_s\big\Vert_{C^r}^p + \mathcal{W}_p\big(\mu_{[0,s]},\nu_{[0,s]}\big)^p\Big) ds,
\]
for some positive constant $C$, so we get from Gronwall lemma the estimate
\begin{equation} \label{EqGronwall}
\mathcal{W}_{p,C_tC^r}\Big(\mcL(u^\mu_{[0,t]}),\mcL(u^{\nu}_{[0,t]})\Big)^p \leq Ce^{CT_0}\int_0^t \mathcal{W}_{p,C_sC^r}\big(\mu_{[0,s]},\nu_{[0,s]}\big)^pds.
\end{equation}
Recall we write $\{0\leq s_1\leq \cdots\leq s_k\leq T_0\}$ for the simplex and we use the shorthand notation $ds$ for $ds_1\dots ds_k$ when integrating over the simplex. Starting from \eqref{EqGronwall}, an elementary iteration gives for any $k\geq 2$
\begin{align*}
\mathcal{W}_{p,C_{T_0}C^r}\big(\Psi_{Z,u_0}^{\circ k}(\mu) , \Psi_{Z,u_0}^{\circ k}(\nu)\big)^p   &\leq  (Ce^{CT_0})^k\int_{\Delta^k_t}
\mathcal{W}_{p,C_{s_k}C^r}\big(\mu_{[0,s_k]} , \nu_{[0,s_k]}\big)^pds   \\
&\leq (Ce^{CT_0})^k \frac{1}{k!} \, \mathcal{W}_{p,C_{T_0}C^r}\big(\mu , \nu\big)^p.
\end{align*}
The map $\Psi_{Z,u_0}^{\circ k}$ is thus a contraction for $k$ sufficiently large, from which it follows that Equation \eqref{EqMcKV0} has a unique solution.

\ssk

Pick now $Z, Z' \in C_{T_0}C^r$ and $u_0 , u_0' \in C^r$. Pick $\mu\in \mathcal{P}_p(C_{T_0}C^r)$ and write $u$ and $u'$ for $u(Z, u_0)$ and $u(Z' , u_0')$, respectively. We can assume without loss of generality that $Z, Z' , u_0, u_0 $ are such that the $p$-th moment of $\norme{u - u'}_{C_{T_0}C^r}$ is equal to the $p$-Wasserstein distance between $\mcL(u(Z,u_0))$ and $\mcL(u(Z' , u_0' ))$. Since
\[
u_s - u'_s = P_s(u_0 - u'_0 ) + Z_s - Z' _s +  \int_0^s  \Big( P_{s-r}\big(g(u_r,\mu_r)\big) - P_{s-r}\big(g(u'_r,\mu_r)\big)\Big) dr,
\]
we have
\begin{equation*}
\sup_{s\in[0,t]} \norme{u_s - u'_s}_{C^r}  \leq  \norme{u_0 - u'_0 }_{C^r} + \norme{Z - Z' }_{C_T C^r} + C \int_0^t \norme{u_s - u'_s}_{C^\alpha} ds
\end{equation*}
and
\begin{align*}  
\bbE\Big[\sup_{s\in[0,t]} \norme{u_s - u'_s}_{C^r}^p\Big]  \lesssim_p  \norme{u_0 - u_0' }_{C^r}^p + \bbE\big[\norme{Z - Z' }_{C_T C^r}^p\big] + \int_0^t \bbE\Big[\sup_{r\in[0,s]}\norme{u_r - u'_r}_{C^r}^p\Big]ds.
\end{align*}
We get the Lipschitz estimate \eqref{EqLispchitzEstimate} from Gronwall lemma.
\end{Dem}

\ssk

Note that we do {\it not} assume that the noise $Z$ and the initial condition $u_0$ are independent.

\ssk

\subsection{Propagation of chaos$\boldmath{.}$   \hspace{0.15cm}}   
\label{SubsectionAdditiveChaos}

Let now $(Z^i,u_0^i)_{i\geq 1}$ be a sequence of independent, identically distributed, random variables with common distribution the law of $(Z,u_0)$. Write $\zeta^i$ for $\mathscr{L}(Z^i)$. Denote by $(\Omega,\mcF,\mathbb{P})$ the probability space on which this sequence of random variables is defined, with $\omega\in\Omega$ a generic element of $\Omega$. Fix $\omega\in\Omega$. For an integer $n\geq 1$ consider the interacting system of fields $\big(u^{1,n}(\omega),\dots,u^{n,n}(\omega)\big)$ with initial conditions $\big(u^1_0(\omega),\dots,u^n_0(\omega)\big)$ and dynamics
\begin{equation} \label{EqFieldSystemSimple} \begin{split}
(\partial_t-\Delta)u^{i,n}(\omega) &\;= \zeta^i(\omega) + g\big(u^{i,n}(\omega),\mu_t^n(\omega)\big),   \\
\mu_t^n(\omega) &\defeq \frac{1}{n}\sum_{k=1}^n \delta_{u_t^{k,n}(\omega)},
\end{split} \end{equation}
for $1\leq i\leq n$. H. Tanaka \cite{Tanaka} was the first to notice that the system \eqref{EqFieldSystemSimple} is actually, {\it for each $\omega\in\Omega$}, an equation of the form \eqref{EqMcKV0} set on the finite probability space $\{1,\dots,n\}$ equipped with the uniform probability measure $\lambda_n$. Following \cite{roughchaos}, we call this observation `{\sl Tanaka's trick}'. Random variables on the space $\{1,\dots,n\}$ are $n$-tuples indexed by $1\leq i\leq n$. Denote by $\mcL_{\lambda_n}(X)$ the law under $\lambda_n$ of an arbitrary random variable $X$ defined on $\{1,\dots,n\}$. Denote also by 
\[
U_n : j\mapsto j
\] 
the canonical random variable on $\{1,\dots,n\}$. Tanaka's trick says that a solution to the system
\[
(\partial_t-\Delta)u^i(\omega) = \zeta^i(\omega) + g\big(u^i(\omega), \mcL_{\lambda_n}(u^{U_n(\cdot)}(\omega))\big), \qquad (1\leq i\leq n)
\]
with parameter $\omega$ and chance element $i\in\{1,\dots,n\}$, is precisely given by the $n$-tuple 
\[
\big(u^{1,n}(\omega),\dots,u^{n,n}(\omega)\big)
\]
of solutions to the field system \eqref{EqFieldSystemSimple}.

\medskip

Recall  that a \textit{sequence $(\mu^n)_{n\geq 1}$ of probability measures on $E^n$}, invariant by the action on $E^n$ of the permutation group of $n$ elements, is said to be \textit{$\mu$-chaotic} if for every $1\leq k\leq n$ and $\phi_1,\dots\phi_k\in C_b(E)$, we have
\[
\mu^n\big(\phi_1 \otimes \dots \otimes \phi_k \otimes {\bf 1}^{\otimes (n-k)} \big) \underset{n\rightarrow\infty}{\longrightarrow} \prod_{i=1}^k \mu(\phi_i).
\]
A well-known criterion of $\mu$-chaoticity is given by the convergence in law of the empirical mean of an independent, identically distributed, $n$-sample of $\mu^n$ to the measure $\mu$ itself -- see for instance Proposition 2.2 in Sznitman's lecture notes \cite{Snitzman89}. Now the law of large numbers tells us that the empirical mean
\[
\frac{1}{n}\sum_{i=1}^n \delta_{(\zeta^i,u_0^i)(\omega)}
\]
converges $\mathbb{P}$-almost surely in $\mcW_{p,C_{T_0}C^{r-2}\times C^r}$ to $\mcL(\zeta,u_0)$; it converges a fortiori in law to the same limit. The following fact is thus a consequence of the Lipschitz estimate \eqref{EqLispchitzEstimate} and Sznitman's criterion. In the next statement we write $u\in \bbL^p(\Omega, C_{T_0}C^r)$ for the solution to Equation \eqref{EqMcKV0Integral}.
  
\medskip

\begin{thm}
For any integer $k\geq 1$ the law of the tuple $(u^{1,n},\dots, u^{k,n})$ converges weakly to $\mcL(u)^{\otimes k}$ when $n $ tends to $+\infty$.
\end{thm}

\medskip

\section{Basics on paracontrolled calculus}
\label{SectionBasics}

The study of Equation \eqref{EqMainEq} with a non-constant diffusivity $f(\cdot)$ requires that we use one of the languages that have been developed in the last ten years for the study of a large class of singular stochastic partial differential equations. The problem involved in this class of equations is best illustrated on the toy example of the parabolic Anderson model equation
\[
(\partial_t-\Delta)u = u\xi
\]
set on $\mathbf{\mathsf{T}}^2$, with $\xi$ a space white noise. Recall $\xi$ has almost surely H\"older space regularity $-1-\epsilon$ for all $\epsilon>0$. One expects from the Schauder estimates satisfied by the resolvent of the heat operator that $u$ has parabolic regularity $(\alpha-2)+2=\alpha$. This regularity is not sufficient for making sense of the product $u\xi$ since $\alpha + (\alpha-2)<0$. There are at least two languages one can use to circumvent this problem and set a robust solution theory for this equation and a whole class of equations involving the same pathology. We choose to work here with the language of paracontrolled calculus first introduced by Gubinelli, Imkeller \& Perkowski in \cite{GIP}. We recall in Section \ref{SubsectionBasics} the notions and results from paracontrolled calculus that we will use; we refer the reader to the works \cite{GPEnsaios, GubinelliPerkowskiRockner, GubinelliICM} of Gubinelli and Perkowski for some accounts of the basics on the subject. The results of the present section are sufficient to deal with the soft, but already interesting, case of a mean field Equation \eqref{EqMainEq} with diffusivity given by the model function \eqref{EqModelDiffusivity} with a $C^2$ kernel $k$. We deal with that case in Section \ref{SectionRegularInteraction} as a warm-up for Section \ref{SectionGeneral}.

On the methodological side, as in Section \ref{SectionAdditiveNoise}, we follow the usual two step fixed point formulation of a mean field equation: Fix a probability measure $\mu$ in an ad hoc space of probability measures, solve an associated equation where $\mu$ stands as a parameter and formulate the initial equation as a fixed point problem on the space of $\mu$. The additional complexity that stems from the singular feature of the dynamics \eqref{EqMainEq} will later be embodied in the structure of the space of $\mu$. In the setting of Section \ref{SectionRegularInteraction} it takes the classical form $\mu \in \mcP_p(\scrC^\alpha_T)$. {\color{black} It takes a different form in Section \ref{SectionGeneral}.}

\ssk

{\color{black} We recall some basic facts about paracontrolled calculus in Section \ref{SubsectionBasics}: The paraproduct and resonant operators $\varolessthan$ and $\odot$, the corrector $\sf C$, some paracontrolled structures on some sets of functions or distributions on ${\sf T}^2$ or $[0,T]\times{\sf T}^2$ and a paralinearization result. All these tools will be used in the sequel.}

\medskip

\subsection{Basics on paracontrolled calculus. \hspace{0.15cm}}
\label{SubsectionBasics}

Let $h_1$ and $h_2$ be two distributions on $\mathbf{\mathsf{T}}^2$. We use the notations $h_1\varolessthan h_2$ and $h_1\odot h_2$ for the paraproduct and the resonant operators defined from the Littlewood-Paley projectors as in \eqref{EqParaproduit} and \eqref{EqResonance}. For any $h_1\in C^{\alpha_1}$ and $h_2\in C^{\alpha_2}$ we have the continuity estimates
\[
\Vert h_1\varolessthan h_2 \Vert_{\alpha_2} \lesssim \Vert h_1\Vert_\infty \Vert h_2\Vert_{\alpha_2}   \qquad   \textrm{if $\alpha_1\geq 0$}
\]
and 
\[
\Vert h_1\varolessthan h_2 \Vert_{\alpha_1+\alpha_2} \lesssim \Vert h_1\Vert_{\alpha_1} \Vert h_2\Vert_{\alpha_2}   \qquad   \textrm{if $\alpha_1<0$}.
\]
The resonance $h_1\odot h_2$ is well defined if and only if $\alpha_1+\alpha_2>0$, in which case it is a continuous bilinear function of $h_1\in C^{\alpha_1}$ and $h_2\in C^{\alpha_2}$. (See for instance Section 2.6 in Bahouri, Chemin \& Danchin's textbook \cite{BCD} for a reference.) The following two continuity results involving the paraproduct and resonant operators will be used below. One has 
\begin{equation} \label{EqRContinuity}
\big\Vert a \varolessthan (b \varolessthan c) - (ab) \varolessthan c \big\Vert_{C^{\alpha_2+\alpha_3}} \lesssim \Vert a\Vert_{L^\infty} \Vert b\Vert_{C^{\alpha_2}} \Vert c\Vert_{C^{\alpha_3}},
\end{equation}
for all $a\in L^\infty, b\in C^{\alpha_2}$ with $\alpha_2$ in $(0,1)$ and $c\in C^{\alpha_3}$ with $-3<\alpha_3<3$. (The regularity exponent $3$ has no particular meaning; it is purely technical. One can find this statement and its proof in Proposition 14 of \cite{BailleulBernicotHighOrder}.) Last we recall from Lemma 2.4 of \cite{GIP} that the {\it corrector}
\[
{\sf C}(a,b,c) \defeq (a\varolessthan b)\odot c - a \, (b\odot c)
\]
has a continuous extension from $C^2\times C^2\times C^2$ to $C^{\alpha_1}\times C^{\alpha_2}\times C^{\alpha_3}$ with values in $C^{\alpha_1+\alpha_2+\alpha_3}$ if $\alpha_2+\alpha_3<0$ and $0<\alpha_1+\alpha_2+\alpha_3<1$.

\medskip

{\it \S1. Spatial paracontrolled structure.} From its definition $h_1 \varolessthan h_2$ is well defined for all distributions $h_1,h_2$ on $\mathbf{\mathsf{T}}^2$ and has high Fourier modes that are some modulations of the high Fourier modes of $h_2$ by the low Fourier modes of $h_1$. On that ground it makes sense to think of $h_1\varolessthan h_2$ as a distribution that `looks like' $h_2$. 

\bigskip

\noindent $\bullet$ {\it Definition --} This brings us to the following definition.

\ssk

\begin{defn} \label{DefnSpacePCStructure}
Pick a reference distribution $\Lambda\in C^\rho$ with $\rho\in\bbR$. A \textbf{distribution $v$ on} \emph{$\mathbf{\textsf{T}}^2$} is said to be \textbf{paracontrolled by $\Lambda$} if there exists a positive regularity exponent $\gamma$ and functions $v'\in C^\gamma$ and $v^\sharp\in C^{\gamma+\rho}$ such that
\[
v = (v'\varolessthan\Lambda) + v^\sharp.
\]
We denote by $\mcD^\gamma(\Lambda)$ the space of all such pairs $(v',v^\sharp)$, equipped with the norm 
\begin{equation*}
\norme{(v', v^\sharp)}_{\mcD^\gamma} \defeq \big\Vert v'\big\Vert_{C^\gamma} + \norme{v^\sharp}_{C^{\gamma+\rho}}.
\end{equation*}
For two reference distributions $\Lambda_1,\Lambda_2\in C^\rho$ and $\mathbf{v}_1=(v_1', v_1^\sharp)\in\mcD^\gamma(\Lambda_1)$ and $\mathbf{v}_2 = (v_2',v_2^\sharp)\in\mcD^\gamma(\Lambda_2)$, we set
\[
d_{\mcD^\gamma}(\mathbf{v}_1,\mathbf{v}_2) \defeq \big\Vert v_1' - v_2'\big\Vert_{C^\gamma} + \big\Vert v_1^\sharp - v_2^\sharp \big\Vert_{C^{\gamma+\rho}}.
\]
\end{defn}

\ssk

The expression `{\sl Gubinelli derivative of $v$}' is sometimes used to talk about $v'$. Note that the exponent $\gamma$ in $\mcD^\gamma(\Lambda)$ does \emph{not} refer to the regularity of $v$ but rather to the regularity exponents of $v'$. Indeed the distribution $v$ is in $C^\rho$. 

\bigskip

\noindent $\bullet$ {\it Paralinearization --} Let $a$ and $b$ be two functions in $C^\alpha$ with $\alpha\in(2/3,1)$ and with $a\in\mcD^\beta(b)$ for $\beta\in(2/3,\alpha)$, with Gubinelli derivative $a'$. {\it Bony's paralinearization} states that if $h$ stands for a $C^3_b$ function from $\bbR$ into itself then 
\[
h(a) = h'(a)\varolessthan a + R_h(a)
\] 
for some `remainder' term $R_h(a)\in C^{\alpha+\beta}$ with 
\[
\norme{R_h(a)}_{C^{\alpha+\beta}}\lesssim 1+\norme{a}_{C^\alpha}^2.
\]
The estimate \ref{EqRContinuity} implies that $h(a)\in\mcD^\beta(b)$, with Gubinelli derivative $h(a)' = h'(a)a'$ and $h'$ the classical derivative of $h$. (See e.g. Section 2.3 of \cite{GIP}.) We denote by $h(a)^\sharp$ the remainder term in the description of $h(a)$ as an element of $\mcD^\beta(b)$. 

The next proposition gives an example of application of this result that is of interest for us. Assume we are given a function $F\in C_b^3(\bbR^2,\bbR)$ and a $C^2_b$ kernel $k(x,y)$ on the torus ${\sf T}^2$. Recall $\beta\in(2/3,\alpha)$. For $a\in C^\alpha$ and $\mu\in\mcP_{p}(C^\alpha)$, with $1\leq p<\infty$, we define here for $x\in {\sf T}^2$
\begin{equation} \label{EqLongRangeF}
f(a,\mu)(x) = \int_{C^\alpha}\int_{{\sf T}^2} F\big(a(x), b(y)\big)k(x,y) \, dy \, \mu(db).
\end{equation}
This is a linear function of its measure argument. We denote below by $\partial_1F$ the partial derivative of $F$ with respect to its first argument.

\ssk

\begin{prop} \label{paralinearisation0}
{\color{black} For any} reference function $\Lambda\in C^\alpha$, any ${\bf a}=(a',a^\sharp)\in \mcD^{\beta}(\Lambda)$ and $\mu\in \mcP_p(C^\alpha)$, one has 
\[
f(a,\mu) = f(a, \mu)' \varolessthan \Lambda + f({\bf a}, \mu)^\sharp
\]
with 
\[
f(a,\mu)'(x) = \int_{C^\alpha}\int_{{\sf T}^2}\partial_1F\big(a(x), b(y)\big) k(x,y) \, dy \, \mu(db),
\]
and
\[
\norme{f({\bf a}, \mu)^\sharp}_{C^{\alpha+\beta}} \lesssim \big(1+\norme{\Lambda}^2_{C^\alpha}\big)   \Big(1+\norme{a'}_{C^\beta} + \norme{a^\sharp}_{C^\alpha} \Big) \Big(1+ \norme{a'}_{C^\beta} + \Vert a^\sharp\Vert_{C^{\alpha+\beta}} \Big).
\]
Furthermore, for $\Lambda_i\in C^\alpha$ and ${\bf a}_i\in \mcD^\beta(\Lambda_i), \mu_i\in\mcP_p(C^\alpha)$, for $i\in\{1,2\}$, one has
\begin{equation} \label{EqLipscthizEstimateFUNu1}
\big\Vert f({\bf a}_1,\mu_1)^\sharp - f({\bf a}_2,\mu_2)^\sharp\big\Vert_{C^{\alpha+\beta}} \lesssim d_{\mcD^\beta}\big({\bf a}_1, {\bf a}_2\big) + \mcW_{p,C^\alpha}\big(\mu_1, \mu_2\big) + \Vert \Lambda_1 - \Lambda_2\Vert_{C^\alpha},
\end{equation}
for an implicit constant that is a polynomial of degree 3 of 
\[
\max_{i=1,2}\Big\{ 1, \Vert {\bf a}_i\Vert_{\mcD^\beta(X^i)}, \mcW_{p,C^\alpha}(\mu_i,\delta_0), \Vert X_t^i\Vert_{C^\alpha}  \Big\}.
\]
\end{prop}

\ssk

\begin{Dem}
We paralinearize with respect to the $x$ variable, with $y$ in the role of a parameter in the paraproducts below. We use the shorthand notations 
\[
k_y(x) \defeq k(x, y), \quad F_{b(y)}(\cdot) \defeq F(\cdot,b(y)).
\]
With these notations one has
\begin{equation*} \begin{split}
F(a, b(y)) &= \partial_1F(a, b(y)) \varolessthan  a  + F_{b(y)}(a)^\sharp   \\
	      &= \big(\partial_1F(a, b(y)) a'\big) \varolessthan \Lambda + \partial_1F(a, b(y)) \varolessthan a^\sharp   \\
  	      &\quad+ \partial_1F\big(a, b(y)) \varolessthan (a' \varolessthan \Lambda) - \big( \partial_1F\big(a, b(y)\big) a' \varolessthan \Lambda \big) + F_{b(y)}(a)^\sharp
\end{split} \end{equation*}
and
\begin{align*}
    f(a, \mu) &= \bigg( a' \int_{{\sf T}^2 \times C^\alpha} \partial_1F\big(a, b(y)\big) k_y \, dy \,\mu(db) \bigg) \varolessthan \Lambda
    \\
    &\quad+ \int_{{\sf T}^2 \times C^\alpha} \bigg(\big( \partial_1 F\big(a, b(y)\big) a' \varolessthan \Lambda\big) k_y - \big\{ k_y\partial_1 F\big(a, b(y)\big) a'\big\} \varolessthan \Lambda \bigg) dy\,\mu(db)   \\
    &\quad+ \int_{{\sf T}^2 \times C^\alpha} k_y \Big\{ \partial_1F(a, b(y)) \varolessthan (a' \varolessthan \Lambda) - \big\{\partial_1F(a, b(y)) a'\big\} \varolessthan \Lambda \Big\} \, dy \,\mu(db)
    \\
    &\quad+ \int_{{\sf T}^2 \times C^\alpha} k_y F_{b(y)}(a)^\sharp \, dy \,  \mu(db) + \int_{C^\alpha}\int_{{\sf T}^2} \big(\partial_1 F(a, b) \varolessthan a^\sharp \big) k_y \, dy \,\mu(db)
    \\
    &\eqdef \bigg( a' \int_{{\sf T}^2 \times C^\alpha} \partial_1F\big(a, b(y)\big) k_y \, dy \, \mu(db)\bigg) \varolessthan \Lambda + f(a,b)^\sharp.
\end{align*}
We estimate each term separately to show that the remainder is regular, using some commutator type estimates when needed. First, since $k_y$ is $C^2_b$ and $\alpha+\beta<2$ we have from \eqref{EqRContinuity} the continuity estimate

\begin{equation*} \begin{split}
\Big\Vert \Big( \big\{\partial_1 F\big(a, b(y)\big) a' \big\} \varolessthan \Lambda \Big) k_y - \big\{ k_{z'}\partial_1 &F\big(a, b(y)\big) a'\big\} \varolessthan \Lambda \Big\Vert_{C^{\alpha+\beta}} 
\\
&\lesssim  \norme{k_y}_{C^{2\alpha}}\norme{\partial_1 F\big(a, b(y)\big) a'}_{C^\beta} \norme{\Lambda}_{C^\alpha}
\\
&\lesssim \norme{k}_{C^2_b}\big(1+\norme{a}_{C^\alpha}\big)\norme{a'}_{C^\beta}\norme{\Lambda}_{C^\alpha}
\\
&\lesssim \big(1+\norme{\Lambda}^2_{C^\alpha}\big) \Big(1+\norme{a'}^2_{C^\beta} + \Vert a^\sharp\Vert^2_{C^\alpha} \Big)
\end{split} \end{equation*}
and
\begin{equation*} \begin{split}
\Big\Vert \partial_1F(a, b(y)) \varolessthan (a' \varolessthan \Lambda) - \big\{ \partial_1F(a, b(y)) a'\big\} \varolessthan \Lambda \Big\Vert_{C^{\alpha+\beta}} &\lesssim \norme{\partial_1 F(a, b(y))}_{C^\beta} \norme{a'}_{C^\beta}\norme{\Lambda}_{C^\alpha}  \\
&\lesssim \big(1+\norme{a}_{C^\alpha}\big) \norme{a'}_{C^\beta}\norme{\Lambda}_{C^\alpha}   \\
&\lesssim \big(1+\norme{\Lambda}^2_{C^\alpha}\big) \Big(1+\norme{a'}^2_{C^\beta} + \Vert a^\sharp\Vert^2_{C^\alpha} \Big)
\end{split} \end{equation*}
and
\begin{equation*} \begin{split}
 \norme{k_y F_{b(y)}(a)^\sharp}_{C^{\alpha+\beta}} \lesssim \norme{F_{b(y)}}_{C^3_b}\big(1+\norme{a}_{C^\alpha}^2\big) \lesssim \big(1+\norme{\Lambda}^2_{C^\alpha}\big)\Big(1+\norme{a'}^2_{C^\beta} + \Vert a^\sharp\Vert_{C^\alpha}^2\Big)
\end{split} \end{equation*}
and 
\begin{equation*} \begin{split}
\norme{ \big(\partial_1 F(a, b(y)) \varolessthan a^\sharp \big) k_y}_{C^{\alpha+\beta}} &\lesssim \big(1+\norme{a}_{C^\alpha}\big) \Vert a^\sharp\Vert_{C^{\alpha+\beta}}   \\
&\lesssim \big(1+\norme{\Lambda}_{C^\alpha}\big)\Big(1+\norme{a'}_{C^\beta} + \Vert a^\sharp\Vert_{C^\alpha}\Big) \Vert a^\sharp\Vert_{C^{\alpha+\beta}}.
\end{split} \end{equation*}
Integrating over $y$ and summing we see that
\[
\norme{f({\bf a},\mu)^\sharp}_{C^{\alpha+\beta}} \lesssim \big(1+\norme{\Lambda}_{C^\alpha}^2 \big)\Big(1+\norme{a'}_{C^\beta} + \Vert a^\sharp\Vert_{C^\alpha} \Big)\Big(1+ \norme{a'}_{C^\beta} + \Vert a^\sharp\Vert_{C^{\alpha+\beta}}  \Big). 
\]
We leave the proof of the Lipschitz estimate \eqref{EqLipscthizEstimateFUNu1} to the reader as it is very similar to what we have just done.
\end{Dem}

\medskip

{\it \S2. Parabolic paracontrolled structure.} We will denote by $k_1\prec k_2$ the `parabolic' paraproduct on {\it spacetime distributions} introduced in Section $5$ of \cite{GIP}. It is defined for $k_1,k_2\in C_T\mcD'({\sf T}^2)$ as 
\[
k_1\prec k_2 \defeq \sum_{i\geq -1} \Delta_{<i-1}(Q_ik_1)\Delta_i(k_2), 
\]
with 
\[Q_ik_1(t) \defeq \int_0^T 2^{2i}\varphi(2^{2i}(t-s))k_1(s)\dd s
\] 
and $\varphi$ some smooth compactly supported function with mass $1$. It is a parabolic version of the space paraproduct operator $\varolessthan$ that has the same analytic properties in the scale of Besov parabolic function spaces as the operator $\varolessthan$ in the scale of spatial Besov function spaces. The parabolic paraproduct almost commutes with the heat operator in the sense that we have the useful estimate
\[
\big\Vert (\partial_t-\Delta)(k_1\prec k_2) - k_1 \prec \big((\partial_t-\Delta)k_2\big)\big\Vert_{C_TC^{\alpha-2+\beta}} \lesssim   \norme{k_1}_{\scrC^\alpha_T} \norme{k_2}_{C_TC^\beta}.
\] 
When applied to some parabolic distributions $k_1\in C_T^{\alpha/2}L^\infty, k_2\in C_TC^\beta$ the two paraproducts $\varolessthan$ and $\prec$ are related by the continuity relation
\begin{equation} \label{continuityrelation}
\big\Vert k_1\varolessthan k_2 - k_1 \prec k_2\big\Vert_{C_TC^{\alpha+\beta}} \lesssim \norme{k_1}_{C_T^{\alpha/2}L^\infty} \norme{k_2}_{C_TC^\beta}.
\end{equation} 
These above two estimates are the content of Lemma 5.1 of \cite{GIP}. We use the $\prec$ paraproduct and a slightly different notion of size to deal with parabolic functions paracontrolled by a reference parabolic function $\Xi$.

\ssk

\begin{defn} \label{DefnParabolicPCStructure}
Pick a reference function $\Xi\in\scrC_T^\rho$ with $\rho>0$. A \textbf{parabolic function $u$ on} \emph{$[0,T]\times\mathbf{\textsf{T}}^2$} is said to be \textbf{paracontrolled by $\Xi$} if there exists a positive regularity exponent $\gamma$ and a function $u'\in\scrC^\gamma_T$ such that
\[
u^\sharp \defeq u-u'\prec \Xi \in \scrC_T^\rho
\]
and 
\[
\sup_{t\in(0,T]}t^{\gamma/2}\big\Vert u^\sharp_t\big\Vert_{C^{\rho+\gamma}} < +\infty.
\]
We denote by {\color{black} $\mcE_T^{\gamma}(\Xi)$} the space of all such pairs $(u',u^\sharp)$; it is equipped with the norm
\[
\big\Vert(u',u^\sharp)\big\Vert_{\mcE^{\gamma}_T} \defeq \big\Vert u'\big\Vert_{\scrC_T^\gamma} + \big\Vert u^\sharp\big\Vert_{\scrC_T^\rho} + \sup_{t\in(0,T]}t^{\gamma/2}\big\Vert u^\sharp_t\big\Vert_{C^{\rho+\gamma}}.
\]
For two reference functions $\Xi_1, \Xi_2\in\scrC^\rho_T$ and $\mathbf{u}_1=(u_1',u_1^\sharp)\in\mcE^{\gamma}_T(\Xi_1)$ and $\mathbf{u}_2 =(u_2',u_2^\sharp)\in\mcE^{\gamma}_T(\Xi_2)$ we set
\[
d_{\mcE^{\gamma}_T}({\bf u}_1,\mathbf{u}_2) \defeq \big\Vert u_1'-u_2'\big\Vert_{\scrC^\gamma} + \big\Vert u_1^\sharp-u_2^\sharp \big\Vert_{\scrC^\rho_T} + \sup_{t\in(0,T]}t^{\gamma/2}\big\Vert u_1^\sharp-u_2^\sharp\big\Vert_{C^{\rho+\gamma}}.
\]
\end{defn}

\ssk

We note that if ${\bf u}\in\mcE^\gamma_T(\Xi)$, with $\Xi\in\scrC^\rho_T$ and $\rho>0$, then it follows from \eqref{continuityrelation} that ${\bf u}_t=(u'_t,u^\sharp_t) \in \mcD^\gamma(\Xi_t)$ for all $0<t\leq T$.

\medskip

\noindent \textit{\textbf{Warning}} -- {\it We will only work in the sequel with some situations where the reference object $\Xi\in\scrC_T^\rho$ in a parabolic paracontrolled structure is a function and $\rho>0$. On the other hand we will sometimes deal with some situations where the reference object $\Lambda\in C^\rho$ in a space paracontrolled structure is only a distribution and $\rho<0$.   }

\medskip

We now recall a kind of Schauder-type estimate first proved by Gubinelli, Imkeller \& Perkowski in Section 5 of \cite{GIP}. The starting point of the next statement is the description for each time $0<t\leq T$ of the right hand side of a parabolic equation as a space paracontrolled distribution, in the sense of Definition \ref{DefnSpacePCStructure}. The statement provides as an outcome a description of the solution of the equation as a parabolically paracontrolled function, in the sense of Definition \ref{DefnParabolicPCStructure}. 

\medskip

\begin{prop} \label{schauder3} 
For $\pi\in C_TC^{\alpha-2}$ we set $\Pi=\scrL^{-1}(\pi)\in\scrC^\alpha_T$. Then for every {\color{black} $w'\in\scrC^\alpha_T$ and $w^\sharp$} such that 
\begin{equation} \label{EqConditionRemainderSize}
\sup_{t\in(0,T]} t^{\beta/2}\big\Vert w_t^\sharp\big\Vert_{(\alpha-2)+\beta} < \infty,
\end{equation}
for every $u_0\in C^\alpha$, the solution $u$ to the equation
\begin{equation} \label{EqSchauderDynamics}
(\partial_t-\Delta)u = w' \varolessthan \pi + w^\sharp, \quad u(0)=u_0,
\end{equation}
belongs to {\color{black} $\mcE^\beta_T(\Pi)$} and $u'=w'$. We further have the estimate 
\[
\norme{(u',u^\sharp)}_{\mcE^\beta_T(\Pi)} \lesssim \norme{u_0}_{C^\alpha} + T^{(\alpha-\beta)/2}\Big( \norme{w'}_{\scrC^\alpha_T}\big(1+\norme{\pi}_{C_TC^{\alpha-2}}\big) + \sup_{t\in(0,T]} t^{\beta/2}\big\Vert w^\sharp_t\big\Vert_{C^{(\alpha-2)+\beta}} \Big).
\]
For $i\in\{1,2\}$ and $w'_i \in \scrC^\alpha_T$ and $w^\sharp_i$ satisfying \eqref{EqConditionRemainderSize}, for different distributions $\pi_i\in C_TC^{\alpha-2}$, set 
\[
m' \defeq \max_{i\in\{1,2\}} \Big\{1,\Vert w'_{i}\Vert_{\scrC^\alpha_T}, \norme{\pi_i}_{C_TC^{\alpha-2}}\Big\}
\]
and denote by $u_1,u_2$ the corresponding solutions to Equation \eqref{EqSchauderDynamics} with associated paracontrolled decomposition ${\bf u}_1, {\bf u}_2$. There is a quadratic polynomial $P$ such that we have
\begin{equation*} \begin{split}
\text{d}_{\mcE^\beta_T}({\bf u}_1, {\bf u}_2) \lesssim P(m') \, T^{(\alpha-\beta)/2}\Big(&\norme{w'_1 - w'_2}_{\scrC^\alpha_T} + \Vert \pi_1-\pi_2\Vert_{C_TC^{\alpha-2}}   \\
&+ \sup_{t\in (0,T]}t^{\beta/2}\big\Vert w_1^\sharp(t) - w_2^\sharp(t)\big\Vert_{C^{(\alpha-2)+\beta}} \Big).
\end{split} \end{equation*}
\end{prop}

\vfill \pagebreak

\subsection{Noise enhancement and product definition. \hspace{0.15cm} }
\label{SectionNoiseEnhancement}

Fix a positive finite time horizon $T$. 

\ssk

{\it \S1. Noise enhancement.} We define the lifting operator
\begin{equation*} \begin{split} 
{\sf L} : C_TC^\infty\times C([0,T],\bbR) &\longrightarrow C_TC^\infty\times C_TC^\infty   \\
		(\ell, c)		&\longrightarrow \big(\ell \, , \, \scrL^{-1}(\ell)\odot\ell - c\big).
\end{split} \end{equation*}
The \textbf{\textit{space ${{\boldsymbol{\mcN}}_T}$ of enhanced noises}} is the closure in $C_TC^{\alpha-2}\times C_TC^{2\alpha-2}$ of the range of $\sf L$. The letter ${{\boldsymbol{\mcN}}_T}$ is chosen for the word `{\sl noise}'. We denote by 
\[
\widehat\zeta = \big(\zeta,\zeta^{(2)}\big)
\]
a generic element of ${{\boldsymbol{\mcN}}_T}$. The natural norm of $\widehat\zeta$ as an element of the ambiant product space is denoted by $\Vert\widehat\zeta\,\Vert$. We use the Schauder estimate \eqref{EqBoundIntegration} to define an element of $\mathscr{C}^\alpha_T$ by setting here
\[
Z  \defeq \scrL^{-1}(\zeta) \in \mathscr{C}^\alpha_T.
\]
The following statement provides a large class of random noises with a natural enhancement defined as some random element of ${{\boldsymbol{\mcN}}_T}$ from a renormalization procedure. The following renormalization statement can be seen as a consequence of a similar statement in Gubinelli \& Perkowski's work \cite{KPZreloaded} or as a direct consequence of the deep results on renormalized models obtained in the recent years. For the reader's convenience we provide a self-contained pedestrian proof of it in Appendix {\sf \ref{SectionEnhancingNoises}}.

\ssk

\begin{thm} \label{ThmNoiseEnhancement}
Let $(\xi_t)_{0\leq t\leq T_0}$ stand for a time-dependent Gaussian random distribution on ${\sf T}^2$ with covariance of the form
\[
\bbE\big[(\xi_t,\phi)(\xi_s,\psi)\big] = c(t,s) \langle \psi\star D , \phi\rangle_{L^2}
\]
for some distribution $D$ on ${\sf T}^2$, for any test functions $\psi,\phi\in C^\infty$, with $\star$ the convolution operator. We assume that the Fourier transform of $D$ satisfies for some $\eta< 1-\alpha$ the condition 
\[
|\widehat{D}(k)| \lesssim |k|^\eta,
\]
and that the function $c$ satisfies the inequality
\[
0\leq c(t,t) + c(s,s) - 2c(s,t) \leq \vert t-s\vert^\delta
\]
for some positive exponent $\delta$. Then one defines a random variable $X\odot \xi\in \bbL^1(\Omega, C_TC^{2\alpha-2})$ setting
\begin{equation} \label{EqDefnXOdotXi}
\big(X\odot\xi)(t) \defeq \int_0^t \Big( P_{t-s}(\xi_s)\odot \xi_t - \bbE[P_{t-s}(\xi_s)\odot \xi_t]\Big)ds 
\end{equation}
One further has $X\odot \xi\in \bbL^p(\Omega, C_TC^{2\alpha-2})$ for all $1\leq r<\infty$, and if $\xi^\epsilon$ stands for a space regularization of $\xi$ by convolution then 
\[
{\sf L}\big(\xi^\epsilon , \bbE[X^\epsilon \odot \xi^\epsilon]\big)
\]
converges in $\bbL^r(\Omega, C_TC^{2\alpha-2})$ to $X\odot \xi$ as $\epsilon>0$ goes to $0$.
\end{thm}

\ssk

As a shorthand notation, for $c\in C([0,T],\bbR)$, we set here for later use
\begin{equation} \label{EqDefnLc}
{\sf L}_c(\cdot) \defeq {\sf L}(\cdot,c).
\end{equation}

\ssk

{\it \S2. Product definition.} The end of this section deals with deterministic enhanced noises. The datum of an element of ${{\boldsymbol{\mcN}}_T}$ allows to give a definition of some otherwise ill-defined product. Recall that $2/3<\beta<\alpha<1$ so we have $\beta>2-2\alpha$.

\ssk

\begin{defn}  \label{produit}
Pick $\widehat\zeta\in{{\boldsymbol{\mcN}}_T}$ and $0<t\leq T$. Pick ${\bf u}_t = (u'_t,u^\sharp_t) \in \mcD^\beta(Z_t)$. We \emph{define} the product ${\bf u}_t\zeta_t$ as the element of $\mcD^\beta(\zeta_t)$ defined by the decomposition
\begin{equation*} \begin{split} 
{\bf u}_t\zeta_t &\defeq u_t \varolessthan \zeta_t + ({\bf u}_t\zeta_t)^\sharp,
\end{split} \end{equation*}
where
\[
({\bf u}_t\zeta_t)^\sharp \defeq \zeta_t \varolessthan u_t + u_t^\sharp\odot \zeta_t + {\sf C}\big(u_t', Z_t, \zeta_t\big) + u_t'\zeta^{(2)}_t
\]
and
\begin{equation} \label{EqEstimateProductSharp} \begin{split}
\norme{({\bf u}_t\zeta_t)^\sharp}_{C^{\alpha-2+\beta}} &\lesssim \norme{{\bf u}}_{\mathcal{D}^\beta(Z_t)}\Big( \norme{\zeta_t}_{C^{\alpha-2}} + \norme{Z_t}_{C^\alpha} \norme{\zeta_t}_{C^{\alpha-2}} + \big\Vert\zeta^{(2)}_t\big\Vert_{C^{2\alpha-2}}\Big)   \\
&{\color{black} \lesssim \norme{{\bf u}}_{\mathcal{D}^\beta(Z_t)} \big(1+\Vert \widehat{\zeta}\Vert_{{{\boldsymbol{\mcN}}_T}}^2\big)}.
\end{split} \end{equation}
\end{defn}

\ssk

For $\widehat\zeta^i = \big(\zeta^i,\zeta^{i(2)}\big)\in{{\boldsymbol{\mcN}}_T}, \,Z^i = \scrL^{-1}(\zeta^i)$ and ${\bf u}^i_t\in\mcD^\beta(Z^i_t)$, with $i\in\{1,2\}$, we set 
\begin{equation*} \begin{split}
{\color{black} m_1} &\defeq \underset{i\in\{1,2\}}{\max}\Big\{ {\big\Vert} {\zeta^i}\big\Vert_{C^{\alpha-2}}, {\big\Vert} {\zeta^{i(2)}}\big\Vert_{C^{2\alpha-2}}, \norme{{\bf u}^i_t}_{\mcE^\beta_T(Z^i_t)} \Big\}   \\
{\color{black} m_2} &{\color{black} \defeq \underset{i\in\{1,2\}}{\max}\Big\{ {\big\Vert} {\zeta^i}\big\Vert_{C^{\alpha-2}}, {\big\Vert} {\zeta^{i(2)}}\big\Vert_{C^{2\alpha-2}}\Big\}   }.
\end{split} \end{equation*}
The proof of the following proposition can be found in \cite{GIP}, Theorem 3.7 therein.

\ssk

\begin{prop} \label{PropLocalLipschitzSharp}
We have the local Lipschitz estimate
\begin{align*}
\Big\Vert({\bf u}^1_t\zeta^1_t)^\sharp - \big({\bf u}^2_t \zeta^2_t\big)^\sharp\Big\Vert_{C^{\alpha-2+\beta}} \lesssim {\color{black} m_1\Big(d_{\mcD^\beta}\big({\bf u}^1_t, {\bf u}^2_t\big) + \big\Vert \widehat\zeta\hspace{1pt}^1 - \widehat\zeta\hspace{1pt}^2\,\big\Vert_{{\boldsymbol{\mcN}}_T} \Big)}
\end{align*}
{\color{black} and if $\widehat\zeta^1=\widehat\zeta^2=\widehat\zeta$ and ${\bf u}^i_t\in\mcD^\beta(Z_t)$ we have}
\begin{align*}
{\color{black} \Big\Vert({\bf u}^1_t\zeta_t)^\sharp - \big({\bf u}^2_t \zeta_t\big)^\sharp\Big\Vert_{C^{\alpha-2+\beta}} \lesssim m_2 \, d_{\mcD^\beta}\big({\bf u}^1_t, {\bf u}^2_t\big)   }.
\end{align*}
Further if {\color{black} ${\bf u}\in\mcE^\beta_T(Z)$} then the function $t\mapsto {\bf u}_t\zeta_t$ is in $C_TC^{\alpha-2}$.
\end{prop}


\section{Mean field equations with a regular interaction in the diffusivity}
\label{SectionRegularInteraction}

We show in this section how to use directly the settings and statements of Section \ref{SectionBasics} to formulate and solve uniquely the mean field type equation \eqref{EqMainEq} in a particular, but still interesting, case. We assume in this section that the diffusivity $f$ in \eqref{EqMainEq} has the form
\begin{equation} \label{EqFormulaFNiceKernel}
f(a,\mu)(x) = \int_{C^\alpha}\int_{{\sf T}^2} F\big(a(x), b(y)\big)k(x,y) \, dy \, \mu(db)
\end{equation}
with $F\in C^{{\color{black} 2}}_b(\bbR^2,\bbR)$ and $k\in C^2_b({\sf T}^2\times{\sf T}^2,\bbR)$. Proposition \ref{paralinearisation0} holds for this function, saying that $f(u,\mu)$ is paracontrolled in space if $u$ is paracontrolled in space. Given an enhanced noise $\widehat{\zeta}\in{{\boldsymbol{\mcN}}_{T_0}}$, {\color{black} for some time horizon $0<T_0<\infty$}, assume we are given ${\bf u}_t\in\mcD^\beta(Z_t)$ and a probability measure $\mu_t\in\mcP_p(C^\alpha)$ for every $0<t\leq T_0$. We denote by $f({\bf u}_t,\mu_t)\in \mcD^\beta(Z_t)$ the space-paracontrolled function given by Proposition \ref{paralinearisation0}. The product $f({\bf u}_t,\mu_t)\zeta_t$ is then defined for each $t$ as in Definition \ref{produit}. {\color{black} (We will proceed differently in Section \ref{SectionGeneral}, where $k(x,y)=\delta_x(y)$, as this situation is outside of the scope of Definition \ref{produit}.)}

\ssk

{\color{black} For the moment} we assume throughout this section that the drift term {\color{black} $g(a,\mu)$ in \eqref{EqMainEq} is of the form
\[
g(a,\mu)(x) = {\sf g}(a(x),\mu)
\]
for some function ${\sf g} : \bbR\times\mcP_p(C^\alpha)\rightarrow\bbR$} that is a $C^1$ function of $z$ and satisfies 
\[
\sup_{z\in\bbR}\sup_{\mu\in \mcP_p(C^\alpha)} |\partial_z{\sf g}(z,\mu)| < \infty
\]
and the Lipschitz condition
\begin{equation} \label{EqConditionLipG}   
{\color{black} \norme{{\sf g}(\cdot,\mu) - {\sf g}(\cdot,\nu)}_{C^1} \lesssim \mcW_{p,C^\alpha}\big(\mu, \nu\big).   }
\end{equation}

\ssk

{\color{black} This condition is stronger than the Lipschitz condition \eqref{EqAssumptionG} from Section \ref{SubsectionAdditive}. The reason for this difference is that $g$ is here considered as a remainder term, hence more regular, while in Section \ref{SubsectionAdditive} we only need $g$ to produce an element of $C^\alpha$ after applying the classical Schauder estimate, so $g$ is controlled in $C^{\alpha-2}$ in that section. The minimal regularity required from $g$ is $C^{(\alpha-2)+\beta}$. We ask the relatively stronger bound \eqref{EqConditionLipG} in $C^1$ norm as this kind of regularity is needed to use some non-explosion result on the solutions of some (gPAM) type equations.

\ssk

We follow the classical formulation of a mean field type equation $X=\Phi(X,\mcL(X))$, for some random variable $X$ and some function $\Phi$. We first fix a probability $\mu$ on the space where $X$ takes its values and show that the fixed point equation $X=\Phi(X,\mu)$ has a unique solution $X^\mu$. A solution to the mean field equation is then defined as a fixed point of the map $\mu\mapsto\mcL(X^\mu)$ in some appropriate space of probability measures. In the present case the first step of the reasoning  is dealt with by the following proposition. Recall from Definition \ref{DefnParabolicPCStructure} the definition of the spaces $\mcE$.   }

\ssk

\begin{prop} \label{SolPAM1}
Fix $0<T_0<\infty$. For every initial condition $u_0\in C^\alpha$, for every enhanced noise $\widehat \zeta\in{{\boldsymbol{\mcN}}_{T_0}}$ and any $\mu\in \mcP_{p}(\scrC^\alpha_{T_0})$, there exists a unique solution in {\color{black} $\mcE^\beta_{T_0}(Z)$} to the equation 
\begin{equation} \label{EqEquationMuFixed}
\scrL u_t = f({\bf u}_t,\mu_t)\zeta_t + g(u_t,\mu_t).
\end{equation}
\end{prop}

\ssk

\begin{Dem}
{\color{black} We separate the analysis of the local in time well-posedness from the non-explosion question.}   \vspace{0.15cm}

 {\it Step 1. Local in time well-posedness.} Rewrite \eqref{EqEquationMuFixed} as the fixed point equation
\[
u_t = P_t (u_0) + \int_0^t P_{t-s}\big(f({\bf u}_s,\mu_s)\zeta_s + g(u_s, \mu_s)\big)\,ds.
\]
We get from Proposition \ref{paralinearisation0} and Proposition \ref{PropLocalLipschitzSharp} that $f({\bf u}_s,\mu_s)\zeta_s + g(u_s, \mu_s)$ is for each $s$ an element of $\mcD^\alpha(\zeta_s)$ with Gubinelli derivative $f(u_s,\mu_s)$ and remainder $(f({\bf u}_s,\mu_s)\zeta_s)^\sharp + g(u_s,\mu_s)$. With the Schauder-type Proposition \ref{schauder3} in mind, we check that $f(u,\mu)\in\mathscr{C}^\alpha_{T_0}$ and $(f({\bf u}_s,\mu_s)\zeta_s)^\sharp + g(u_s,\mu_s)$ satisfies \eqref{EqConditionRemainderSize}. 

Take ${\bf u}\in\mcE^\beta_{T_0}(Z)$. First, for $(s,x),(t,y)\in [0,T_0]\times {\sf T}^2$ with $s\leq t$, one has
\begin{align*}
\big|f(u_t,\mu_t)(y) - f(u_s &,\mu_s)(x)\big|   \\
&= \Big|\int_{{\sf T}^2\times\scrC_{T_0}^\alpha} F\big(u_t(y),v_t(z)\big)k(y,z) - F\big(u_s(x),v_s(z)\big) k(x,z) \, dz\,\mu(dv)\Big|   \\
&\leq \int_{{\sf T}^2\times\scrC_{T_0}^\alpha} \bigg( \big|F\big(u_t(y),v_t(z)\big) \big(k(y,z)-k(x,z)\big)\big|
\\
&\hspace{3cm}+ \big|F\big(u_t(y),v_t(y)\big)-F\big(u_s(x),v_s(x)\big)\big| \, |k(x,z)| \bigg) \, dz\,\mu(dv)
\\
&\lesssim \int_{{\sf T}^2\times\scrC_{T_0}^\alpha} \Big( |x-y| +\big(\norme{u}_{\scrC^\alpha_t} + \norme{v}_{\scrC^\alpha_t} \big) \big(|x-y|^\alpha+|t-s|^{\alpha/2}\big)\Big) \, dz\,\mu(dv)
\\
&\lesssim \big(1+\norme{u}_{\scrC^\alpha_T}+\mcW_{p,\scrC^\alpha_T}(\mu,\delta_{\bf 0})  \big)\big(|x-y|^\alpha+|t-s|^{\alpha/2}\big),
\end{align*}
so we have the norm estimate
\[
\norme{f(u,\mu)}_{\scrC^\alpha_{T_0}}\lesssim \big(1 + \norme{Z}_{\scrC_{T_0}^\alpha}\big) \Big(1+ \Vert {\bf u}\Vert_{\mcE^\beta_{T_0}} + \mcW_{p,\scrC^\alpha_{T_0}}(\mu,\delta_{\bf 0})\Big).
\]
Second, for $0<T\leq T_0$, we can see $\mu$ as an element of $\mcP_p(\scrC^\alpha_T)$ and rewrite the above estimate on $\big|f(u_t,\mu_t)(y) - f(u_s , \mu_s)(x)\big|$, for $0\leq s\leq t\leq T$, with $\scrC^\alpha_T$ in place of $\scrC^\alpha_{T_0}$  in the integrals. We thus have
\begin{equation}\label{estimeerhs1}
\sup_{t\in(0,T]} t^{\beta/2}\big\Vert \big(f({\bf u}_t,\mu_t)\lambda_t\big)^\sharp + g(u_t,\mu_t)\big\Vert_{\alpha-2+\beta} \lesssim \big( 1 + \Vert\widehat\zeta\hspace{2pt}\Vert_{{\boldsymbol{\mcN}}_T}^{{\color{red} 4}} \big) \Big(1 + \Vert {\bf u}\Vert^2_{\mcE^\beta_T} + \mcW_{p,\scrC^\alpha_T}\big(\mu,\delta_{\bf 0}\big)\Big).
\end{equation}
from Proposition \ref{paralinearisation0} and Proposition \ref{PropLocalLipschitzSharp} {\color{black} and the Lipschitz assumption \eqref{EqConditionLipG} on $g$}. It follows from Proposition \ref{schauder3} that the map 
\[
\Phi_{\widehat\zeta,\mu} : \mcE^\beta_T(Z) \rightarrow\mcE^\beta_T(Z)
\]
which associates to ${\bf u}\in \mcE^\beta_T(Z)$ the solution {\color{black} paracontrolled} $w$ of the equation 
\[
\scrL w = f({\bf u},\mu) \zeta + g(u,\mu),
\] 
with initial condition $w_0=u_0$, is well defined and satisfies the estimate
\begin{align}\label{estimee111}
\big\Vert\Phi_{\widehat\zeta, \mu}({\bf u})\big\Vert_{\mcE^\beta_T(Z)} &\lesssim  \norme{u_0}_{C^\alpha} + T^{\frac{\alpha-\beta}{2}}\Big(1+\big\Vert\widehat\zeta\hspace{2pt}\big\Vert_{{\boldsymbol{\mcN}}_T}^{{\color{red} 4}}\Big)\Big(1 + \norme{{\bf u}}_{\dab(X)}^2 + \mcW_{p,\scrC_T^\alpha}\big(\mu,\delta_{\bf 0}\big)\Big).
\end{align}
One can then find 
\[
M = M\Big(\norme{u_0}_\alpha+\hspace{1pt} \Vert\widehat\zeta\hspace{2pt}\Vert_{{\boldsymbol{\mcN}}_T} + \hspace{1pt}\mcW_{p,\scrC^\alpha_T}\big(\mu,\delta_{\bf 0}\big)\Big)
\]
and 
\[
T=T\Big(\norme{u_0}_\alpha+ \hspace{1pt}\Vert\widehat\zeta\hspace{1pt}\Vert_{{\boldsymbol{\mcN}}_T}+ \hspace{1pt} \mcW_{p,\scrC^\alpha_T}\big(\mu,\delta_{\bf 0}\big)\Big)
\]
such that the map $\Phi_{\widehat\zeta,u_0,\mu}$ sends the ball $\big\{ {\bf u}\in \mcE^\beta_T(Z) \,; \norme{{\bf u}}_{\mcE^\beta_T}\leq M \big\}$ into itself. For $\norme{{\bf u}}_{\mcE^\beta_T(Z)}\leq M$, Proposition \ref{schauder3} tells us that
\[
d_{\mcE^\beta_T(Z)}\big(\Phi_{\widehat\zeta,\mu}({\bf u}_1),\Phi_{\widehat\zeta,\mu}({\bf u}_2)\big) \lesssim_M T^{(\alpha-\beta)/2} d_{\dab}\big({\bf u}_1, {\bf u}_2\big),
\]
so choosing $T$ small enough ensures that the map $\Phi_{\widehat\zeta,\mu}$ has a unique fixed point in $\mcE^\beta_T(Z)$.

\medskip

{\it Step 2: Long-time well-posedness.} {\color{black} Following some recent long work of Chandra, Feltes \& Weber \cite{CFW} formulated in the language of regularity structures Shen, Zhu \& Zhu gave in \cite{ShenZhuZhu} a simple and short proof of the non-explosion of the solution to the equation
\[
(\partial_t - \Delta)v_t = f({\bf v}_t) \zeta_t
\]
on paracontrolled functions, with initial condition $u_0\in C^\alpha$, for $f\in C^2_b$. It is elementary to follow their proof when
\[
(\partial_t - \Delta)v_t = {\sf f}({\bf v}_t,\cdot) \zeta_t + g(v_t)
\]
for some ${\sf f}\in C^2_b(\bbR\times{\sf T}^2)$ and $g\in C^1_b$, and obtain the same a posteriori estimates that guarantee the non-explosion of the solution. (See Appendix \ref{section_appendix_nonexplosion} for a more complex situation. The regularity of $g$ is used for the local in time well-posedness result and the bounded character of $g$ is used for the a posteriori global in time estimates. With the notations of \cite{ShenZhuZhu}, insert $g(v_t)$ inside the right hand side of the dynamics of $u_1^\sharp$.) A close inspection of the proof of the main result of \cite{ShenZhuZhu}, making explicit the dependence of all the implicit multiplicative constants as some functions of ${\sf f}\in C^2_b(\bbR\times{\sf T}^2)$ and $g\in C^1_b(\bbR)$, shows that the analysis of \cite{ShenZhuZhu} also works for some equations
\begin{equation} \label{EqTimeDpdtEq}
(\partial_t - \Delta)v_t = {\sf f}({\bf v}_t, \cdot, t) \zeta_t + g(v_t,t)
\end{equation}
where ${\sf f}$ and $g$ are time-dependent with finite $C_TC^2(\bbR\times{\sf T}^2)$ and $C_TC^1$-norms respectively. Equation \eqref{EqTimeDpdtEq} is globally well-posed and there are some exponents $\kappa_1,\kappa_2, \kappa_3\in [1,\infty)$ independent of the time horizon $T$ from Step 1 such that 
\begin{equation} \label{EqSizeEstimateSolutionTimeDpdtEq}
\Vert {\bf u}\Vert_{\mcE_T^\beta(Z)} \lesssim \big(1 + \Vert\widehat\zeta\Vert_{{\boldsymbol{\mcN}}_T}\big)^{\kappa_1} \Big( 1 + \Vert {\sf f}\Vert_{C_TC^2} + \Vert g\Vert_{C_TC^1}\Big)^{\kappa_2} \big(1+\Vert u_0\Vert_{C^\alpha}\big)^{\kappa_3}.
\end{equation}
(We note that they use a time weight to measure $u$ in $\scrC^\alpha$. This is only used to deal with an initial condition in $L^\infty$, and one can work with an unweighted norm for $u$ if the initial condition is in $C^\alpha$, as here. Then the norms on the space of paracontrolled functions used in \cite{ShenZhuZhu} and the present work are equivalent.) The local in time well-posedness result from Step 1 and the a posteriori estimate \eqref{EqSizeEstimateSolutionTimeDpdtEq} together prove Proposition \ref{SolPAM1}.   }
\end{Dem}

\ssk

{\color{black} We note that \cite{CFW} gives some bound on the uniform norm of $u$ on $[0,T]\times {\sf T}^2$ that is linear in $\Vert u_0\Vert_{C^\alpha}$. The bound above is polynomial in $\Vert u_0\Vert_{C^\alpha}$ but the norm $\Vert {\bf u}\Vert_{\mcE_T^\beta(Z)}$ is stronger.   }

\ssk

We now specialize the result of Proposition \ref{SolPAM1} to the case where $\widehat\zeta$ is the random enhancement $\widehat\xi(\omega)$ of a random noise $\xi$ provided by Theorem \ref{ThmNoiseEnhancement}. Given $\epsilon>0$ we set 
\[
c^\epsilon(t) \defeq \bbE[X^{\epsilon}_t \odot \xi^{\epsilon}_t].
\]
Let us write $f'$ for the derivative $\partial_1f$ of $f$ with respect to its first argument. It is a consequence of Step 2 that the equation 

\begin{equation} \label{EqRenormalizedEqMuFixed}
(\partial_t - \Delta) u^{\epsilon_k}_t = f(u^{\epsilon_k}_t,\mu_t)\zeta^{\epsilon_k}_t - c^{\epsilon_k} (f'f)(u^{\epsilon_k}_t , \mu_t) + g(u^{\epsilon_k}_t ,\mu_t)
\end{equation}
is globally well-posed on $[0,T_0]\times {\sf T}^2$.

\ssk

\begin{lem} \label{LemASconvergenceSubsequence}
There is a sequence $\epsilon_k>0$ converging to $0$ such the $u^{\epsilon_k}$ converge almost surely in $\scrC^\alpha_{T_0}$ to $u$.
\end{lem}

\ssk

\begin{Dem}
The enhanced noise $\widehat\xi$ is the almost sure limit in ${\boldsymbol{\mcN}_{T_0}}$ of the sequence of enhanced smooth noises $\widehat\xi^{\epsilon_k} \defeq (\xi^{\epsilon_k} , \xi^{\epsilon_k}\odot X^{\epsilon_k} - c^{\epsilon_k})$, for an appropriate choice of sequence $(\epsilon_k)_{k\geq 1}$. It follows from Proposition \ref{SolPAM1} that the function $u$ is the limit in $\scrC^\alpha_{T_0}$ of the sequence $\widetilde u^{\epsilon_k}$ where $\widetilde{\bf u}^{\epsilon_k}$ is the solution to \eqref{EqEquationMuFixed} with the {\color{black} enhanced} noise $\widehat\xi^{\epsilon_k}$ in place of $\widehat{\zeta}$. But we have
\begin{align*}
f({\widetilde{\bf u}^{\epsilon_k}, \mu})\xi^{\epsilon_k} + g(\widetilde u^{\epsilon_k},\mu) &= f(\widetilde u^{\epsilon_k} ,\mu)\varolessthan \xi^{\epsilon_k} + \xi^{\epsilon_k} \varolessthan f(\widetilde u^{\epsilon_k} ,\mu)+ f(\widetilde u^{\epsilon_k} ,\mu)^\sharp\odot \xi^{\epsilon_k} 
\\
&\quad+{\sf C} \big( f(\widetilde {\bf u}^{\epsilon_k} ,\mu)',X_n ,\xi^{\epsilon_k} )+ (f'f)(\widetilde u^{\epsilon_k} ,\mu) \big(\xi^{\epsilon_k}\odot X^{\epsilon_k} - c^{\epsilon_k}\big) + g(\widetilde u^{\epsilon_k} ,\mu)
\\
&= f(\widetilde u^{\epsilon_k},\mu)\xi^{\epsilon_k} - c^{\epsilon_k}(f'f)(\widetilde u^{\epsilon_k},\mu) + g(\widetilde u^{\epsilon_k},\mu),
\end{align*}
so $\widetilde u^{\epsilon_k}$ is a solution of the equation
\[
(\partial_t - \Delta) \widetilde u^{\epsilon_k}_t = f(\widetilde{\bf u}^{\epsilon_k}_t,\mu_t)\xi^{\epsilon_k}_t - c^{\epsilon_k}(t)(f'f)( \widetilde u^{\epsilon_k}_t,\mu_t) + g(\widetilde u^{\epsilon_k}_t,\mu_t),
\]
and one has indeed $u^{\epsilon_k} = \widetilde u^{\epsilon_k}$ {\color{black} since they have the same initial condition}.
\end{Dem}

\ssk

{\color{black} We now proceed with the second step of the resolution of the mean field equation \eqref{EqMainEq} in the special case considered in this section, where $f$ is given by \eqref{EqFormulaFNiceKernel} in terms of $F$ and an interaction kernel $k\in C^2_b$.}

\ssk

\begin{thm} \label{ThmLongRange2ptfixe}
Fix $0<T_0<\infty$ and $1\leq p<\infty$. Assume $F\in C^3_b$ and $g$ satisfies the Lipschitz assumption \eqref{EqConditionLipG}. There exists $0<T\leq T_0$ with the following property. For every $u_0\in C^\alpha$ there exists a unique paracontrolled solution $\bf u$ to the mean field Equation \eqref{EqMainEq} in $\bbL^p(\Omega \,; \scrC^\alpha_{{\color{black} T}}\big)$. It is a locally Lipschitz continuous function of the enhanced noise $\widehat\xi\in \bbL^{{\color{black} r_1}}(\Omega\,;\boldsymbol{\mcN}_T)$, {\color{black} for some appropriate choice of exponent $r_1\in(0,\infty)$}. Furthermore {\color{black} the function} $u$ {\color{black} associated with $\bf u$} is the limit in $\bbL^p\big(\Omega \,; \scrC^\alpha_{{\color{black} T}}\big)$ of the solutions $u^\eps$ of the renormalized equations
\[
(\partial_t-\Delta)u^\eps = f\big(u^\eps,\mcL(u^\eps_t)\big)\zeta^\eps - c^\eps(t) (f'f)\big(u^\eps,\mcL(u^\eps_t)\big) + g\big(u^\eps,\mcL(u^\eps_t)\big)
\]
{\color{black} with initial condition $u_0$.}
\end{thm}

\ssk

{\color{black} We only obtain here a local in time well-posedness result for \eqref{EqMainEq} as we use a naive fixed point and small time contraction strategy to prove this statement. Since $u_0$ is fixed we do not record in the notations in the proof the dependence of the different objects on $u_0$. The precise value of the exponent $r_1$ does not matter so much in so far as we consider some noises that are actually in all the $\bbL^q(\Omega\,; \boldsymbol{\mcN}_{T_0})$ spaces.

We note as a preliminary remark to the proof of Theorem \ref{ThmLongRange2ptfixe} that if one writes $f$ in \eqref{EqFormulaFNiceKernel} in the form $f(a,\mu)(x)={\sf f}(a(x),x,\mu)$, with some obvious notation, then for all path $(\mu_t)_{0\leq t\leq T}$ of probability measures on $ \scrC^\alpha_T$, the quantity
\begin{equation} \label{EqMeasureUniformBoundNonlinearity}
\sup_{0\leq t\leq T} \Vert{\sf f}(\cdot,\cdot,\mu_t)\Vert_{C^2(\bbR\times{\sf T}^2)} \leq \Vert F\Vert_{C^2} \Vert k\Vert_{C^2}
\end{equation}
has an upper bound independent of $(\mu_t)_{0\leq t\leq T}$. On the other hand the regularity assumption \eqref{EqConditionLipG} on $g$ gives the path-dependent bound
\begin{equation} \label{EqBoundDriftTerm}
\sup_{0\leq t\leq T} \Vert{\sf g}(\cdot,\mu_t)\Vert_{C^1(\bbR)} \leq \sup_{0\leq t\leq T} \Vert \mu_t\Vert_{\mcP_p(C^\alpha)}.
\end{equation}   }

\ssk

\begin{Dem}
Pick $0<T\leq T_0$. Write ${\bf u}_{\widehat\xi}^{ \mu}$ for the {\color{black} paracontrolled} solution of \eqref{EqEquationMuFixed} {\color{black} with $\widehat{\xi}$ in place of $\widehat{\zeta}$ and denote by $u^\mu_{\widehat\xi}$ its associated function. We first check that $u^\mu_{\widehat\xi} \in \bbL^p(\Omega \,; \scrC^\alpha_T)$ using the pathwise bounds
\[
\norme{u^\mu_{\widehat\xi}}_{\scrC^\alpha_T} \lesssim (1+\norme{X}_{\scrC^\alpha_T})\norme{{\bf u}^\mu_{\widehat\xi}}_{\dab(X)},
\]
and \eqref{EqSizeEstimateSolutionTimeDpdtEq}. They imply with Cauchy-Schwarz inequality and the uniform upper bound \eqref{EqMeasureUniformBoundNonlinearity} and \eqref{EqBoundDriftTerm} that

\[
\Vert u^\mu_{\widehat\xi}\Vert_{\bbL^p(\Omega\,;\,\scrC^\alpha_T)}\lesssim \Vert \widehat{\xi} \, \Vert_{\bbL^{2p\kappa_1}(\Omega\,;\,{\boldsymbol{\mcN}}_T)}^{1+\kappa_1} \Big(1+\sup_{0\leq t\leq T} \Vert \mu_t\Vert_{\mcP_p(C^\alpha)}\Big)^{\kappa_2}.
\]   
The exponent $r_1$ in the statement of the theorem is chosen larger than $2p\kappa_1$. Recall we denote by $\mcL(Z)$ the law of a random variable $Z$. It follows from the previous estimate that we define a map $\Psi_{\widehat\xi}$ from $\mcP_p(\scrC^\alpha_T)$ into itself by setting}
\[
\Psi_{\widehat\xi} (\mu) \defeq \mcL\big(u_{\widehat\xi}^{ \mu}\big).
\]
Pick $A>0$. For $M$ sufficiently big and $T=T(M)$ even smaller, the map $\Psi_{\widehat\xi}$ sends the ball
\[
\Big\{ \mu \in \mcP_p(\scrC^\alpha_T) \,;\, \mcW_{p,\scrC^\alpha_T}(\mu,\delta_{\bf 0})\leq M\Big\}
\]
into itself. Now pick $\widehat\xi_1, \widehat\xi_2$ {\color{black} in $\bbL^{r_1}(\Omega \,; {{\boldsymbol{\mcN}}_T})$ and $\mu_1, \mu_2$ in $\mcP_p(\scrC^\alpha_T)$} such that
\begin{equation} \label{EqBallCondition}
\mcW_{p,\scrC^\alpha_T}(\mu_i , \delta_{\bf 0}) \leq M,
\end{equation}
for $i\in\{1,2\}$. Write ${\bf u}_i$ for $\Phi_{\widehat\xi_i}(\mu_i)$ and define the random variable
\[
A \defeq \big\Vert \widehat\xi_1\big\Vert_{{{\boldsymbol{\mcN}}_T}} + \big\Vert \widehat\xi_2\big\Vert_{{{\boldsymbol{\mcN}}_T}}.
\]
Combining Proposition \ref{schauder3}, Proposition \ref{paralinearisation0} and Proposition \ref{PropLocalLipschitzSharp}, we have
\begin{equation*} \begin{split}
d_{\mcE_T^\beta} \big({\bf u}_1, {\bf u}_2\big)  &\lesssim_A  T^{\frac{\alpha-\beta}{2}} \bigg\{ \big\Vert \widehat\xi_1 - \widehat\xi_2\big\Vert_{{{\boldsymbol{\mcN}}_T}} + d_{\mcE^\beta_T}({\bf u}_1, {\bf u}_2) + \mcW_{p,\scrC^\alpha_T}{\big({\mu}^1, {\mu}^2\big)} \bigg\}
\end{split} \end{equation*}
for some implicit positive multiplicative constant that is a polynomial of $A$. {\color{black} The exponent $r_1\in[1,\infty)$ in the statement of the theorem is chosen larger than two times the degree of that polynomial.} Integrating and using Cauchy-Schwarz inequality we get
\begin{align*}
\bbE\big[d_{\mcE^\beta_T}\big({\bf u}_1, {\bf u}_2\big)^{2p} \big]^2 \lesssim T^{2p(\alpha-\beta)} \Big\{   \bbE\big[ \big\Vert \widehat\xi_1 - \widehat\xi_2\big\Vert_{{{\boldsymbol{\mcN}}_T}}^{4p}\big] +  \bbE\big[ d_{\mcE^\beta_T}\big({\bf u}_1, {\bf u}_2\big)^{2p} \big]^2 + \mcW_{p,\scrC^\alpha_T}(\mu_1,\mu_2)^{4p} \Big\}
\end{align*}
so taking $T>0$ deterministic, small enough ensures that
\begin{align*}
\bbE\big[ d_{\mcE_T^\beta}\big({\bf u}_1, {\bf u}_2\big)^{2p} \big]^2 \lesssim T^{2p(\alpha-\beta)} \Big(\bbE\big[ \Vert \widehat\xi_1 - \widehat\xi_2\Vert^{4p}\big] + \mcW_{p,\scrC^\alpha_T}(\mu_1,\mu_2)^{4p}\Big).
\end{align*}
Now since
\[
\norme{u_1-u_2}_{\scrC^\alpha_T} \lesssim  \big(1+\norme{X_1}_{\scrC_T^\alpha}\big) \, d_{\dab}({\bf u}_1,{\bf u}_2)  + \norme{X_1-X_2}_{\scrC^\alpha_T} \norme{{\bf u}_2}_{\dab(X_2)}
\]
we obtain from Cauchy-Schwarz inequality the estimate
\[
\bbE\big[\norme{u_1-u_2}_{\scrC^\alpha_T}^p\big]^2 \lesssim  \big(1+\bbE[\norme{X_1}_{\scrC_T^\alpha}^{2p}]\big) \, \bbE\big[d_{\dab}({\bf u}_1,{\bf u}_2)^{2p}\big]  + \bbE\big[\norme{X_1-X_2}_{\scrC^\alpha_T}^{2p}\big] \, \bbE\big[\norme{{\bf u}_2}_{\dab(X_2)}^{2p}\big].
\]
{\color{black} Since
\[
\bbE\big[\norme{{\bf u}_2}_{\dab(X_2)}^{2p}\big] \lesssim \bbE\big[ \Vert\widehat{\xi}\Vert_{{\boldsymbol{\mcN}}_T}^{2p\kappa_1}\big] \Big(1+\sup_{0\leq t\leq T} \Vert \mu_t\Vert_{\mcP_p(C^\alpha)}\Big)^{2p\kappa_2}
\]
we see in the end that}
\begin{align*}
\mcW_{p,\scrC^\alpha_T}\big( \Psi_{\widehat\xi_1}(\mu_1) , \Psi_{\widehat\xi_2}(\mu_2) \big) \lesssim  \Vert \widehat{\xi}\,\Vert_{\bbL^{2p\kappa_1}(\Omega\,;\,{\boldsymbol{\mcN}}_T)}^{\kappa_1} (1+M)^{\kappa_2} \, \bbE\big[ \Vert \widehat\xi_1 - \widehat\xi_2\Vert_{{{\boldsymbol{\mcN}}_T}}^{4p}\big] + T^{\frac{\alpha-\beta}{2}}  \mcW_{p,\scrC^\alpha_T}(\mu_1,\mu_2).
\end{align*}
{\color{black} To make this estimate useful requires that we ask $r_1$ to be larger than or equal to $4p$.} In the end \eqref{EqMainEq} is well-posed in small time {\color{black} and its unique fixed point defines a locally Lipschitz function of $\widehat{\xi}\in \bbL^{r_1}(\Omega \,; {\boldsymbol{\mcN}}_T)$. The convergence part of the statement of Theorem \ref{ThmLongRange2ptfixe} is a consequence of Lemma \ref{LemASconvergenceSubsequence}.   }
\end{Dem}

\medskip

\section{Mean field equations with pointwise interaction in the diffusivity}
\label{SectionGeneral}

{\color{black} We dealt in Section \ref{SectionRegularInteraction} with a situation where the singular and mean field features of the equation \eqref{EqMainEq} did not interact. This is why we could use the classical tools first introduced by \cite{GIP} for the study of singular equations from the paracontrolled point of view.} We deal now with a large family of mean field type equations \eqref{EqMainEq} where the interaction in the diffusivity $f$ can be pointwise. {\color{black} In the example \eqref{EqFormulaFNiceKernel} of diffusivity $f$ it means working with $k(x,y)=\delta_x(y)$ with $\delta_x$ the Dirac distibution at $x$.} The enhancement of the noise needed to make sense of \eqref{EqMainEq} {\color{black} in that case} is specific to the mean field setting and described in Section \ref{SubsectionEnhancedNoise}. {\color{black} The term $f(u_t,\mcL(u_t))\xi_t$ then has a particular structure whose definition involves a paracontrolled  structure for $f(u_t,\mcL(u_t))$ different from the paracontrolled structures from Definition \ref{DefnSpacePCStructure} or Definition \ref{DefnParabolicPCStructure}. This specific structure is described in Section \ref{SubsectionParacStructure}. The random variable $u$ itself is still the function associated with an element $\bf u$ of the random space $\mcE^\beta(X)$.} The proper statement and proof of item {\it (a)} of Theorem \ref{ThmMainSketch} is done in Section \ref{SubsectionSolvingEquation}. {\color{black} The proof of item {\it (b)} in Theorem \ref{ThmMainSketch} is the object of Section \ref{SectionChaos}.   }

\ssk

{\color{black} {\it -- Assumptions on $f$ and $g$ in \eqref{EqMainEq}.}} We assume in this section that
\begin{equation} \label{InteractionPolynomial}
f(u,\mu)(x) = \int F\big(u(x), v_1(x), \dots, v_m(x)\big) \, \mu^{\otimes m}(dv_1\dots dv_m)
\end{equation}
for some integer $m\geq1$, for a function $F : \bbR^{m+1}\rightarrow\bbR$ of class {\color{black} $C^3_b$} -- or we assume that $f$ is a linear combination of such monomials. Notice that in \eqref{InteractionPolynomial} we evaluate the $v_i$ at the same point $x\in{\sf T}^2$ as $u$. With $m=1$, and compared to the regular interaction \eqref{EqLongRangeF} studied in Section \ref{SectionRegularInteraction}, the function \eqref{InteractionPolynomial} corresponds to a pointwise, non-regular, Dirac kernel 
\[
k(x,x')=\delta_x(x').
\]
{\color{black} As in Section \ref{SectionRegularInteraction} we assume that 
\begin{equation} \label{EqStructuralAssumptionG}
g(u,\mu)(z) = {\sf g}(u(z),\mu) \qquad (\forall\,z\in [0,T_0]\times {\sf T}^2)
\end{equation}
satisfies the Lipschitz-type condition \eqref{EqConditionLipG}.   }

\ssk

\subsection{Mean field enhancement of the noise$\boldmath{.}$   \hspace{0.15cm}}   
\label{SubsectionEnhancedNoise}

As above we work with the class of random Gaussian noises specified in Theorem \ref{ThmNoiseEnhancement}. The random field $\xi$ is initially defined on a probability space $(\Omega,\mcF,\mathbb{P})$. We extend it canonically as a random variable defined on the probability space $\big(\Omega^2, \mcF^{\otimes 2}, \mathbb{P}^{\otimes 2}\big)$ setting 
\[
\xi(\omega,\varpi) = \xi(\omega).
\] 
We also define 
\[
\overline{\xi}(\omega, \varpi) \defeq \xi(\varpi);
\]
this is under $\mathbb{P}^{\otimes 2}$ an independent copy of $\xi$. For a distribution $\Lambda$ on $\mathsf{T}^2$ and a positive regularization parameter $\epsilon$ we set 
\[
\Lambda^\epsilon \defeq \Lambda\circ P_\epsilon \in C^\infty
\]
{\color{black} with $P_\epsilon$ the heat operator at time $\epsilon$.} Recall $T_0$ stands for the time horizon that we use in our definition of the space of enhanced noises ${{\boldsymbol{\mcN}}_{T_0}}$ -- the interval $[0,T_0]$ is our maximal interval of time. Pick $1\leq p<\infty$. We define on $\big(\Omega^2, \mcF^{\otimes 2}, \mathbb{P}^{\otimes 2}\big)$ the {\color{black} $\scrC^\alpha_{T_0}$-valued} random variable
\[
\overline{X} \defeq \scrL^{-1}(\overline\xi).
\]
Note that $\bbE[(\xi^\epsilon\odot \overline{X}^\epsilon)(z)]=0$ for all $0<\epsilon\leq 1$ and all $z\in [0,T_0]\times{\sf T}^2$. The proof of Theorem \ref{ThmNoiseEnhancement} can be repeated verbatim to see that the random variables $\xi^\epsilon(\omega) \odot \scrL^{-1}(\overline{\xi}\hspace{1pt}^\epsilon(\varpi))$ converge in some $\bbL^{\color{black} r}(\bbP^{\otimes 2})$ to some limit random variable 
\[
\xi\odot\overline{X} \in \bbL^{\color{black} r}\big(\mathbb{P}^{\otimes 2}\,; C_{T_0}C^{2\alpha-2}\big),
\] 
as $\epsilon>0$ goes to $0$. We have for $\mathbb{P}$-almost every $\omega\in\Omega$ and $\mathbb{P}$-almost every $\varpi\in\Omega$
\[
\big\Vert \big(\xi\odot\overline{X}\big)(\omega,\cdot) \big\Vert_{\bbL^{\color{black} r}(\Omega \,; \,C_{T_0}C^{2\alpha-2})} < \infty
\]
and 
\[
\big\Vert \big(\xi\odot\overline{X}\big)(\cdot,\varpi) \big\Vert_{\bbL^{ \color{blue} r}(\Omega \,; \,C_{T_0}C^{2\alpha-2})} < \infty.
\]
We use below the notation $\overline{\bbE}$ to denote the (partial) expectation operator with respect to $\varpi$ on the product probability space.

\medskip

\begin{defn}\label{MeanFieldEnhancement}
The \textbf{mean field enhancement of the random noise $\xi$} is the random variable
\[
{\widehat{\xi}}^+(\omega, \varpi) \defeq \Big(\xi(\omega) , \big(\xi\odot X\big)(\omega), \overline{\xi}(\varpi), \big(\xi\odot\overline{X}\big)(\omega,\varpi)\Big) \in {{\boldsymbol{\mcN}}_{T_0}} \times {{\boldsymbol{\mcN}}_{T_0}},
\]
defined on $\big(\Omega^2, \mcF^2,\mathbb{P}^{\otimes 2}\big)$. We define on $(\Omega,\mcF,\bbP)$ the $\bbL^{16p}(\Omega \,; \bbR)$ random variable
\begin{equation} \label{EqDefnNormDoublyEnhancedNoise} \begin{split}
\llparenthesis \widehat\xi^+\rrparenthesis(\omega) &\defeq \norme{\xi(\omega)}_{C_{T_0}C^{\alpha-2}} + \big\Vert \xi^{(2)}(\omega)\big\Vert_{C_{T_0}C^{2\alpha-2}}    \\
&\quad+ \overline\bbE\Big[ \norme{\overline\xi(\omega,\cdot)}_{C_{T_0}C^{\alpha-2}}^4 \Big]^{\frac{1}{4}}+\overline\bbE\Big[ \norme{(\xi\odot\overline X)(\omega,\cdot)}_{C_{T_0}C^{2\alpha-2}}^4 \Big]^{\frac{1}{4}}.
\end{split} \end{equation}
\end{defn}

\ssk

Below we talk of ``our class of noises'' as the class of Gaussian random distributions specified in Theorem \ref{ThmNoiseEnhancement}.

\ssk

\subsection{Paracontrolled structure for mean field singular SPDEs}   \hspace{0.15cm}
\label{SubsectionParacStructure}

\ssk

{\it \S1. The mean field paracontrolled structure {\color{black} in space}.} The appropriate notion of space paracontrolled structure {\color{black} for the term $f(u_t,\mcL(u_t))$ in the right hand side of \eqref{EqMainEq}} is captured by the following definition.

\ssk

\begin{defn} \label{DefnMeanFieldPC}
Let $\Lambda : \Omega\rightarrow C^\alpha$ be an $\bbL^2$ random variable. A $C^\alpha$\textbf{-valued random variable} $\phi$ on $\Omega$ is said to be \textbf{$\omega$-paracontrolled by $\Lambda$} if there are some random variables 
\[
\delta_z\phi : \Omega \rightarrow C^\beta
\]
and 
\[
\delta_\mu \phi : \Omega \rightarrow \bbL^{\frac{4}{3}}\big(\Omega \,; C^\beta\big)
\] 
and 
\[
\phi^\sharp : \Omega \rightarrow C^{\alpha+\beta}
\]
such that one has
\begin{equation} \label{EqMeanFieldPCStructure}
\phi(\omega) = (\delta_z \phi)(\omega) \varolessthan \Lambda(\omega) + \overline{\bbE}\Big[(\delta_\mu \phi)(\omega,\cdot) \varolessthan \overline{\Lambda}(\cdot)\Big] + \phi^\sharp(\omega)
\end{equation}
for $\bbP$-almost all $\omega\in\Omega$, and 
\[
\big\Vert \delta_z\phi \big\Vert_{\bbL^2(\Omega \,;\, C^\beta)} + \big\Vert \delta_\mu \phi \big\Vert_{\bbL^2(\Omega \,;\, \bbL^{4/3}(\Omega\,;\,C^\beta)} + \big\Vert \phi^\sharp \big\Vert_{\bbL^2(\Omega \,;\,C^{\alpha+\beta})} < \infty.
\]
\end{defn}

\ssk

{\color{black} Let $\widehat\xi^+$ be the mean field enhancement} of the random noise $\xi$ {\color{black} and $\phi_t$ be $\omega$-paracontrolled by $X_t$. We set $\Phi_t = \big(\phi_t, \delta_z \phi_t, \delta_\mu \phi_t, \phi_t^\sharp\big)$ and we define}
\begin{equation}\label{eq_produit111}
(\Phi_t \xi_t)(\omega) \defeq \phi_t(\omega) \varolessthan \xi_t(\omega) + (\Phi_t \xi_t)^\sharp(\omega)
\end{equation}
with
\begin{equation}\label{eq_produit2} \begin{split}
(\Phi_t \xi_t)^\sharp(\omega) \defeq  &\; \xi_t(\omega) \varolessthan \phi_t(\omega) + \phi_t^\sharp(\omega)\odot \xi_t(\omega)   \\
&+ {\sf C}\big((\delta_z \phi_t)(\omega), X(\omega),\xi_t(\omega)\big) + \overline{\bbE}\Big[{\sf C}\big((\delta_\mu \phi_t)(\omega,\cdot), \overline{X}(\cdot),\xi_t(\omega)\big)\Big]   \\
&+ (\delta_z\phi_t)(\omega) \xi^{(2)}_t(\omega) + \overline{\bbE}\Big[(\delta_\mu \phi_t)(\omega,\cdot)\,\big(\xi\odot\overline{X}\big)(\omega,\cdot)\Big].
\end{split} \end{equation}
The proof of the next statement comes from some standard continuity estimates on the paraproduct and corrector operators and from H\"older inequality in the expectation $\overline\bbE$; it is left to the reader.

\ssk

\begin{prop} \label{PropProduct}
{\color{black} Let $\widehat{\xi}^+$ stand for the mean field enhancement of a noise $\xi$ in our class of noises. Assume we are given some random variables $\delta_z\phi, \delta_\mu \phi$ and $\phi^\sharp$ and $1\leq p<\infty$ such that
\begin{equation} \label{EqConditionDerivativesRemainder}
{\color{black} \Vert \delta_z \phi\Vert_{\bbL^p(\Omega\,;\,C_TC^\beta)} + \Vert \delta_\mu \phi\Vert_{\bbL^p(\Omega\,;\,C_T\bbL^{4/3}(\Omega;\,C^\beta))} + \Vert \phi^\sharp \Vert_{\bbL^p(\Omega\,;\,t^{\beta/2}C_TC^{\alpha+\beta})} < \infty. }
\end{equation}
We define $\phi$ from \eqref{EqMeanFieldPCStructure} with the random variables $\xi$ and $\overline{\xi}$ in the roles of $\Lambda$ and $\overline{\Lambda}$.} One has $\bbP$-almost surely $(\Phi \xi)(\omega)\in C_TC^{\alpha-2}$ and 
\[
\big\Vert (\Phi_t \xi_t)^\sharp(\omega)\big\Vert_{C^{\alpha-2+\beta}} \lesssim \Big(1 + \llparenthesis \widehat\xi^+\rrparenthesis(\omega)^2 \Big) \bigg( \norme{(\delta_z \phi_t)(\omega)}_{C^\beta} + \overline{\bbE}\big[ \big\Vert \delta_\mu \phi_t(\omega,\cdot)\big\Vert_{C^\beta}^{\frac{4}{3}}\big]^{\frac{3}{4}} + \big\Vert \phi_t^\sharp(\omega)\big\Vert_{C^{\alpha+\beta}}\bigg).
\]
Pick two noises $\xi^1, \xi^2$ in our class, with respective mean field enhancements $\widehat{\xi}^1{}^+, \widehat{\xi}^2{}^+$, {\color{black} and two sets of random variables $\delta_z \phi^i, \delta_\mu \phi^i, \phi^{i\sharp}$, for $i\in\{1,2\}$, both satisfying \eqref{EqConditionDerivativesRemainder}. One has $\bbP$-almost surely}

\begin{align*}
\big\Vert(&\Phi^1_t \xi^1_t)^\sharp(\omega) - (\Phi^2_t \xi^2_t)^\sharp(\omega)\big\Vert_{C^{\alpha-2+\beta}}   \\
&\lesssim (\star)_{12}(t, \omega)   \\
&\qquad\times \bigg(\big\Vert \delta_z \phi_t^1 - \delta_z \phi_t^2\big\Vert_{C^\beta} + \overline{\bbE}\Big[ \big\Vert\delta_\mu \phi_t^1 - \delta_\mu \phi_t^2\big\Vert_{C^\beta}^{\frac{4}{3}}\Big]^{\frac{3}{4}} + \big\Vert \phi_t^{1\sharp} - \phi_t^{2\sharp}\big\Vert_{C^{\alpha+\beta}} + \llparenthesis \widehat\xi^1{^+} - \widehat\xi^2{^+} \rrparenthesis(\omega) \bigg),
\end{align*}
for all $0<t\leq T$, where $(\star)_{12}(t, \omega)$ is a quadratic polynomial of 
\[
\underset{i\in\{1,2\}}{\max} \Big\{ \llparenthesis \widehat\xi^i{^+} \rrparenthesis(\omega), \Vert\delta_z \phi_t^i\Vert_{C^\alpha} , \overline\bbE\big[\Vert\delta_\mu \phi_t^i\Vert_{C^\alpha}^{\frac{4}{3}} \big]^{\frac{3}{4}} , \Vert \phi_t^{i\sharp}\Vert_{C^{\alpha+\beta}}  \Big\}.
\]
\end{prop}

\ssk

{\it \S2. A paralinearization formula.} {\color{black} We lift the function $f : C^\alpha\times\mcP_p(C^\alpha) \rightarrow C^\alpha$ from \eqref{InteractionPolynomial} into a $C^\alpha$-valued function on $C^\alpha\times \bbL^p(\Omega,\overline\bbP \,; C^\alpha)$ setting
\[
f_o(a\vert A) \defeq f\big(a,\mcL(A)\big),
\]
for $a\in C^\alpha, A\in \bbL^p\big(\Omega,\overline\bbP \,; C^\alpha\big)$. Given that $f_o$ is polynomial in its second argument} it will be useful below to work on the probability space 
\[
\big(\Omega^{m+1},\mcF^{\otimes(m+1)},\bbP^{\otimes(m+1)}\big)
\]
and write 
\[
(\omega,\omega_1,\dots,\omega_m)
\]
for an element of $\Omega^{m+1}$. We denote by $\overline \bbE^i$ the expectation operator with respect to the variable $\omega_i$, and for $I=(i_1,\dots,i_k)$ a subset of the integer interval $[\![1,m]\!]$ we write $\overline \bbE^I$ for the expectation operator with respect to the variables $(\omega_{i_1},\dots,\omega_{i_k})$. In those terms, and for $a\in C^\alpha, A\in \bbL^p(\Omega,\overline\bbP \,;  C^\alpha)$ and $\mu=\mcL(A)$, one has
\[
f(a, \mu)(x) = f_o\big(a \vert A\big)(x) = \overline\bbE^{[\![1,m]\!]}\Big[F\Big(a(x), A(\omega_1)(x), \dots, A(\omega_m)(x)\Big)\Big].
\]
As $F\in {\color{black} C^2_b}\subset C^1_b$ one has 
\[
\big\Vert F(a, A(\omega_1), \dots,A(\omega_m)\big\Vert_{C^\alpha}\lesssim 1 + \Vert a\Vert_{C^\alpha} + \sum_{j=1}^m\Vert A(\omega_j)\Vert_{C^\alpha}.
\]
Also, since $A\in C^\alpha$ is integrable the function $f_o(a\vert A)$ on ${\sf T}^2$ is indeed $\bbP^{\otimes (m+1)}$-almost surely an element of $C^\alpha$. For $i\in [\![1,m]\!]$ we set
\[
\partial_i f_o(a\vert  A)(x) \defeq \overline\bbE^{[\![1,m]\!]}\Big[(\partial_i F)\Big(a(x), A(\omega_1)(x),\dots,A(\omega_m)(x)\Big)\Big].
\]

{\color{black}
From the structural assumption \eqref{InteractionPolynomial} on $f$ there also exists a constant $L$ such that for every $a_1, a_2$ in $C^\alpha$ and $b_1, b_2$ in $\bbL^p(\Omega \,; C^\alpha)$ we have
\begin{equation} \label{EqLipConditionF} \begin{split}
 \big\Vert f_o(a_1 \vert b_1) - f_o(a_2 \vert b_2) \big\Vert_{C^\alpha} \leq L\Big(\Vert a_1 - a_2\Vert_{C^\alpha} + \overline\bbE\big[ \Vert b_1 - b_2\Vert_{C^\alpha}^p \big]^{\frac{1}{p}}\Big).   
\end{split} \end{equation} }

\ssk

\begin{prop} \label{PropNonlinearity}
Fix $t>0$ and assume we are given two random variables ${\bf h}_t = (h_t', h_t^\sharp)$ and ${\bf k}_t = (k'_t, k_t^\sharp)$ in $\bbL^{{\color{black} p}}(\Omega \,; \mcD^{{\color{black} \beta}}(X_t))$, with corresponding $C^\alpha$ functions $h_t, k_t$ on $\textsf{\emph{T}}^2$. Then $f_o({\bf h}_t \vert {\bf k}_t)$ is $\omega$-paracontrolled by $X_t$ in the sense of Definition \ref{DefnMeanFieldPC}, with 
\[
\big(\delta_zf_o\big)(h_t \vert k_t)(\omega) = (\partial_1f_o)\big(h_t(\omega) \vert k_t\big)h '(\omega)
\]
and 
\begin{equation*} \begin{split}
\big(&\delta_\mu f_o\big)({\bf h}_t \vert {\bf k}_t)(\omega, \varpi)   \\
&= \sum_{j=1}^m \overline\bbE^{[\![1,m]\!]\backslash\{j\}}\bigg[\big(\partial_{j+1}F\big)\Big(h_t(\omega), k_t(\omega_1), \dots, k_t(\omega_{j-1}), k_t(\varpi), k_t(\omega_{j+1}), \dots, k_t(\omega_m)\Big)\bigg] k_t'(\varpi),
\end{split} \end{equation*}
and
\begin{align*}
\big\Vert f_o({\bf h}_t(\omega) \vert {\bf k}_t)^\sharp\big\Vert_{C^{\alpha+\beta}}&\lesssim \bigg(1 + \norme{X_t(\omega)}^2_{C^\alpha} + \overline\bbE\big[\norme{\overline X_t}^4_{C^\alpha}\big]^{\frac{1}{2}}\bigg)
\\ &\quad\times\bigg(1+\norme{h_t'(\omega)}_{C^\beta} + \norme{h_t^\sharp(\omega)}_{C^\alpha} + \overline{\bbE}\big[\norme{k_t'}^4_{C^\beta}\big]^{\frac{1}{4}} + \overline{\bbE}\big[\norme{k_t^\sharp}^4_{C^\alpha}\big]^{\frac{1}{4}}    \bigg)
\\
&\quad \times \bigg(\hspace{-0.05cm}1+ \norme{h_t'(\omega)}_{C^\beta} + \norme{h_t^\sharp(\omega)}_{C^{\alpha+\beta}}+\overline{\bbE}\big[\norme{k_t'}_{C^\beta}^4\big]^{\frac{1}{4}} \hspace{-0.05cm}+ \overline{\bbE}\big[\norme{k_t^\sharp}_{C^{\alpha+\beta}}^4\big]^{\frac{1}{4}} \bigg) .
\end{align*}
For {\color{black} $X_t^1, X_t^2$ associated with two enhanced noises $\widehat\xi^{1+}, \widehat\xi^{2+}$}, for ${\bf h}_t^i$ and ${\bf k}_t^i$ in $\bbL^{4}(\Omega \,; \mcD^\alpha(X^i_t))$, for $i\in\{1,2\}$, we have 
\begin{equation} \label{EqfSharpIncrement} \begin{split}
&\big\Vert f_o\big({\bf h}_t^1(\omega) \vert {\bf k}_t^1\big)^\sharp - f_\circ\big({\bf h}_t^2(\omega) \vert {\bf k}_t^2\big)^\sharp \big\Vert_{C^{\alpha+\beta}} \lesssim (\star)_{12}(t, \omega) \times  \\
&\quad\bigg\{ \norme{X_t^1(\omega) - X_t^2(\omega)}_{C^\alpha} + \overline\bbE\big[\big\Vert\overline X_t^1 - \overline X_t^2\big\Vert_{C^\alpha}^4\big]^{\frac{1}{4}} + d_{\mcD^\beta}\big({\bf h}_t^1(\omega) , {\bf h}_t^2(\omega)\big) + \overline\bbE\big[d_{\mcD^\beta}\big({\bf k}_t^1 , {\bf k}_t^2\big)^4\big]^{\frac{1}{4}}\bigg\},
\end{split} \end{equation}
where $(\star)_{12}(t, \omega)$ is a polynomial function of 
\[
\underset{i\in\{1,2\}}{\max} \Big\{ \norme{X_t^i(\omega)} _{C^\alpha}, \, \overline\bbE\big[\Vert X_t^i\Vert_{C^\alpha}^4\big]^{\frac{1}{4}}, \, \big\Vert {\bf h}_t^i(\omega)\big\Vert_{\mcD^\alpha} , \, \overline\bbE\big[\Vert {\bf k}_t^i\Vert_{\mcD^\alpha}^4\big]^{\frac{1}{4}} \Big\}.
\]
\end{prop}

\ssk

{\color{black} It follows from Proposition \ref{PropProduct} and Proposition \ref{PropNonlinearity} that the random variable $f_o({\bf h}_t \vert {\bf k}_t)\xi_t$ is well-defined and is an element of $\bbL^p(\Omega, \bbP\,;\mcD(\xi_t))$. }

\ssk

\begin{Dem}
One has from paralinearisation
\begin{equation*} \begin{split}
F\big(h_t&(\omega), k_t(\omega_1),\dots, k_t(\omega_m)\big)   \\
&=\partial_1F\big(h_t(\omega), k_t(\omega_1),\dots, k_t(\omega_m)\big) \varolessthan h_t(\omega) + \sum_{j=1}^m\partial_{j+1}F\big(h_t(\omega), k_t(\omega_1),\dots, k_t(\omega_m)\big) \varolessthan k_t(\omega_j)
\\
&\quad+ R_F\big(h_t(\omega), k_t(\omega_1),\dots, k_t(\omega_m)\big)
\end{split} \end{equation*}
\begin{equation*} \begin{split}
&=\big( \partial_1F\big(h_t(\omega), k_t(\omega_1),\dots, k_t(\omega_m)\big)h_t'(\omega) \big) \varolessthan X_t(\omega)   \\
&\quad+  \sum_{j=1}^m\Big(\partial_{j+1} F\big(h_t(\omega), k_t(\omega_1),\dots, k_t(\omega_m)\big)k_t'(\omega_j)\Big) \varolessthan \overline X_t(\omega_j) + R_F + R_0 + \sum_{j=1}^m R_j
\end{split} \end{equation*}
where $R_F = R_F\big(h_t(\omega), k_t(\omega_1),\cdots, k_t(\omega_m)\big)\in C^{\alpha+\beta}$ and 
\begin{align*}
R_0&\defeq \Big\{\partial_1F\big(h_t(\omega), k_t(\omega_1), \dots, k_t(\omega_m)\big) \varolessthan\big(h'(\omega) \varolessthan X_t(\omega)\big)   \\
&\hspace{5.2cm}- \big(\partial_1F\big(h_t(\omega), k_t(\omega_1), \dots, k_t(\omega_m)\big)h_t'\big) \varolessthan X_t(\omega) \Big\}   \\
&\quad+ \partial_1F\big(h_t(\omega), k_t(\omega_1),\cdots, k_t(\omega_m)\big) \varolessthan h_t^{\sharp}(\omega),
\end{align*}
and $R_j$ is given by the explicit formula

\vfill \pagebreak

\begin{align*}
R_j&\defeq \Big\{\partial_{j+1} F\big(h_t(\omega), k_t(\omega_1),\cdots, k_t(\omega_m)\big) \varolessthan \big( k_t'(\omega_j) \varolessthan \overline X_t(\omega_j)\big)
\\
&\hspace{5.2cm}-\Big(\partial_{j+1} F\big(h_t(\omega), k_t(\omega_1),\cdots, k_t(\omega_m)\big) k_t'(\omega_j)\Big) \varolessthan \overline X_t(\omega_j) \Big\}
\\
&\quad+ \partial_{j+1} F\big(h_t(\omega), k_t(\omega_1),\cdots, k_t(\omega_m)\big) \varolessthan k_t^\sharp(\omega_j).
\end{align*}
From classical results in paradifferential calculus we have
\begin{equation*} \begin{split}
\norme{R_F}_{C^{\alpha+\beta}} &\lesssim \norme{F}_{C^2}\Big(1+\norme{h_t(\omega)}_{C^\alpha}^2+\sum_{j=1}^m \norme{k_t(\omega_j)}^2_{C^\alpha}\Big)   \\
&\lesssim \Big(1 + \norme{X_t(\omega)}^2_{C^\alpha} + \sum_{j=1}^m\norme{\overline X_t(\omega_j)}^2_{C^\alpha}\Big)   \\
&\hspace{2cm}\times\Big(1+\norme{h_t'(\omega)}_{C^\beta}^2 + \norme{h_t^\sharp(\omega)}_{C^\alpha}^2 + \sum_{j=1}^m\norme{k_t'(\omega_j)}_{C^\beta}^2 + \norme{k_t^\sharp(\omega_j)}_{C^\alpha}^2\Big),
\end{split} \end{equation*}
and
\begin{equation*} \begin{split}
\norme{R_0}_{C^{\alpha+\beta}} &\lesssim \big\Vert\partial_1F\big(h_t(\omega), k_t(\omega_1),\cdots, k_t(\omega_m)\big)\big\Vert_{C^\alpha}\Big(\norme{h_t'(\omega)}_{C^\beta} \norme{X_t(\omega)}_{C^\alpha} + \norme{h_t^\sharp(\omega)}_{C^{\alpha+\beta}}\Big)   \\
&\lesssim \Big(1+\norme{h_t(\omega)}_{C^\alpha} + \sum_{j=1}^m\norme{k_t(\omega_j)}_{C^\alpha}\Big)\Big(\norme{h_t'}_{C^\beta}\norme{X_t(\omega)}_{C^\alpha} + \norme{h_t^\sharp(\omega)}_{C^{\alpha+\beta}}\Big)   \\
&\lesssim \Big(1+\norme{X_t(\omega)}^2_{C^\alpha} + \sum_{j=1}^m\norme{\overline{X}_t(\omega_j)}^2_{C^\alpha}\Big)   \\
&\qquad\times \Big(1+\norme{h_t'(\omega)}_{C^\beta} + \norme{h_t^\sharp(\omega)}_{C^\alpha} + \sum_{j=1}^m \norme{k_t'(\omega_j)}_{C^\beta} + \norme{k_t^\sharp(\omega_j)}_{C^\alpha}    \Big)
\\
&\qquad\times \Big(1+ \norme{h_t'(\omega)}_{C^\beta} + \norme{h_t^\sharp(\omega)}_{C^{\alpha+\beta}} + \sum_{j=1}^m \norme{k_t(\omega_j)}_{C^\beta} + \norme{k_t^\sharp(\omega_j)}_{C^\alpha} 
  \Big),
\end{split} \end{equation*}
and, for $1\leq i\leq m$, we have for $\norme{R_i}_{C^{\alpha+\beta}}$ the upper bound
\begin{align*}
&\Big(1+\norme{X_t(\omega)}^2_{C^\alpha}+\sum_{j=1}^m\norme{\overline{X}_t(\omega_j)}^2_{C^\alpha}\Big) 
\\
&\quad\times \bigg\{1+\norme{h_t'(\omega)}_{C^\beta} + \norme{h_t^\sharp(\omega)}_{C^\alpha} + \sum_{j=1}^m \norme{k_t'(\omega_j)}_{C^\beta} + \norme{k_t^\sharp(\omega_j)}_{C^\alpha} \bigg\}
\\
&\quad \times \bigg\{1+ \norme{h_t'(\omega)}_{C^\beta} + \norme{h_t^\sharp(\omega)}_{C^\alpha} + \norme{k_t^\sharp(\omega_i)}_{C^{\alpha+\beta}} + \sum_{j=1}^m\norme{k_t'(\omega_j)}_{C^\beta} + \norme{k_t^\sharp(\omega_j)}_{C^\alpha} \bigg\}.
\end{align*}
So we have for $\big\Vert R_F+\sum_{j=0}^m R_j\big\Vert_{C^{\alpha+\beta}}$ the bound
\begin{align*}
 &\Big(1+\norme{X_t(\omega)}^2_{C^\alpha} + \sum_{j=1}^m\norme{\overline{X}_t(\omega_j)}^2_{C^\alpha}\Big) \bigg(1+\norme{h_t'}_{C^\beta}+\norme{h^\sharp}_{C^\alpha}+\sum_{j=1}^m \norme{k_t'(\omega_j)}_{C^\beta} + \norme{k_t^\sharp(\omega_j)}_{C^\alpha} \bigg)
\\
&\times \bigg(1+ \norme{h_t'}_{C^\beta} + \norme{h_t^\sharp}_{C^{\alpha+\beta}} + \sum_{j=1}^m\norme{k_t'(\omega_j)}_{C^\beta} + \norme{k_t^\sharp(\omega_j)}_{C^{\alpha+\beta}}
\bigg).
\end{align*}
Taking the $\overline\bbE^{[\![1,m]\!]}$ expectation one gets
\begin{align*}
&f_o({\bf h}_t(\omega) \vert {\bf k}_t)   \\
&=\big(\partial_1f_o(h_t(\omega) \vert k_t) h_t'(\omega)\big) \varolessthan X_t(\omega) 
\\
&\quad+  \overline\bbE^{[\![1,m]\!]}\bigg[ \sum_{j=1}^m \Big(\big(\partial_{j+1}F\big)(h_t(\omega), k_t(\omega_1), \dots, k_t(\omega_m)) k_t'(\omega_j)\Big) \varolessthan \overline X_t(\omega_j) \bigg] + f_o({\bf h}_t(\omega) \vert {\bf k}_t)^\sharp   \\
&= \big(\partial_1f_o(h_t(\omega) \vert k_t) h_t(\omega)'\big) \varolessthan X_t(\omega) \, + \sum_{j=1}^m   \\
&\overline\bbE\bigg[ \overline\bbE^{[\![1,m]\!]\backslash\{j\}}\bigg[\big(\partial_{j+1}F\big)\Big(h_t(\omega) , k_t(\omega_1),\cdots, k_t(\omega_{j-1}), k_t(\varpi), k_t(\omega_{j+1}),\dots, k_t(\omega_m)\Big) k_t'(\varpi)\bigg] \varolessthan \overline X_t(\varpi) \bigg]   \\
&\quad+ f_o({\bf h}_t(\omega) \vert {\bf k}_t)^\sharp,
\end{align*}
with
\begin{equation*} \begin{split}
\big\Vert f_o({\bf h}_t(\omega) \vert {\bf k}_t)^\sharp \big\Vert_{C^{\alpha+\beta}} &\lesssim \Big(1 + \norme{X_t(\omega)}^2_{C^\alpha} + \overline\bbE\big[\norme{\overline X_t}^4_{C^\alpha}\big]^{\frac{1}{2}}\Big)   \\ 
&\quad\times\bigg(1+\norme{h_t'(\omega)}_{C^\beta} + \norme{h_t^\sharp(\omega)}_{C^\alpha} + \overline{\bbE}\big[\norme{k_t'}^4_{C^\beta}\big]^{\frac{1}{4}} + \overline{\bbE}\big[\norme{k_t^\sharp}^4_{C^\alpha}\big]^{\frac{1}{4}} \bigg)   \\
&\quad \times \bigg(1+ \norme{h_t'(\omega)}_{C^\beta} + \norme{h_t^\sharp(\omega)}_{C^{\alpha+\beta}} + \overline{\bbE}\big[\norme{k_t'}_{C^\beta}^4\big]^{\frac{1}{4}}  + \overline{\bbE}\big[\norme{k_t^\sharp}_{C^{\alpha+\beta}}^4\big]^{\frac{1}{4}} \bigg).
\end{split} \end{equation*}
This concludes the proof of the proposition.
\end{Dem}

\medskip

\subsection{Solving Equation \eqref{EqMainEq}$\boldmath{.}$   \hspace{0.15cm}}
\label{SubsectionSolvingEquation}
    
{\color{black} Fix a time horizon $0<T_0<\infty$, an integrability exponent $4\leq p<\infty$ and $u_0\in C^\alpha$.} Recall from \eqref{EqDefnLc} the definition of the lifting maps ${\sf L}_c$, for $c\in C([0,T_0],\bbR)$, and the existence of a {\color{black} sequence $(\epsilon_k)$ decreasing to $0$ and some} functions $c^{\epsilon_k}\in C([0,T_0],\bbR)$ such that the random variables ${\sf L}_{c^{\epsilon_k}}(\xi^{\epsilon_k})$ are converging in {\color{black} $\bbL^{r}(\Omega,C_{T_0}C^{2\alpha-2})$} to the random variable $\scrL^{-1}(\xi)\odot\xi$ {\color{black} for all $1\leq r<\infty$. }

\ssk

\begin{prop} \label{PropSolPAM}
For every ${\bf v}\in \bbL^p\big(\Omega, \bbP \,; \mcE^\beta_{T_0}(X)\big)$ {\color{black} there exists} a unique solution {\color{black} ${\bf u}_{\widehat\xi^+,{\bf v}}\in\mcE^\beta_{T_0}(X)$} to the equation 
\begin{equation} \label{EqPAM}
(\partial_t-\Delta)u_t = f_o({\bf u}_t \vert {\bf v}_t)\,\xi_t + g_o(u_t \vert v_t)
\end{equation}
{\color{black} with initial condition $u_0$}. It satisfies almost surely the local Lipschitz continuity property
\begin{equation} \label{EqContinuityProperty}
d_{\mcE^\beta_{T_0}}\big({\bf u}_{\widehat\xi^+_1,{\bf v}_1}(\omega) , {\bf u}_{\widehat\xi^+_2,{\bf v}_2}(\omega)\big) \leq M'(\omega) \Big( \overline\bbE\big[\Vert {\bf v}_1 - {\bf v}_2\Vert_{\bbL^p(\Omega \,; \,\mcE^\beta_{T_0})}\big]  + \llparenthesis \widehat\xi^+_1 - \widehat\xi^+_2\rrparenthesis(\omega)\Big)
\end{equation}
for $M'(\omega)$ a non-increasing function of $\max_{i=1,2}\Big\{\overline\bbE\big[ \Vert{\bf v}_i\Vert_{\mcE^\beta_{T_0}}^p\big] , \llparenthesis \widehat\xi_i^+\rrparenthesis(\omega) \Big\}$. The random variable $u_{\widehat\xi^+,{\bf v}}(\cdot)\in\scrC^{\alpha}_{T_0}$ on $\Omega$ associated with ${\bf u}_{\widehat\xi^+,{\bf v}}(\cdot)$ is the limit in probability of the solutions $u^{\epsilon_k}$ of the equations
\begin{equation} \label{EqRenormalizedEquation}
(\partial_t-\Delta) u_t^{\epsilon_k} = f_o(u_t^{\epsilon_k} \vert {\bf v}_t)\,\xi^{\epsilon_k}_t + g_o(u_t^{\epsilon_k} \vert {\bf v}_t) - c^{\epsilon_k}(t) \big(f'_o f_o\big)(u_t^{\epsilon_k} \vert {\bf v}_t)
\end{equation}
with initial condition $u_0$.
\end{prop}

\ssk

\begin{Dem}
We proceed in three steps.

\ssk

{\it Step 1: Local in time well-posedness.} Rewrite \eqref{EqPAM} as the fixed point equation
\[
u_t = P_t(u_0) + \int_0^t P_{t-s}\big(f_o({\bf u}_s \vert {\bf v}_s)\xi_s + g(u_s \vert v_s)\big) ds.
\]
We get from Proposition \ref{PropProduct} and Proposition \ref{PropNonlinearity} that $f_o({\bf u}_s \vert {\bf v}_s)\xi_s + g_o(u_s, v_s)$ is for each $s$ an element of $\mcD^\alpha(\xi_s)$ with Gubinelli derivative {\color{black} $f(u_s \vert v_s)$} and remainder $\big(f_o({\bf u}_s \vert {\bf v}_s)\xi_s\big)^\sharp + g_o(u_s \vert v_s)$. 

$\bullet$ We check first that {\color{black} for any ${\bf u}\in\mcE^\beta_T\big(X(\omega)\big)$ the quantities} $f_o(u\vert v)\in\mathscr{C}^\alpha_T$ and $\big(f_o({\bf u}_s \vert {\bf v}_s)\xi_s\big)^\sharp + g_o(u_s\vert v_s)$ satisfy the bound \eqref{EqConditionRemainderSize} involved in the Schauder type estimate of Proposition \ref{schauder3}. Recall from \eqref{EqDefnNormDoublyEnhancedNoise} the definition of the mixed pathwise/averaged random variable $\llparenthesis\widehat\xi^+\rrparenthesis(\omega)$. Take {\color{black} any $0<T\leq T_0$ and} ${\bf u}\in\mcE^\beta_T(X)$. First, {\color{black} from the structural assumption \eqref{InteractionPolynomial} on $f$ where $F\in C^2_b$, and from the parabolic paracontrolled structure of $u$ and $v$ from Definition \ref{DefnParabolicPCStructure}, one has the estimate}
\begin{align*}
\norme{f_o(u \vert v)}_{\scrC^\alpha_T} &\lesssim 1+ \norme{u}_{\scrC^\alpha_T} + \overline{\bbE}\big[ \norme{v}_{\scrC^\alpha_T} \big] \lesssim \big(1+\llparenthesis\widehat\xi^+\rrparenthesis\big) \Big(1+ \norme{{\bf u}}_{\dab} + \overline\bbE\big[\norme{{\bf v}}_{\dab}^2 \big]^{\frac{1}{2}}   \Big).
\end{align*}
Second, combining the estimates from Lemmas \ref{PropProduct} and \ref{PropNonlinearity}, one gets at some fixed time $t$ the estimates
\begin{align*}
\big\Vert (f_o({\bf u}_t \vert {\bf v}_t) &\xi_t)^\sharp\big\Vert_{C^{\alpha-2+\beta}}   \\
&\lesssim \big(1 + \llparenthesis\widehat\xi^+\rrparenthesis^2\big)  \Big(\hspace{-2pt}\norme{\delta_zf(u_t, v_t)}_{C^\beta} + \overline{\bbE}\big[\norme{\delta_\mu f(u_t, v_t)}_{C^\beta}^{\frac{4}{3}}\big]^{\frac{3}{4}} + \norme{f({\bf u}_t, {\bf v}_t)^\sharp}_{C^{\alpha+\beta}} \hspace{-2pt} \Big)
\\
&\lesssim \big(1 + \llparenthesis\widehat\xi^+\rrparenthesis^2\big) \bigg\{ \Big( 1+\norme{u_t}_{C^\alpha}+\overline{\bbE}\big[\norme{v_t}_{C^\alpha} \big]   \Big)\norme{u'_t}_{C^\beta}
\\
&\hspace{1.3cm}+ \Big( 1+\norme{u_t}_{C^\alpha} + \overline{\bbE}\big[\norme{v_t}_{C^\alpha}^2 \big]^{\frac{1}{2}}\Big) \overline{\bbE}\big[\norme{v'_t}^4_{C^\beta}\big]^{\frac{1}{4}} 
+ \norme{f({\bf u}_t, {\bf v}_t)^\sharp}_{C^{\alpha+\beta}} \bigg\}
\\
&\lesssim \big(1+\llparenthesis\widehat\xi^+\rrparenthesis^4\big) \bigg(1+\norme{u'_t}_{C^\beta}\hspace{-1pt} + \norme{u^\sharp_t}_{C^\alpha}\hspace{-2pt} + \overline{\bbE}\big[\norme{v'_t}^4_{C^\beta}\big]^{\frac{1}{4}}\hspace{-2pt} + \overline{\bbE}\big[\norme{v^\sharp_t}^4_{C^\alpha}\big]^{\frac{1}{4}} \bigg)
\\
&\hspace{1.3cm}\times \bigg(1+ \norme{u'_t}_{C^\beta} + \norme{u^\sharp_t}_{C^{\alpha+\beta}} + \overline{\bbE}\big[\norme{v'_t}_{C^\beta}^4\big]^{\frac{1}{4}} + \overline\bbE\big[ \norme{v^\sharp_t}_{C^{\alpha+\beta}}^4\big]^{\frac{1}{4}} \bigg).
\end{align*}
It follows from this bound that we have
\[
\sup_{t\in(0,T]}t^{\beta/2} \norme{(f_o({\bf u}_t \vert {\bf v}_t)\xi_t)^\sharp}_{C^{\alpha-2+\beta}} \lesssim \big(1 + \llparenthesis\widehat\xi^+\rrparenthesis^4\big)\Big(1+\norme{{\bf u}}^2_{\dab}+\overline\bbE\big[ \norme{{\bf v}}^4_{\dab}\big]^{\frac{1}{2}}\Big).
\]
We also have {\color{black} from the structural assumption \eqref{EqStructuralAssumptionG} on $g$ and the assumption \eqref{EqConditionLipG} on $\sf g$ that}
\begin{align*}
\sup_{t\in(0,T]}t^{\beta/2}\norme{g_o(u_t \vert v_t)}_{C^{\alpha-2+\beta}}&\lesssim {\color{black} \sup_{t\in(0,T]}t^{\beta/2}\norme{g_o(u_t \vert v_t)}_{L^\infty}}   \\
&\lesssim \sup_{t\in(0,T]}t^{\beta/2} \Big(1 + \norme{u_t}_{C^\alpha} + \overline{\bbE}\big[\norme{v_t}_{C^\alpha}^2\big]^{\frac{1}{2}}\Big)
\\
&\lesssim \big(1 + \llparenthesis\widehat\xi^+\rrparenthesis\big) \Big(1 + \norme{{\bf u}}_{\dab} + \overline\bbE\big[\norme{{\bf v}}_{\dab}^2\big]^{\frac{1}{2}}\Big).
\end{align*}
We have in the end the pathwise estimate
\[
\sup_{t\in(0,T]}t^{\beta/2} \big\Vert (f_o({\bf u}_t \vert {\bf v}_t)\xi_t)^\sharp + g_o(u_t \vert v_t) \big\Vert_{C^{\alpha-2+\beta}} \lesssim \big(1 + \llparenthesis\widehat\xi^+\rrparenthesis^4\big)\Big(1 + \norme{{\bf u}}^2_{\dab} + \overline\bbE\big[ \norme{{\bf v}}^4_{\dab}\big]^{\frac{1}{2}}\Big).
\]
$\bullet$ It follows from Proposition \ref{schauder3} that the map 
\[
\Phi_{\widehat\xi^+,{\bf v}} : \mcE^\beta_T\big(X(\omega)\big) \rightarrow\mcE^\beta_T\big(X(\omega)\big)
\]
which associates to ${\bf u}\in\dab(X(\omega))$ the solution $w$ of the equation {\color{black} on the time interval $[0,T]$}
\[
(\partial_t - \Delta) w = f_o({\bf u}_t \vert {\bf v}_t) \xi_t + g_o(u_t \vert v_t)
\]
with initial condition $w_0=u_0 \in C^\alpha$, is well defined and satisfies the bound
\begin{align*}
\big\Vert\Phi_{\widehat\xi^+,{\bf v}}({\bf u})\big\Vert_{\mcE^\beta_T} &\lesssim  \norme{u_0}_{C^\alpha} + T^{(\alpha-\beta)/2}\big(1 + \llparenthesis\widehat\xi^+\rrparenthesis^4\big) \Big(1+\norme{{\bf u}}_{\dab}^2+ \overline\bbE\big[\Vert{\bf v}\Vert^4_{\dab}\big]^{\frac{1}{2}}\Big).
\end{align*}
Recall $4\leq p<\infty$. One can then find some positive random variables
\[
M(\omega) = M\Big(\overline\bbE\big[ \Vert{\bf v}\Vert_{\dab}^p\big] , \llparenthesis\widehat\xi^+\rrparenthesis(\omega)\Big),  \qquad   T_1(\omega) = T_1\Big(\overline\bbE\big[ \Vert{\bf v}\Vert_{\dab}^p\big] , \llparenthesis\widehat\xi^+\rrparenthesis(\omega)\Big)
\] 
such that the map $\Phi_{\widehat\xi^+,{\bf v}}$ sends the ball 
\[
\Big\{ {\bf u}\in \mcE^\beta_T\big(X(\omega)\big) \,; \norme{{\bf u}}_{\mcE^\beta_T} \leq M(\omega) \Big\}
\]
into itself. Now, given $\widehat\xi_1^+,\widehat\xi_2^+$ in $\bbL^{8p}(\Omega^2 \,; {{\boldsymbol{\mcN}}_T}\times{{\boldsymbol{\mcN}}_T})$ and ${\bf v}_1, {\bf v}_2$ in $\bbL^p\big(\Omega \,; \mcE^\beta_{T_0}(X(\omega))\big)$, we define a random variable
\[
M'(\omega) = M'\Big(\max_{i=1,2}\Big\{\overline\bbE\big[ \Vert{\bf v}_i\Vert_{\mcE^\beta_{T_0}}^p\big] , \llparenthesis \widehat\xi_i^+\rrparenthesis(\omega) \Big\}\Big) \geq M(\omega).
\]
For $\norme{{\bf u}}_{\dab}\leq M'(\omega)$ Proposition \ref{schauder3} tells us that
\begin{equation*} \begin{split}
d_{\mcE^\beta_T}\Big(\Phi_{\widehat\xi^+_1,{\bf v}_1}({\bf u}_1) &, \Phi_{\widehat\xi^+_2,{\bf v}_2}({\bf u}_2)\Big)   \\
&\lesssim_{M'(\omega)} T^{(\alpha-\beta)/2}\Big\{ d_{\dab}({\bf u}_1, {\bf u}_2)  + \Vert {\bf v}_1-{\bf v}_2\Vert_{\bbL^p(\Omega \,; \mcE^\beta_T)} + \llparenthesis \widehat\xi^+_1 - \widehat\xi^+_2\rrparenthesis(\omega)\Big\}.
\end{split} \end{equation*}
So choosing $0<T(\omega) = T\Big(\sum_{i\in\{1,2\}}\Big\{\Vert {\bf v}_i\Vert_{\bbL^p(\Omega \,; \mcE^\beta_T)}  + \llparenthesis \widehat\xi_i^+\rrparenthesis(\omega) \Big\}\Big)\leq T_1(\omega)$ small enough ensures that the map $\Phi_{\widehat\xi^+,\mu}$ has a unique fixed point ${\bf u}_{\widehat\xi^+,\mu}(\omega)$, which further satisfies
\begin{equation} \label{EqLocalLipschitzDependence}
d_{\mcE^\beta_T}\Big({\bf u}_{\widehat\xi^+_1, {\bf v}_1}(\omega) , {\bf u}_{\widehat\xi^+_2, {\bf v}_2}(\omega)\Big) \lesssim_{M'(\omega)}  \Vert {\bf v}_1 - {\bf v}_2\Vert_{\bbL^p(\Omega \,; \mcE^\beta_T)}  + \llparenthesis \widehat\xi^+_1 - \widehat\xi^+_2\rrparenthesis(\omega).
\end{equation}

\ssk

{\it Step 2: Long time well-posedness.} { \color{blue} We follow the analysis of \cite{ShenZhuZhu} as in the proof of Proposition \ref{SolPAM1} in Section \ref{SectionRegularInteraction}. 

\ssk

\begin{thm}\label{thm_nonexplosionA} One has a polynomial bound
\begin{equation} \label{EqNonExplosionMeanField}
\big\Vert {\bfu}_{\widehat\xi^+ , {\bfv}}\big\Vert_{\mcE^\beta_T} \lesssim \big(1 + \llparenthesis \widehat\xi^+ \rrparenthesis \big)^{\kappa_1} \big(1+\overline{\bbE}\big[\norme{\bfv}^2_{\mcE^\beta_T}\big]\big)^{\kappa_2}(1+\norme{u_0}_{C^\alpha})^{\kappa_3}
\end{equation}
for some positive finite constants $\kappa_1,\kappa_2,\kappa_3$.
\end{thm}

\ssk

The proof of this statement is given in Appendix \ref{section_appendix_nonexplosion}.  }

\ssk

{\it Step 3: Renormalized equation.} Recall that $(\xi,X\odot\xi)\in\boldsymbol{\mcN}_{T_0}$ is the limit in any $\bbL^q(\Omega)$ space, $1\leq q<\infty$, of the sequence of enhanced noises $\big(\xi^{\epsilon_k} , \xi^{\epsilon_k}\odot X^{\epsilon_k} - c^{\epsilon_k}\big)$, for some appropriate diverging functions $c^{\epsilon_k}$, and that $\xi\odot\overline X$ is the limit in $\bbL^q(\Omega^2,\bbP^{\otimes 2})$ of $\xi^{\epsilon_k}\odot \overline{X}^{\epsilon_k} $. We then have 
\begin{align*}
f_o({\bf u}^{\epsilon_k} \vert {\bf v}_t)\xi^{\epsilon_k}_t + g_o(u^{\epsilon_k}  v_t) &= f(u^{\epsilon_k} \vert v_t) \varolessthan \xi^{\epsilon_k}_t + \xi^{\epsilon_k}_t \varolessthan f_o(u^{\epsilon_k} \vert v_t) + f_o({\bf u}^{\epsilon_k} \vert v_t)^\sharp\odot \xi^{\epsilon_k}_t 
\\
&\quad+ {\sf C} \big(\delta_z f_o(u^{\epsilon_k} \vert v_t), X^{\epsilon_k} ,\xi^{\epsilon_k}_t) + \overline{\bbE}\Big[ {\sf C} \big(\delta_\mu f_o(u^{\epsilon_k} \vert v_t), \overline{X}^{\epsilon_k} ,\xi^{\epsilon_k}_t \big)\Big]
\\
&\quad+ \delta_z f_o(u^{\epsilon_k} \vert v_t) \big(\xi^{\epsilon_k}\odot X^{\epsilon_k} - c^{\epsilon_k}\big) + \overline{\bbE}\Big[ \delta_\mu f_o(u^{\epsilon_k} \vert v_t) \big(\xi^{\epsilon_k}_t \odot \overline{X}^{\epsilon_k}  \big) \Big]
\\
&\quad+ g_o(u^{\epsilon_k} \vert v_t)
\\
&= f_o(u^{\epsilon_k} \vert v_t) \xi^{\epsilon_k}_t - c^{\epsilon_k} \big(f_o'f_o\big)(u^{\epsilon_k} \vert v_t) + g_o(u^{\epsilon_k} \vert v_t),
\end{align*}
so the function $u^{\epsilon_k}$ is a solution of the renormalized equation
\[
(\partial_t-\Delta)u^{\epsilon_k} = f_o(u^{\epsilon_k} \vert v_t) \xi^{\epsilon_k}_t - c^{\epsilon_k} \big(f_o'f_o\big)(u^{\epsilon_k} \vert v_t) + g_o(u^{\epsilon_k} \vert 
\]
As we know from \eqref{EqLocalLipschitzDependence} and Step 2 that the solution ${\bf u}_{\widehat\xi^+,{\bf v}}\in \mcE^\beta_{T_0}(X)$ is a continuous function of $\widehat\xi^+\in\boldsymbol{\mcN}_{T_0}\times\boldsymbol{\mcN}_{T_0}$, and since $\widehat\xi^{\epsilon_k+}$ converges to $\widehat\xi^+$ in probability, we see that ${\bf u}_{\widehat\xi^+,{\bf v}}$ is the limit in probability of the sequence of elements of $\mcE^\beta_{T_0}(X^{\epsilon_k})$ associated with $u^{\epsilon_k}$.
\end{Dem}

\ssk

The proof of well-posedness of Equation \eqref{EqMainEq} requires a second fixed point which is the object of the next statement. We fix as above $4\leq p<\infty$.

\ssk

\begin{thm} \label{SolMcKeanSingular}
Assume that $g$ satisfies the Lipschitz assumption \eqref{EqConditionLipG}. There exists a deterministic positive time $T\leq T_0$ with the following property.
\begin{itemize}
	\item[--] For every $u_0\in C^\alpha$ there exists a unique solution ${\bf u}=(u',u^\sharp)$ to \eqref{EqMainEq} in $\bbL^p(\Omega \,; \mcE^\beta_T (X))$. The law $\mcL({\bf u}) \in \mcP_p(\mcE^\beta_T (X))$ of $\bf u$ depends continuously on $\widehat\xi^+$ {\color{black} in some fixed ball of} $\bbL^{r_1}(\Omega^2 \,; {{\boldsymbol{\mcN}}_T}\times{{\boldsymbol{\mcN}}_T})$  {\color{black} for some appropriate choice of exponent $r_1\in(0,\infty)$}.   \vspace{0.15cm}
	
	\item[--] The function $u\in\scrC^\alpha_T$ associated with $\bf u$ is the limit in probability of the family of solutions of the renormalized equations
\[
(\partial_t - \Delta)u^{\epsilon_k} = f\big(u^{\epsilon_k} , \mcL(u^{\epsilon_k}_t)\big) \xi^{\epsilon_k}_t - c^{\epsilon_k}_t \big(f'f\big)(u^{\epsilon_k} , \mcL(u^{\epsilon_k}_t)) + g(u^{\epsilon_k},\mcL(u^{\epsilon_k}_t)).
\]
\end{itemize}
\end{thm}

\medskip

{\color{black} The proof makes clear the dependence of the time horizon $T$ on the size of $\widehat\xi^+\in\bbL^{r_1}(\Omega^2 \,; {{\boldsymbol{\mcN}}_T}\times{{\boldsymbol{\mcN}}_T})$.}

\medskip

\begin{Dem}
{\color{black} We proceed in two steps.

-- First we check that the map $\Psi : {\bf v}\mapsto {\bf u}_{\widehat\xi^+,{\bf v}}$ sends $\bbL^p(\Omega\,;\,\mcE^\beta_T(X))$ into itself for an appropriate choice of $T=T\big(\overline{\bbE}[\Vert {\bf v}\Vert^p_{\mcE^\beta_{{T}}}]\big)$.} One has {\color{black} from the fact that ${\bf u}_{\widehat\xi^+,{\bf v}}$ solves the equation \eqref{EqPAM} and from Proposition \ref{schauder3}
\begin{align*}
\big\Vert {\bf u}_{\widehat\xi^+,{\bf v}} \big\Vert_{\mcE^\beta_T} &\lesssim \norme{u_0}_{C^\alpha} + T^{\frac{\alpha-\beta}{2}} \big(1 + \llparenthesis\widehat\xi^+\rrparenthesis(\omega)^4 \big) \Big(1+ \norme{{\bf u}_{\widehat\xi^+,{\bf v}}}^2_{\mcE^\beta_T} + \overline{\bbE}\big[ \norme{{\bf v}}^4_{\mcE^\beta_T}\big]^{1/2}\Big)
\end{align*}
Using Theorem \ref{thm_nonexplosionA}, we get the bound 

\begin{align*}
\norme{\bfu_{\widehat\xi^+,\bfv}}_{\mcE^\beta_T}&\lesssim 1 + T^{\frac{\alpha-\beta}{2}} \Big(1 + \llparenthesis\widehat\xi^+\rrparenthesis(\omega)^{4+\frac{3}{2}\kappa_1} \Big) \Big(1 + \overline{\bbE}\big[ \norme{{\bf v}}^4_{\mcE^\beta_T}\big]^{\frac{3}{2}\kappa_2} \Big) \norme{{\bf u}_{\widehat\xi^+,{\bf v}}}^{\frac{1}{2}}_{\mcE^\beta_T}.
\end{align*}
Integrating and using Hölder inequality we get for some appropriate exponents $r_1,r_2$
\begin{align*}
\bbE\Big[ &\big\Vert {\bf u}_{\widehat\xi^+,{\bf v}} \big\Vert_{\mcE_T^\beta}^p\Big]^2 \lesssim 1 +  T^{p(\alpha-\beta)} \hspace{-0.03cm} \Big(1 \hspace{-0.03cm}+ \bbE\big[\llparenthesis\,\widehat\xi^+\rrparenthesis^{r_1}\big]\Big) \hspace{-0.03cm} \hspace{-0.03cm}  \Big(1 \hspace{-0.03cm}+\hspace{-0.03cm} \overline{\bbE}\big[\norme{{\bf v}}^{4}_{\dab}\big]^{r_2}\Big) \bbE\big[\norme{{\bf u}_{\widehat\xi^+,\bfv}}_{\mcE_T^{\beta}}^{p}\big].
\end{align*}}
So for $T = T\big(\overline{\bbE}\big[\Vert{\bf v}\Vert_{\dab}^4\big]\big)$ sufficiently small we have 
\[
\bbE\big[\norme{{\bf u}}_{\mcE_T^{\alpha,\beta}}^{p}\big]^{\frac{1}{p}} \lesssim 1 + T^{\frac{\alpha-\beta}{2}}\Big(1 + \bbE\big[\llparenthesis\widehat\xi^+\rrparenthesis^{r_1}\big]^{\frac{1}{2p}} \Big)\, \bbE\big[\norme{{\bf v}}_{\mcE_T^\beta}^{p}\big]^{\frac{r_2}{p}}.
\]
Pick 
\[
A > \max\big(C_0^2 , 2\bbE\big[\llparenthesis\widehat\xi^+\rrparenthesis^{16p}\big]\big).
\]
For $M$ sufficiently big and $T=T(M,A)$ even smaller the map $\Psi_{\widehat\xi^+}$ sends the ball
\[
\left\{ {\bf v} \in \bbL^p\big(\Omega \,; \dab(X)\big) \,;\, \norme{{\bf v}}_{\bbL^p(\Omega \,; \dab)}\leq M\right\}
\]
into itself. Now pick ${\bf v}_1, {\bf v}_2$ in $\bbL^p(\Omega \,; \dab(X))$ and $\widehat\xi^+_1, \widehat\xi^+_2$ in $\bbL^{r_1}(\Omega^2 \,; {{\boldsymbol{\mcN}}_T}\times{{\boldsymbol{\mcN}}_T})$ such that
\[
\bbE\big[\llparenthesis\widehat\xi^+_i\rrparenthesis^{r_1}\big] \leq A, \qquad \overline{\bbE}\big[\norme{{\bf v}_i}_{\dab}^p\big] \leq M,
\]  
for $1\leq i\leq 2$. Write ${\bf u}_i$ for $\Phi_{\widehat\xi^+_i}({\bf v}_i)$ and define the random variable
\[
R(\omega) \defeq \llparenthesis \widehat\xi^+_1\rrparenthesis(\omega) + \llparenthesis \widehat\xi^+_2\rrparenthesis(\omega).
\] 
We have
\begin{equation*} \begin{split}
d_{\mcE_T^\beta}\big({\bf u}_1, {\bf u}_2\big) &\lesssim_{R(\omega)} T^\delta \bigg\{ \llparenthesis \widehat\xi^+_1 - \widehat\xi^+_2\rrparenthesis(\omega) + d_{\mcE^\beta_T}({\bf u}_1, {\bf u}_2) + \overline{\bbE}\big[d_{\dab}{\big({\bf v}^1, {\bf v}^2\big)^4}\big]^{\frac{1}{4}} \bigg\}
\\ 
& \lesssim_{R(\omega)}  T^\delta \bigg\{ \llparenthesis \widehat\xi^+_1 - \widehat\xi^+_2\rrparenthesis(\omega) + d_{\mcE^\beta_T}({\bf u}_1, {\bf u}_2)^{\frac{1}{2}}\hspace{-2pt} + \overline{\bbE}\big[d_{\dab}{\big({\bf v}^1,{\bf v}^2\big)^4}\big]^{\frac{1}{4}} \bigg\},   
\end{split} \end{equation*}
for some implicit positive multiplicative constant that is a polynomial of $R(\omega)$, combining Propositions \ref{schauder3}, \ref{PropNonlinearity} and \ref{PropProduct}. Integrating and using Cauchy-Schwarz inequality we obtain the estimate
\begin{align*}
 \bbE\Big[d_{\mcE_T^\beta}\big({\bf u}_1, {\bf u}_2\big)^p \Big]^2 \lesssim \bbE\big[ \llparenthesis \widehat\xi^+_1 - \widehat\xi^+_2\rrparenthesis^{2p}\big] + T^{2p\delta} \bigg\{ \bbE\big[ d_{\mcE_T^\beta}\big({\bf u}_1, {\bf u}_2\big)^{p} \big] + \overline{\bbE}\big[d_{\dab}\big({\bf v}_1, {\bf v}_2\big)^4\big]^{\frac{p}{2}} \bigg\},
\end{align*}
so taking $T>0$ deterministic small enough ensures that we have
\begin{align*}
\bbE\big[ d_{\mcE_T^\beta}\big({\bf u}_1, {\bf u}_2\big)^{p} \big]^2 \lesssim \bbE\big[ \llparenthesis \widehat\xi^+_1 - \widehat\xi^+_2\rrparenthesis^{2p}\big] + T^{2p\delta} \, \overline{\bbE}\big[d_{\dab}\big({\bf v}_1, {\bf v}_2\big)^4\big]^{\frac{p}{2}}.
\end{align*}
As $4\leq p<\infty$, we conclude that Equation \eqref{EqMainEq} has a unique local solution $\bf u$ in $\mcP_p(\dab(X))$, and that the law $\mcL({\bf u}) \in \mcP_p(\dab(X))$ of $\bf u$ depends continuously on $\widehat\xi^+\in \bbL^{r_1}(\Omega^2 \,; {{\boldsymbol{\mcN}}_T}\times{{\boldsymbol{\mcN}}_T})$.

{\color{black} The second statement in Theorem \ref{SolMcKeanSingular} follows from Step 3 in the proof of Proposition \ref{PropSolPAM}.}
\end{Dem}

\section{Propagation of chaos}
\label{SectionChaos}

Let now $(\xi^i,u_0^i)$ be a sequence of independent and identically distributed random variables with common law $\mcL(\xi,u_0)$, defined on some probability space $(\Omega,\mcF,\bbP)$. We fix $\omega\in\Omega$ and an integer $n\geq 1$ and study the dynamics
\begin{equation}\label{EqSystemSingular} \begin{split}
(\partial_t-\Delta)u^{i,n}(\omega) &= f\big(u^{i,n}(\omega) , \mu_t^n\big)\xi^i(\omega) +g\big(u^{i,n}(\omega) , \mu_t^n(\omega)\big), \qquad  (1\leq i\leq n)   \\
\mu_t^n(\omega) &\defeq \frac{1}{n}\sum_{i=1}^n\delta_{u^{i,n}_t(\omega)},
\end{split} \end{equation}
with initial conditions $\big(u_0^1(\omega),\dots,u_0^n(\omega)\big)$. We suppose here that $f$ is of the form \eqref{InteractionPolynomial} with $F : \bbR^{m+1}\rightarrow\bbR$ of class $C^3_b$, and that $g$ satisfies the Lipschitz condition \eqref{EqConditionLipG}. 

System \eqref{EqSystemSingular} can either be understood as a multidimensional singular stochastic PDE driven by a multidimensional (enhanced) noise or as a mean field singular stochastic PDE. 

\ssk

\begin{thm} \label{ThmInterpretation}
The multi-dimensional and mean field interpretations of the system \eqref{EqSystemSingular} coincide.
\end{thm}

\ssk

\begin{Dem}
To lighten the notations we consider here the case that the  diffusivity $f$ is linear in its measure argument -- see \eqref{EqLinearSimpleCase} below. The polynomial case is treated similarly. One can see Equation \eqref{EqSystemSingular} as a single multidimensional singular stochastic equation 
\[
(\partial_t-\Delta){\sf u}^{(n)} = {\sf f}({\sf u})\xi^{(n)} + {\sf g}({\sf u}^{(n)})
\]
with multi-dimensional unknown ${\sf u}^{(n)} = \big(u^{1,n},\dots,u^{n,n}\big)$ and noise $\xi^{(n)} = \big(\xi^1,\dots,\xi^n\big)$, and where ${\sf f}$ is $(f^1,\dots,f^n)$ with
\[
f^i : \big(u^{1,n},\dots,u^{n,n}\big)\mapsto f\bigg(u^{i,n},\frac{1}{n}\sum_{j=1}^n\delta_{u^{j,n}}\bigg) \eqdef f(u^{i,n},\mu^n),
\]
with a similar definition of $\sf g$. The noise $\xi^{(n)}$ needs to be enhanced to make sense of the equation. The solution will be a tuple of paracontrolled functions 
\[
u^{i,n} = (u^{i,n})' \varolessthan X^i + (u^{i,n})^\sharp = f^i(u^{1,n},\dots,u^{n,n}) \varolessthan X^i + (u^{i,n})^\sharp
\]
so we will have from paralinearisation
\[
f^i\big(u^{1,n},\cdots,u^{n,n} \big) = \sum_{j=1}^n \Big(\partial_j f^i\big(u^{1,n},\dots,u^{n,n}\big)(u^{j,n})'\Big) \varolessthan X^j + f^i\big(u^{1,n},\dots,u^{n,n}\big)^\sharp, 
\]
with
\[
\partial_j f^i\big(u^{1,n},\dots,u^{n,n}\big)=\delta_{i,j}\partial_1 f\big(u^{i,n},\mu^n\big) + \frac{1}{n} \partial_2 F\big(u^{i,n},\mu^n\big),
\]
since 
\begin{equation} \label{EqLinearSimpleCase}
f(u^{i,n}, \mu^n) = \frac{1}{n}\sum_{j=1}^n F\big(u^{i,n},u^{j,n}\big).
\end{equation}
The singular product in \eqref{EqSystemSingular} then reads

\begin{equation} \label{EqFormulation1} \begin{split}
f\big(u^{i,n},\mu^n\big)\xi^i &= f\big(u^{i,n},\mu^n\big) \varolessthan \xi^i + \xi^i\varolessthan f\big(u^{i,n},\mu^n\big)+ f\big(u^{i,n},\mu^n\big)^\sharp\odot\xi^i 
\\
&\quad+{\sf C}\Big( \partial_1 f\big(u^{i,n},\mu^n\big)(u^{i,n})' , X^i , \xi^i\Big) + \frac{1}{n}\sum_{j=1}^n {\sf C}\Big(\partial_2 F(u^i,\mu^n)(u^{j,n})' , X^j , \xi^i \Big)
\\
&\quad +  \partial_1 f\big(u^{i,n},\mu\big)(u^{i,n})' \big(\xi^i\odot X^j\big) + \frac{1}{n}\sum_{j=1}^n \partial_2 F(u^i,\mu)(u^{j,n})' \big(\xi^i \odot X^j\big) .
\end{split} \end{equation}
Our task is now to prove that \eqref{EqSystemSingular} may also be understood as a mean field singular stochastic PDE with a suitable enhancement of the noise and that the two interpretations coincide. With the notations of Section \ref{SubsectionAdditiveChaos}, Tanaka's trick gives an interpretation of \eqref{EqSystemSingular} as the mean field type equation
\begin{equation}\label{MeanFieldTanakaSingular}
(\partial_t-\Delta)u^{i,n}(\omega) = f_\circ\Big(u^{i,n}(\omega) \big\vert u^{U_n(\cdot),n}(\omega) \Big)\xi^i(\omega) + g_\circ\Big(u^{i,n}(\omega) \big\vert u^{U_n(\cdot),n}(\omega) \Big)
\end{equation}
studied in Section \ref{SectionGeneral}, but now set on the finite probability space $\big([\![ 1,n]\!],2^{[\![ 1,n]\!]},\lambda_n\big)$, with generic chance element $i$. The enhanced noise from Definition \ref{MeanFieldEnhancement} is then 
\[
\Big\{ \xi^i,\enskip \xi^i\odot X^i,\enskip \xi^j, \enskip \xi^j \odot X^i \Big\}_{1\leq i,j\leq n}, 
\]
where the index $i$ plays the role of $\omega$ and $j$ the role of $\varpi$. Let us now clarify the meaning of the singular product. We have
\[
\delta_z f_\circ\big(u^{i,n} \vert u^{u^{\epsilon_k}(\cdot)}\big) = \partial_1 f\big(u^{i,n} \vert u^{U_n(\cdot),n} \big)\big(u^{i,n}\big)' ,
\]
and 
\[
\delta_\mu f_\circ\big(u^{i,n} \vert u^{U_n(\cdot),n}\big) = \partial_2 F\big(u^{i,n} , v^{U_n(\cdot),n}\big)\big(u^{U_n(\cdot),n}\big)'.
\]
In the sense of Section \ref{SubsectionParacStructure} the singular product in Equation \eqref{MeanFieldTanakaSingular} is defined as
\begin{equation} \label{EqFormulation2} \begin{split}
f_\circ\big(u^{i,n} \vert u^{U_n(\cdot),n} \big)\xi^i &= f\big(u^{i,n} \vert u^{U_n(\cdot),n} \big)  \varolessthan \xi^i + \xi^i \varolessthan f_\circ\big(u^{i,n} \vert u^{U_n(\cdot),n} \big) + f\circ\big(u^{i,n} \vert u^{U_n(\cdot),n} \big)^\sharp\odot\xi^i   \\
&\quad+{\sf{C}}\Big( \partial_1 f_\circ\big(u^{i,n} \vert u^{U_n(\cdot),n} \big)\big(u^{i,n}\big)',X^i,\xi^i \Big) + \partial_1 f_\circ\big(u^{i,n} \vert u^{U_n(\cdot),n} \big)\big(u^{i,n}\big)'\big(\xi\odot X\big)^i   \\
&\quad+\frac{1}{n}\sum_{j=1}^n {\sf C}\Big( \partial_2 F\big(u^{i,n},u^{j,n}\big)\big(u^{j,n}\big)', X^j,\xi^i  \Big)   \\
&\quad+\frac{1}{n}\sum_{j=1}^n \partial_2 F\big(u^{i,n},u^{U_n(\cdot),n}\big)\big(u^{U_n(\cdot),n}\big)'\big(\xi^i\odot X^j\big).
\end{split} \end{equation}
We conclude from \eqref{EqFormulation1} and \eqref{EqFormulation2} that the two formulations coincide as they amount to solving the same classical PDE for the remainders $(u^{i,n})^\sharp$.
\end{Dem}

\ssk

We know from the continuity result of Theorem \ref{SolMcKeanSingular} that the almost sure convergence of 
\[
\mathcal{W}_p \bigg(\frac{1}{n}\sum_{i=1}^n \delta_{(\widehat\xi^{i,+},u_0^i)(\omega)} , \mcL\big( \widehat\xi^+, u_0\big) \bigg)
\]
to $0$ granted by the law of large numbers entails the convergence
\[
\mcW_{p,C_TC^\alpha}\bigg(\frac{1}{n}\sum_{i=1}^n \delta_{u^{i,n}} , \mcL(u)\bigg) \underset{n\rightarrow+\infty}{\longrightarrow} 0
\]
where $u$ is the function associated with the solution $\bf u$ of the mean field dynamics \eqref{EqMainEq}. It follows then from Sznitman's Proposition 2.2 in \cite{Snitzman89} that there is propagation of chaos for the system \eqref{EqSystemSingular} of interacting fields to the mean field limit dynamics \eqref{EqMainEq}.

\ssk

\begin{thm} \label{ThmPropagationChaos}
For any fixed integer $k$ the law of $\big(u^{1,n},\dots, u^{k,n}\big)$ converges to $\scrL(u)^{\otimes k}$ when $n $ tends to $+\infty$.
\end{thm}

\ssk

{\color{black} The setting of Section \ref{SectionGeneral} is more general than the setting of Section \ref{SectionRegularInteraction}. Theorem \ref{ThmPropagationChaos} thus holds for the two settings at once.}

\medskip

\appendix
\section{Enhancing some random noises}
\label{SectionEnhancingNoises}

We prove Theorem \ref{ThmNoiseEnhancement} in this section.  Recall from \eqref{EqDefnXOdotXi} the definition of the random variable $X\odot\xi$. Write $e_k$ for the function $x\mapsto \exp (i(k,x))$ and $\widehat\xi(k)$ for $(\xi,e_k)$. Our noises satisfy the identity
\begin{equation}\label{observation}
\bbE\big[\,\widehat\xi_t(k)\widehat\xi_s(-k')\big] = \mathbf{1}_{k=k'} \, c(t,s) \,\widehat\eta(k).
\end{equation}
We denote below by $\textsc{Var}(A)$ the variance of a random variable $A$.

\ssk

\begin{lem}\label{LemE}
There exists a positive constant $\kappa$ such that on has for all $\ell\in \bbN$, $s,t,a,b\in \bbR_+$ and $x\in {\sf T}^2$, the estimate
\[
\textsc{\emph{Var}}\Big(\Delta_\ell \big(P_t(\xi_s)\odot \xi_a\big)(x)\Big) \lesssim \frac{2^{2\ell}2^{2\ell\eta}}{t} \, e^{- \kappa t2^{2\ell}}\big( c(s,s)\,c(a,a) + c(s,a)^2 \big)
\]
and 
\[
\textsc{\emph{Var}}\Big(\Delta_\ell \Big(\big((\emph{Id}-P_b)P_t(\xi_s)\big)\odot \xi_a\Big)(x)\Big) \lesssim b\,\frac{2^{2\ell}2^{2\ell\eta}}{t}\,e^{-\kappa t2^{2\ell-1}}\big( c(s,s)\,c(a,a) + c(s,a)^2 \big).
\]
\end{lem}

\ssk

\begin{Dem}
The proof follows closely the proof of Lemma 5.2 in \cite{GIP}. We have
\begin{align*}
\Delta_\ell\big(P_t(\xi_s)\odot &\xi_a\big)(x)   \\
&= (2\pi)^{-2}\sum_{k\in\bfZ^2}e^{i(k,x)}\rho_\ell(k) \mathcal{F}\big(P_t(\xi_s)\odot\xi_a\big)(k)   \\
&=(2\pi)^{-4} \hspace{-0.2cm} \sum_{k_1,k_2\in\bfZ^2}\sum_{|i-j|\leq 1} \hspace{-0.1cm} \rho_\ell(k_1+k_2)\rho_i(k_1) \, e^{-t|k_1|^2} \, \widehat\xi_s(k_1) \,  \rho_j(k_2) \, \widehat\xi_a(k_2) \, e_{k_1+k_2}(x),
\end{align*}
then $\textsc{{Var}}\Big(\Delta_\ell(P_t(\xi_s)\odot\xi_a)(x)\Big)$ is equal to

\begin{align*}
&(2\pi)^{-8}\sum_{k_1,k_2,k_1',k_2'}\sum_{|i-j|\leq 1}\sum_{|i'-j'|\leq 1}\rho_\ell(k_1+k_2) \, \rho_i(k_1) \, e^{-t|k_1|^2}\rho_j(k_2)
\\
&\quad\times \rho_\ell(k_1'+k_2') \, \rho_{i'}(k_1') \, e^{-t|k_1'|^2} \, \rho_{j'}(k_2') \, \textsc{Cov}\Big(\widehat\xi_s(k_1)\widehat\xi_a(k_2) \,,\, \widehat\xi_s(k_1')\widehat\xi_a(k_2')\Big)e_{k_1+k_2+k_1'+k_2'}(x).
\end{align*}
Using Wick theorem and the identity \ref{observation} one gets

\begin{align*}
\text{Cov}\Big(\widehat\xi_s(k_1) &\widehat\xi_a(k_2) \,,\, \widehat\xi_s(k_1')\widehat\xi_a(k_2')\Big)   \\
&= \bbE\Big[\widehat\xi_s(k_1)\widehat\xi_a(k_2)\widehat\xi_s(k_1')\widehat\xi_a(k_2')\Big]-\bbE\Big[\widehat\xi_s(k_1)\widehat\xi_a(k_2)\Big]\bbE\Big[\widehat\xi_s(k_1')\widehat\xi_a(k_2')\Big]   \\ &= \bbE\Big[\widehat\xi_s(k_1)\widehat\xi_s(k_1')\Big]\bbE\Big[\widehat\xi_a(k_2)\widehat\xi_a(k_2')\Big] + \bbE\Big[\widehat\xi_s(k_1)\widehat\xi_a(k_2')\Big]\bbE\Big[\widehat\xi_s(k_1')\widehat\xi_a(k_2)\Big]    \\
 &= (2\pi)^4 \, \widehat\eta(k_1) \, \widehat\eta(k_2)\Big( \mathbf{1}_{k_1=-k_1', k_2=-k_2'}\,c(s,s) \, c(a,a) + \mathbf{1}_{k_1=-k_2', k_2=-k_1'}c(s,a)^2\Big),
\end{align*}
consequently
\begin{align*}
\textsc{{Var}}\Big(&\Delta_\ell(P_t(\xi_s)\circ\xi_a)(x)\Big)   \\
&= \sum_{k_1,k_2} \sum_{|i-j|\leq 1}\sum_{|i'-j'|\leq 1} (2\pi)^4 \, \widehat\eta(k_1) \, \widehat\eta(k_2) \, \rho_\ell(k_1+k_2)^2 \, \rho_i(k_1) \, \rho_j(k_2)
\\
&\hspace{1.5cm}\times\Big( c(s,s) \, c(a,a) \, \rho_{i'}(k_1) \, \rho_{j'}(k_2)e^{-2t|k_1|^2} + c(s,a)^2\rho_{i'}(k_2)\rho_{j'}(k_1) \, e^{-t|k_1|^2-t|k_2|^2} \Big).
\end{align*}
The factors $ \rho_i(k_1)\,\rho_{i'}(k_1)$ and $\rho_i(k_1)\,\rho_{j'}(k_1) $ ensure that one can restrict the sum on $i$ and $i'$ to couples $(i,i')$ such that $\frac{1}{\mu}|i|\leq |i'| \leq \mu |i|$ for some constant $\mu$, which will be denoted by $i\sim i'$. Likewise the factor $\rho_\ell(k_1+k_2)$ enables us to restrict the sum to $|i|\geq \frac{1}{\mu'}l$ for some $\mu'$. There exists some $\kappa_0>0$ such that $e^{-2t|k|^2}\lesssim e^{-t\kappa_0 2^{2i}}$ for $k\in\text{supp}(\rho_i)$, so that for some $\kappa>0$ 
\begin{align*}
\textsc{{Var}}\Big(\Delta_\ell&(P_t(\xi_s)\circ\xi_a)(x)\Big)
\\
&\lesssim \big(c(s,s)\,c(a,a) + c(s,a)^2\big) \sum_{i,i',j,j'} \mathbf{1}_{\ell\lesssim i}\mathbf{1}_{i\sim i' \sim j \sim j'} \sum_{k_1,k_2} \mathbf{1}_{\text{supp}(\rho_\ell)}(k_1+k_2) \\
&\quad \times \mathbf{1}_{\text{supp}(\rho_i)}(k_1)\mathbf{1}_{\text{supp}(\rho_j)}(k_2)2^{2i\eta} e^{-2t\kappa 2^{2i}}
\end{align*}

\begin{align*}
&\lesssim \big(c(s,s)\,c(a,a) + c(s,a)^2\big) \sum_{i, l\lesssim i}2^{2i}2^{2\ell}2^{2i\eta}e^{-2t\kappa 2^{2i}} \
\\
&\lesssim \big(c(s,s) \, c(a,a) + c(s,a)^2\big) \frac{2^{2\ell}2^{2\ell\eta}}{t}e^{-2t\kappa 2^{2\ell}},
\end{align*}
hence the first estimate. For the second estimate we notice that the $e^{-t|k_1|^2}$ is replaced by $(1-e^{-b|k_1|^2})\,e^{-t|k_1^2|}$ and that
\[
(1-e^{-b|k_1|^2}) \, e^{-t|k_1^2|}\leq b|k_1|^2 e^{-t|k_1|^2} \lesssim v e^{-t |k_1|^2/2}
\]
The remainder of the proof is the same as for the first estimate.
\end{Dem}

\ssk

We can now prove Theorem \ref{ThmNoiseEnhancement}. We estimate $\bbE\big[ \big\Vert\big(X\odot\xi\big)(t)-\big(X\odot\xi\big)(s)\big\Vert_{B^{2\alpha-2}_{2p,2p}}^{2p}\big]$, and use then Kolmogorov continuity criterion and  Besov embedding. For $0<s\leq t$ we write
\begin{align*}
\int_0^t P_{t-a}&(\xi_a)\odot\xi_t \, da - \int_0^s P_{s-a}(\xi_a)\odot \xi_s \, da    \\
&= \int_0^s \hspace{-0.1cm} \big( (P_{t-s} - \textrm{Id})P_{s-a}(\xi_a) \big)\odot\xi_t \,da + \int_0^s \hspace{-0.1cm}P_{s-a}(\xi_a) \odot (\xi_t-\xi_s) \,da + \int_s^t \hspace{-0.1cm}P_{t-a}(\xi_a)\odot\xi_t \,da   \\ 
& =: \int_0^s \widetilde A_1(a)\,da +\int_0^s \widetilde A_2(a)\,da +\int_s^t \widetilde A_3(a)\,da,
\end{align*}
and set 
\[
A_i \defeq \widetilde A_i - \bbE\big[\widetilde A_i\big] \qquad (i\in [\![1,3 ]\!]).
\]
The quantity 
\[
\bbE\Big[\norme{\big(X\odot\xi\big)(t)-\big(X\odot\xi\big)(s)}_{B^{2\alpha-2}_{2p,2p}}^{2p}\Big]
\]
is equal to  
\begin{align*}
\sum_{\ell\geq -1} &2^{2p\ell(2\alpha-2)} \int_{{\sf T}^2}\bbE\Big[ \Big|\Delta_\ell \Big(\big(X\odot\xi\big)(t)-\big(X\odot\xi\big)(s)\Big) \Big|^{2p}\Big]   \\
&= \sum_{\ell\geq -1} 2^{2p\ell(2\alpha-2)}\int_{{\sf T}^2} \bbE\bigg[ \bigg| \int_0^s \Delta_\ell A_1(a)\,da + \int_0^s \Delta_\ell A_2(a)\,da + \int_s^t \Delta_\ell A_3(a)\,da \bigg|^{2p}\bigg].
\end{align*}
From Gaussian hyper-contractivity we have
\begin{align*}
&\bbE\bigg[ \bigg| \int_0^s \Delta_\ell A_1(a)\,da + \int_0^s \Delta_\ell A_2(a)\,da + \int_s^t \Delta_\ell A_3(a)\,da \bigg|^{2p}\bigg]
\\
&\leq\bbE\bigg[  \int_0^s \big|\Delta_\ell A_1(a)\big|\,da + \int_0^s \big|\Delta_\ell A_2(a)\big|\,da + \int_s^t \big|\Delta_\ell A_3(a)\big|\,da \bigg]^{2p}
\\
&\lesssim \bigg(\hspace{-0.03cm}\int_0^s\bbE\big[ \big|\Delta_\ell A_1(a)\big|^2\big]^{1/2}\,da\bigg)^{2p} \hspace{-0.1cm}+\hspace{-0.03cm} \bigg(\hspace{-0.03cm}\int_0^s\bbE\big[ \big|\Delta_\ell A_2(a)\big|^2\big]^{1/2}\,da\bigg)^{2p} \hspace{-0.1cm}+\hspace{-0.03cm} \bigg(\hspace{-0.03cm}\int_s^t\bbE\big[ \big|\Delta_\ell A_3(a)\big|^2 \big]^{1/2}\,da\bigg)^{2p},
\end{align*}
so that the bounds for 
\[
\bbE\bigg[\big\Vert\big(X\odot\xi\big)(t)-\big(X\odot\xi\big)(s)\big\Vert_{B^{2\alpha-2}_{2p,2p}}^{2p}\bigg]^{1/(2p)}
\] 
becomes
\begin{align*}
&\left( \sum_{\ell\geq -1} 2^{2p\ell(2\alpha-2)}\int_{{\sf T}^2} \bbE\bigg[ \bigg| \int_0^s \Delta_\ell A_1(a)\,da + \int_0^s \Delta_\ell A_2(a)\,da + \int_s^t \Delta_\ell A_3(a)\,da \bigg|^{2p}\bigg] \right)^{\frac{1}{2p}}   \\
\end{align*}
\begin{align*}
&\lesssim  \sum_{\ell\geq -1} 2^{\ell(2\alpha-2)}\left(\int_0^s \bbE\big[\big| \Delta_\ell A_1(a)\big|^2 \big]^{\frac{1}{2}} \,da + \int_0^s \bbE\big[\big| \Delta_\ell A_2(a)\big|^2 \big]^{\frac{1}{2}} \,da + \int_s^t \bbE\big[\big| \Delta_\ell A_3(a)\big|^2 \Big]^{\frac{1}{2}} \,da \right)   \\
&\eqdef S_1+S_2+S_3
\end{align*}
We use Lemma \ref{LemE} to estimate the $\bbE\big[\big| \Delta_\ell A_i(a)\big|^2 \big]$. First we have
\begin{align*}
    \bbE\big[\big| \Delta_\ell A_1(a)\big|^2 \big] &= \textsc{{Var}}\Big( \Delta_\ell \Big(\big( (P_{t-s} - \textrm{Id})P_{s-a}(\xi_a) \big)\odot\xi_t\Big)\Big)
    \\
    &\lesssim (t-s)\frac{2^{2\ell}2^{2\ell\eta}}{s-a}\,e^{-\kappa(s-a)2^{2\ell}}\big(c(s,s)\,c(a,a) + c(s,a)^2 \big)
\end{align*}
and
\begin{align*}
    \bbE\big[\big| \Delta_\ell A_2(a)\big|^2 \big] &= \textsc{{Var}}\Big( \Delta_\ell \Big(\big( P_{s-a}(\xi_a) \big)\odot\big(\xi_t-\xi_s)\Big)\Big)
    \\
    &\lesssim \frac{2^{2\ell}2^{2\ell\eta}}{s-a} \, e^{-\kappa(s-a)2^{2\ell}} \Big( c(a,a)\big(c(t,t) + c(s,s) - 2c(s,t)\big) \Big)
\end{align*}
and
\begin{align*}
    \bbE\big[\big| \Delta_\ell A_3(a)\big|^2 \big] &= \textsc{{Var}}\Big( \Delta_\ell \Big(\big( P_{t-a}(\xi_a) \big)\odot\xi_t\Big)\Big)
    \\
    &\lesssim \frac{2^{2\ell}2^{2\ell\eta}}{t-a}e^{-\kappa(t-a)2^{2\ell}}\big( c(t,t) c(a,a) + c(t,a)^2 \big).
\end{align*}
So, writing $c_{st}$ for $c(t,t)+c(s,s)-2c(s,t)$, we get
\begin{align*}
\int_0^s \bbE\Big[\big| \Delta_\ell A_1(a)\big|^2 \Big]^{1/2} \,da &\lesssim  (t-s)^{\frac{1}{2}} 2^\ell 2^{2\ell\eta} \int_0^s e^{-\kappa(s-a)2^{2\ell-1}}\frac{\,da}{(s-a)^{1/2}}
\\
\int_0^s \bbE\Big[\big| \Delta_\ell A_2(a)\big|^2 \Big]^{1/2} \,da &\lesssim  c_{st}^{\frac{1}{2}} \, 2^\ell 2^{2\ell\eta}\int_0^s e^{-\kappa(s-a)2^{2\ell-1}}\frac{\,da}{(s-a)^{1/2}}
\\
\int_s^t \bbE\Big[\big| \Delta_\ell A_3(a)\big|^2 \Big]^{1/2} \,da &\lesssim   2^\ell 2^{2\ell\eta}\int_s^t e^{-\kappa(t-a)2^{2\ell-1}}\frac{\,da}{(t-a)^{1/2}}.
\end{align*}
We have
\begin{align*}
    S_1&\lesssim (t-s)^{\frac{1}{2}}\sum_{\ell\geq -1} 2^{\ell(2\alpha+2\eta-1)} \int_0^s e^{-\kappa(s-a)2^{2\ell-1}}\frac{\,da}{(s-a)^{1/2}}
    \\
    &\lesssim (t-s)^{\frac{1}{2}}\int_0^s \int_{-1}^{+\infty} 2^{x(2\alpha+2\eta-1)}e^{-\kappa(s-a)2^{2x-1}}\frac{dx\,da}{(s-a)^{1/2}}
    \\
    &\lesssim (t-s)^{\frac{1}{2}} \int_0^s \int_0^{+\infty} (s-a)^{-\alpha-\eta} y^{2\alpha+2\eta-2}e^{-\kappa y^2/2} \,dy\,da
\end{align*}
and similarly
\begin{align*}
    S_2&\lesssim c_{st}^{\frac{1}{2}}\sum_{\ell\geq -1} 2^{\ell(2\alpha-1)} \int_0^s e^{-\kappa(s-a)2^{2\ell-1}}\frac{\,da}{(s-a)^{1/2}}
    \\
    &\lesssim c_{st}^{\frac{1}{2}}  \int_0^s \int_0^{+\infty} (s-a)^{-\alpha-\eta} y^{2\alpha+2\eta-2}e^{-\kappa y^2/2} \,dy\,da
\end{align*}
and
\begin{align*}
    S_3&\lesssim \sum_{\ell\geq -1} 2^{\ell(2\alpha-1)} \int_s^t e^{-\kappa(t-a)2^{2\ell-1}}\frac{\,da}{(t-a)^{1/2}}
    \\
    &\lesssim \int_s^t \int_0^{+\infty} (t-a)^{-\alpha-\eta} y^{2\alpha+2\eta-2}e^{-\kappa y^2/2} \dd y\,da.
\end{align*}
Finally we see that
\begin{align*}
\bbE\Big[\big\Vert\big(X\odot\xi\big)(t)-\big(X\odot\xi\big)(s)\big\Vert_{B^{2\alpha-2}_{2p,2p}}^{2p}\Big]&\lesssim \Big( (t-s)^{1/2}+c_{st}^{1/2}+(t-s)^{1-\alpha-\eta}\Big)^{2p}
\\
&\lesssim |t-s|^{2pm},
\end{align*}
with 
\[
m\defeq \min\Big\{ \frac{1}{2}, \, \frac{\delta}{2},\, 1-\alpha-\eta\Big\}.
\] 
From Kolmogorov continuity criterion and Besov embedding, for every $\alpha<1$ and $1\leq p<\infty$, the process $X\odot\xi$ is almost surely an element of $C^{m-1/p}_TC^{2\alpha-2 - 1/p}$.

\ssk

The mollifier approximation result in the statement of Theorem \ref{ThmNoiseEnhancement} comes from the same arguments and calculations writing 
\begin{align*}
(X\odot\xi)(t)-\Big((X^\eps\odot\xi^\eps)(t) &- \bbE\big[(X^\eps\odot\xi^\eps)(t) \big]\Big) 
\\
&= \int_0^t \Big(P_{t-a}(\xi_a - \xi_a^{\eps})\odot \xi_t  - \bbE\big[P_{t-a}(\xi_a - \xi_a^{\eps})\odot \xi_t\big]\Big) da
\\
&\qquad+ \int_0^t \Big(P_{t-a}(\xi_a^\eps)\odot (\xi_t-\xi_t^\eps)  - \bbE\big[P_{t-a}(\xi_a^\eps) \odot (\xi_t-\xi_t^\eps)\big]\Big) da.
\end{align*}
If $\varphi$ is the fourier transform of the mollifier, we have
\[
\widehat\xi^\eps(k)=\varphi(k\eps) \, \widehat\xi(
\]
and the same calculations as in the proof of Lemma \ref{LemE} give
\begin{align*}
&\textsc{{Var}}\Big(\Delta_\ell(P_{t-a}(\xi_a - \xi_a^\eps)\odot\xi_t)(x)\Big)
\\
&\lesssim \sum_{i,i',j,j'} \mathbf{1}_{l\lesssim i}\mathbf{1}_{i\sim i' \sim j \sim j'} \sum_{k_1,k_2} (1-\varphi(k_1\eps))\mathbf{1}_{\text{supp}(\rho_\ell)}(k_1+k_2) \mathbf{1}_{\text{supp}(\rho_i)}(k_1)\mathbf{1}_{\text{supp}(\rho_j)}(k_2)2^{2 i \eta}e^{-2t\kappa 2^{2i}}
\\
&\lesssim
\sum_{i} \mathbf{1}_{\ell\lesssim i} \sum_{k_1,k_2} (1-\varphi(k_1\eps))\mathbf{1}_{\text{supp}(\rho_\ell)}(k_1+k_2) \mathbf{1}_{\text{supp}(\rho_i)}(k_1)\mathbf{1}_{\text{supp}(\rho_j)}(k_2)2^{2 i \eta}e^{-2t\kappa 2^{2i}}.
\end{align*}
For some integer $N=N(\eps)$, one can decompose the last sum as 
\begin{align*}
&\sum_{i\leq N} \mathbf{1}_{\ell\lesssim i} \sum_{k_1,k_2} (1-\varphi(k_1\eps))\mathbf{1}_{\text{supp}(\rho_\ell)}(k_1+k_2) \mathbf{1}_{\text{supp}(\rho_i)}(k_1)\mathbf{1}_{\text{supp}(\rho_j)}(k_2)2^{2 i \eta}e^{-2t\kappa 2^{2i}} 
\\
&+\sum_{i> N} \mathbf{1}_{\ell\lesssim i} \sum_{k_1,k_2} (1-\varphi(k_1\eps))\mathbf{1}_{\text{supp}(\rho_\ell)}(k_1+k_2) \mathbf{1}_{\text{supp}(\rho_i)}(k_1)\mathbf{1}_{\text{supp}(\rho_j)}(k_2)2^{2 i \eta}e^{-2t\kappa 2^{2i}}
\\ 
&\lesssim \sup_{|x|\leq N}\big(1-\varphi(x\eps)\big) \frac{2^{2\ell}2^{2\ell\eta}}{t-a} \, e^{-\kappa(t-a)2^{2\ell}} + \frac{2^{2N}2^{2N\eta}}{t-a} \, e^{-\kappa(t-a)2^{2N}}
\end{align*}
So, choosing $N(\eps)$ such that $N(\eps)\to \infty$ and $\eps N(\eps)\to 0$ as $\eps$ goes to zero, one gets
\[
\textsc{{Var}}\Big( \Delta_\ell\Big( P_{t-a}(\xi_a-\xi_a^{\eps})\odot \xi_t  (x)  \Big)  \Big) \lesssim \psi_\ell(\eps) \, \frac{2^{2\ell(1+\eta)}}{t-a} \, e^{-\kappa(t-a)2^{2\ell}},
\]
where $0\leq \psi_\ell(\eps)\leq 1$ tends to $0$ as $\epsilon>0$ goes to $0$. Likewise one has
\[
\textsc{{Var}}\Big( \Delta_\ell\Big( P_{t-a}(\xi_a^\eps)\odot (\xi_t-\xi_t^\eps)  (x)  \Big)  \Big) \lesssim \psi_\ell(\eps) \, \frac{2^{2\ell(1+\eta)}}{t-a} \, e^{-\kappa(t-a)2^{2\ell}}.
\]
The same calculations as above give for 
\[
\bbE\Big[\big\Vert (X\odot\xi)(t) - \big\{ (X^\eps\odot\xi^\eps)(t)-\bbE\big[(X^\eps\odot\xi^\eps)(t) \big]\big\}\big\Vert_{B^{2\alpha-2}_{2p,2p}}^{2p}\Big] 
\]
the bound
\[
\sum_{\ell\geq 0} \psi_\ell(\eps) \, 2^{\ell(2\alpha+2\eta-2)} \int_0^t \bbE\big[\big| \Delta_\ell A_3(a)\big|^2 \big]^{\frac{1}{2}} \, da.
\]
The result then follows from dominated convergence argument, as the series 
\[
\sum_{\ell\geq 0} 2^{\ell(2\alpha+2\eta-2)} \int_0^t \bbE\big[\big| \Delta_\ell A_3(a)\big|^2 \big]^{\frac{1}{2}} da
\] 
is seen to be convergent.

{\color{black} 
\section{Proof of Theorem \ref{thm_nonexplosionA}} 
\label{section_appendix_nonexplosion}

We follow closely Shen, Zhu \& Zhu's work \cite{ShenZhuZhu}.

\ssk

\subsection{Setup and strategy\boldmath{.} \hspace{0.15cm} }

Take $2 / 3 <
\alpha < \alpha_{\circ} <\alpha+\eta< 1$, where $\eta$ is given in the assumptions of Theorem \ref{ThmNoiseEnhancement}. For $\gamma \geq 0$ and any $0 < T < \infty$ we
set
\[ \| a \|_{\gamma, \alpha} = \sup_{0 < t < s \leq T} t^{\gamma} \left( \norme{ a_t
   }_{C^\alpha} + \frac{\norme{ a_t - a_s }_{L^{\infty}}}{| t - s |^{\alpha / 2}}
   \right),
   \]
and
\[ | b |_{\gamma, 2 \alpha} = \sup_{0 < t  \leq T} t^{\gamma} \norme{ b_t }_{C^{2
   \alpha}} 
   \qquad \textit{and} \quad
   |b|_{\gamma, 0} = \sup_{0 < t   \leq T} t^{\gamma} \norme{ b_t }_{L^{\infty}}. \]
We write also
\[
 | u |_{t^{\alpha / 2} C^{\alpha / 2}_t L^{\infty}} = \sup_{0 < t < s \leq T} t^{\alpha/2} \frac{\norme{ a_t - a_s }_{L^{\infty}}}{| t - s |^{\alpha / 2}}
\]
To lighten the notations we denote in this appendix by $u$ the function associated with ${\bf u}_{\widehat\xi^+,{\bf v}}$. It will be useful consider it as an element of the larger space of paracontrolled functions $\{ \mathbf{w} = (w, w^{\sharp}) \}$ with $\| w \|_{\alpha / 2, \alpha} < \infty$ and $| w^{\sharp} |_{\alpha / 2, 2 \alpha} <
\infty$.

\medskip

We will see in Proposition \ref{prop_AAAAAA3} below that the paracontrolled structure of $u$ that the explosion time of $(u, u^{\sharp})$ is given by the explosion time of $u^{\sharp} \in t^{\alpha / 2} C_t C^{2 \alpha}$. We will obtain an a posteriori bound on $u^{\sharp}$ by looking at its dynamics. We have first
\begin{equation} \label{EqDynamicsuSharp}
\scrL u^{\sharp} = f (u) \xi - L (f (u) \varolessthan X_{> n})   \hspace{3em}   (u^{\sharp} (0) = u_0) .
\end{equation}
We choose to decompose the right hand side of \eqref{EqDynamicsuSharp} into the sum of a term with positive regularity at positive times and another term. We split accordingly $u^{\sharp} = u^{\sharp}_2 + u^{\sharp}_1$ and identify each term separately.
Given another integer parameter $n_1$ we write $\xi^{(2)}_{> n_1}$ for the high frequencies in space of $\xi^{(2)}$ in a Littlewood-Paley decomposition. We define $u_1^{\sharp}$ and $u_2^{\sharp}$ from their initial conditions, $u_1^\sharp=u_0$ and $u_2^\sharp=0$, and
the equations
\begin{equation}\begin{split}\label{eq_defu1dash}
    \scrL u_1^\# &=  f_o(u_t \vert v_t)\varolessthan \xi_{>n} - \scrL(f_o(u_t \vert v_t)\varolessthan X_{>n}) + (\xi_{>n} \varolessthan f(u_t \vert v_t)) 
    \\
    &\quad+ \delta_z f_o(u_t\vert v_t)u_t' \, (\xi\odot X)_{>n_1}  +  \delta_\mu f_o(u_t \vert v_t)v_t'(\xi\odot\overline{X})_{>n_1} + g(u_t\vert v_t).
    \end{split}\end{equation}
and
\begin{equation} \label{eq_defu2dash}
\begin{split}
\scrL u_2^\# &= f_o(\bfu_t\vert \bfv_t)\odot \xi_{>n}  + f_o(u_t\vert v_t)\xi_{\leq n}  - \delta_zf_o (u_t\vert v_t) f_o(u_t\vert v_t) (\xi\odot X)_{>n_1}   \\
&\quad+ \overline\bbE\big[ \delta_\mu f_o(u_t\vert v_t) v'_t (\xi\odot\overline{X})_{>n_1} \big].
\end{split} \end{equation}

\medskip

\begin{lem} The right hand side of \ref{eq_defu2dash} has positive
regularity at any positive time.
\end{lem}
\medskip

\begin{Dem} Fix a positive time. Using the corrector
${\sf C}$ we have
\begin{align*} 
&f_o (\bfu \vert \bfv ) \odot \xi_{> n}  - \delta_zf_o (u_t\vert v_t) u_t' (\xi\odot X)_{>n_1} - \delta_\mu f_o(u_t\vert v_t)v'_t (\xi\odot\overline{X})_{>n_1} 
\\
&=  {\sf C}\big(\delta_zf_o (u_t\vert v_t) u_t', X,\xi_t\big) + \overline{\bbE}\Big[{\sf C}\big(\partial_2 F(u,v)v', \overline{X},\xi\big)\Big]
   + f (u)^{\sharp} \odot \xi_{> n} 
   \end{align*}
which has indeed positive regularity.
\end{Dem}

\ssk

The cut-off parameters $n$ and $n_1$ will be taken later as some functions of $\widehat{\xi}^+$ and $u$ itself. We will see that
\[ 
| u_1^{\sharp} |_{\alpha / 2, 2 \alpha} {\lesssim 2^{- n (\alpha_{\circ} - \alpha)}}   \| u \|_{\alpha / 2, \alpha} + 1 
\]
and
\[ 
{\color[HTML]{000000}\| u \|_{\alpha / 2, \alpha}} \lesssim \| u_2^{\sharp} \|_{{\color[HTML]{000000}\alpha / 2, \alpha}} + 1 
\] 
and we will finally obtain a control $\| u_2^{\sharp} \|_{\alpha / 2, 2 - \epsilon_1} \lesssim 1$ of $u_2^{\sharp}$ in a norm stronger than the $| \cdot
|_{\alpha / 2, 2 \alpha}$ norm ($2 - 2 \alpha_{\circ} < \epsilon_1 < 2 - 2
\alpha$).

\subsection{Everything depends on $u_2^{\sharp}$\boldmath{.} \hspace{0.15cm}}

We use the dynamics \ref{eq_defu1dash} of $u_1^{\sharp}$ to show that $u_1^{\sharp}$ is controlled in terms of some norms of $u_2^{\sharp}$.

\medskip

\begin{prop} \label{prop_AAAAAA2}
An adequate choice of $n_1$ depending on $\| u \|_{\alpha / 2, \alpha}$, a parameter $0 < \epsilon_1 < 3 \alpha - 2$ and $\widehat{\xi}^+$ ensures that one has both
\begin{equation}\label{eq_propA2}
\begin{split}
\| u_1^{\sharp} \|_{C_T C^{\alpha}}    &\lesssim 1+ \overline\bbE[\norme{\bfv}^2_{\mcE^\beta_T}]   \\
| u_1^{\sharp} |_{\alpha / 2, 2 \alpha} &\lesssim 2^{-n (\alpha_{\circ} - \alpha)}   \| u \|_{\alpha / 2, \alpha}  + \overline{\bbE}[\norme{\bfv}_{\mcE^\beta_T}] .
\end{split} \end{equation}
\end{prop}

\medskip

The first estimate will be used in \ref{eq_conditionn1} below. What matters presently is the second estimate. The proof actually shows that we also have $| u_1^{\sharp}
|_{\gamma, 2 \alpha} {\lesssim 2^{- n (\alpha_{\circ} - \alpha)}}   \| u \|_{\gamma, \alpha} + 1$ for $0 < \gamma < 1$; we will use this fact in the sequel.

\medskip

\begin{Dem} Recall the continuity estimate ($0 < \beta < 1$ and $0
\leq a \leq \gamma$)
\begin{equation} \label{eq_schauder0}
  \| \scrL^{- 1} (h) \|_{\gamma - a, \beta} \lesssim | h |_{\gamma, \beta - 2 + 2
  a}
\end{equation}
-- \textit{First estimate}. With $P_t$ the free propagation operator, write for any fixed time
\begin{equation} \label{eq_defu1dash_bis} \begin{split}
    u^{\sharp}_1(t) & =  P_t (u_0) + \scrL^{- 1}\Big(  f_o(u_t\vert v_t)\varolessthan \xi_{>n} + (\xi_{>n} \varolessthan f_o(u_t\vert v_t))  + \delta_z f_o(u_t\vert v_t)f_o(u_t\vert v_t) \, (\xi\odot X)_{>n_1}
    \\
    &\hspace{3cm}  + \overline{\bbE}\big[ \delta_\mu f_o(u_t\vert v_t)v_t'(\xi\odot\overline{X})_{>n_1}\big] + g_\circ(u_t\vert v_t)\Big) - f_\circ(u_t\vert v_t)\varolessthan X_{>n}
    \\
    & =  P_t(u_0) + \scrL^{- 1} ((1) + (2) + (3) + (4)+(5)) - (6).
  \end{split}\end{equation}
and estimate $u_1^{\sharp}$ in $C^{\alpha}$ at a fixed time. Since $\norme{P_t(u_0) }_{\alpha}$ is of order $\norme{u_0 }_{\alpha}$ independently of the time (in a fixed time interval) it is useless to have some small additive constants in the estimates. We estimate $(1)$ and $(2)$ in $C_TC^{\alpha - 2}$ and use \ref{eq_schauder0} with $\gamma = a = 0$ for these two terms. Since $F$ is bounded they are of size $\| f_o (u\vert v) \|_{\infty} \leq \| F \|_{\infty}$, uniformly in $n$. This is also true for the terms $(5)$ and $(6)$. Since $(\xi\odot X)_{>n_1},(\xi\odot \overline{X})_{>n_1} \in C^{2 \alpha - 2}$ we introduce a parameter $0 < \epsilon_1 < 3 \alpha- 2$ and use \ref{eq_schauder0} to write
\begin{align*} 
\| \scrL^{- 1} (3) \|_{0, \alpha} \lesssim | (3) |_{1 - \alpha + \epsilon_1 / 2, - \alpha + \epsilon_1} 
   &
   \leq \Big( \| u \|_{1 - \alpha + \epsilon_1 / 2, \alpha} + \overline\bbE [ \| v \|_{1 - \alpha + \epsilon_1 / 2, \alpha} ] \Big)  |  (\xi\odot X)_{> n_1} |_{- \alpha + \epsilon_1} 
   \\
   & \lesssim  \big( \| u \|_{1 - \alpha + \epsilon_1 / 2, \alpha} + \overline\bbE\big[ \| \bfv \|_{\mcE^\beta_T}\big] \big) \, 2^{- n_1 (3 \alpha - 2 - \epsilon_1)} . 
\end{align*}
Likewise one has
\[ 
\| \scrL^{- 1} (4) \|_{0, \alpha} \lesssim  \big( \| u \|_{1 - \alpha + \epsilon_1 / 2, \alpha} + \overline\bbE\big[ \| \bfv \|^2_{\mcE^\beta_T}\big] \big) \, 2^{- n_1 (3 \alpha - 2 - \epsilon_1)} .
\]
We choose $n_1$ here as a function of $u$ to ensure that
\begin{equation} \label{eq_conditionn1}
  2^{- n_1 (3 \alpha - 2 - \epsilon_1)} \| u \|_{1 - \alpha + \epsilon_1 / 2,
  \alpha} \simeq 1.
\end{equation}

\ssk

-- \textit{Second estimate}. The quantity $| u_1^{\sharp} |_{\alpha / 2, 2
\alpha} \leq \| u_1^{\sharp} \|_{\alpha / 2, 2 \alpha}$ is bounded from above by a
constant multiple of
\begin{align*} 
\norme{ u_0 }_{\alpha} &+ \big\| \scrL^{- 1} (f_o (u\vert v) \varolessthan \xi_{> n}) - f_o(u\vert v) \varolessthan X_{> n} \big\|_{\alpha / 2, 2 \alpha}   \\
   &+ \Big(| u|_{\alpha / 2, \alpha}  +  \overline\bbE[|v|_{\alpha / 2, \alpha}] \Big) \Big(| \xi_{> n} |_{0, \alpha - 2} + | (\xi\odot X)_{> n_1} |_{0, 2 \alpha - 2} \Big)   \\
   &+ \overline\bbE\Big[|u|_{\alpha / 2, \alpha}+| v |_{\alpha / 2,\alpha} | (\xi\odot \overline{X})_{> n_1} |_{0, 2 \alpha - 2} \Big] . 
\end{align*}
A commutator lemma gives an upper bound of the second term by a constant
multiple of 
\[
(\| u \|_{\alpha / 2, \alpha} |+ \overline{\bbE}[\| v \|_{\alpha / 2, \alpha}] )| \xi_{> n} |_{0, \alpha - 2} \lesssim (\| u \|_{\alpha / 2, \alpha}+\overline\bbE [\| \bfv \|_{\mcE^\beta_T} ] ) 2^{- n (\alpha_{\circ} - \alpha)}.
\]
Now since $\alpha > 3 \alpha - 2 - \epsilon_1$ and $2 \alpha - 2 \leq - \alpha + \epsilon_1$, because $0 < \epsilon_1 < 3 \alpha - 2$, the condition \ref{eq_conditionn1} implies that
\[
| u |_{\alpha / 2, \alpha} | (\xi\odot X)_{> n_1} |_{0, 2 \alpha - 2} \lesssim \| u \|_{1 - \alpha + \epsilon_1 / 2, \alpha} 2^{- n_1 (3 \alpha - 2 - \epsilon_1)} \lesssim 1,
\]
idem for the term involving $(\xi\odot\overline X )$. This gives in the end an upper bound $2^{- n (\alpha_{\circ} - \alpha)} \|u\|_{\alpha / 2, \alpha} + \overline{\bbE}\big[\norme{\bfv}_{\mcE^\beta_T}\big]$.
\end{Dem}

\medskip

\begin{prop} \label{prop_AAAAAA3}
For the same parameter $\epsilon_1$ as in Proposition \ref{prop_AAAAAA2} one has
\begin{equation} \label{eq_propA3}
\| u \|_{\alpha / 2, \alpha} \lesssim \| u_2^{\sharp} \|_{\alpha / 2, \alpha} + 1.
\end{equation}
\end{prop}
\medskip

The proof shows that we have more generally $\| u \|_{\gamma, \alpha} \lesssim \| u_2^{\sharp} \|_{\gamma, \alpha} + 1$ for $0 < \gamma < 1$. We will use this fact in the sequel. It follows from \ref{eq_propA2} and \ref{eq_propA3} that $| u_1^{\sharp} |_{\alpha / 2, 2 \alpha} {\lesssim 2^{- n (\alpha_{\circ} - \alpha)}}   \| u_2^{\sharp} \|_{\alpha / 2, \alpha} + 1$, so $| u_2^{\sharp} |_{\alpha / 2, 2 \alpha} \lesssim \| u_2^{\sharp} \|_{\alpha / 2, \alpha} + 1$.

\medskip

\begin{Dem} We first look at the size in $t^{\alpha / 2} C^{\alpha /
2}_t L^{\infty}$ of a number of quantities. One has first from the
paracontrolled structure of $u$ that
\[
| u |_{t^{\alpha / 2} C^{\alpha / 2}_t L^{\infty}} \lesssim | f_o (u\vert v) \varolessthan X_{> n} |_{t^{\alpha / 2} C^{\alpha / 2}_t L^{\infty}} + | u^{\sharp} |_{t^{\alpha / 2} C^{\alpha / 2}_t L^{\infty}}
\]
with
\begin{equation}\label{eq_choixaftern} \begin{split}
  | f_o (u\vert v) \varolessthan X_{> n} |_{t^{\alpha / 2} C^{\alpha / 2}_t
  L^{\infty}} \lesssim \| F \|_{C^1} \big( | u |_{t^{\alpha / 2} C^{\alpha / 2}_t 
  L^{\infty}} + \overline\bbE[\norme{v}_{\mcE^\beta}]\big) | \xi_{> n} |_{C_T C^{\alpha - 2}} + \| F \|_{L^{\infty}} |
  \xi_{> n} |_{C_T C^{\alpha - 2}} .
\end{split} \end{equation}
Since $\| \xi_{> n} \|_{C_T C^{\alpha - 2}} \lesssim 2^{- n (\alpha_{\circ} - \alpha)} | \xi |_{C_T C^{\alpha_{\circ} - 2}}$, choosing $n$ large enough as a function of $\xi$ ensures that 
\[
| u |_{t^{\alpha / 2} C^{\alpha / 2}_t L^{\infty}} \lesssim | u^{\sharp} |_{t^{\alpha / 2} C^{\alpha / 2}_t L^{\infty}} +\overline\bbE\big[\norme{\bfv}_{\mcE^\beta_T}\big] +  1.
\]
We now estimate $| u_1^{\sharp} |_{t^{\alpha / 2} C^{\alpha / 2}_t L^{\infty}}$ by looking at its different pieces in \ref{eq_defu1dash_bis}. First we have $| P_t (u_0) |_{t^{\alpha / 2} C^{\alpha / 2}_t
L^{\infty}} \lesssim | u_0 |_{C^{\alpha}}$. From \ref{eq_schauder0} we have 
\[ 
\big| \scrL^{- 1} ((1) +\dots + (5)) \big|_{t^{\alpha / 2} C^{\alpha / 2}_t L^{\infty}} \lesssim \big| (1) + \dots + (5) \big|_{\alpha / 2, \alpha - 2}.
\]
One has $| (1) + (2) |_{\alpha / 2, \alpha - 2} \lesssim ( | u |_{t^{\alpha / 2} C_T L^{\infty}}  + \overline\bbE[\norme{\bfv}_{\mcE^\beta_T}] ) \| \xi_{> n} \|_{C_T C^{\alpha - 2}}$. For $(3)$, as we have
from \ref{eq_conditionn1} $ \| u \|_{1 - \alpha + \epsilon_1 / 2, \alpha}  | (\xi\odot X)_{> n_1} |_{0, - \alpha + \epsilon_1} \lesssim 1$, with $1 - \alpha + \epsilon_1 / 2 < \alpha / 2$, we see that $| (1) + \dots+ (5) |_{\alpha / 2, \alpha - 2}$ is bounded by $1+\overline{\bbE}[\norme{v}^2_{\mcE^\beta_T}]$. We thus have 
\begin{align*} 
| u_1^{\sharp} |_{t^{\alpha / 2} C^{\alpha / 2}_t L^{\infty}} &\lesssim 1 + \overline{\bbE}\big[\norme{\bfv}^2_{\mcE^\beta_T}\big]
   +| u |_{t^{\alpha / 2} C^{\alpha / 2}_t L^{\infty}} | \xi_{> n} |_{C_T
   C^{\alpha - 2}} 
   \\
   &\lesssim 1 + \overline{\bbE}\big[\norme{\bfv}^2_{\mcE^\beta_T}\big] + | \xi_{> n} |_{C_T C^{\alpha - 2}} \big( | u_1^{\sharp} |_{t^{\alpha / 2} C^{\alpha / 2}_t L^{\infty}} + | u_2^{\sharp} |_{t^{\alpha / 2} C^{\alpha / 2}_t L^{\infty}}\big). 
\end{align*}

Taking $n$ even larger if necessary enventually gives
\[ | u_1^{\sharp} |_{t^{\alpha / 2} C^{\alpha / 2}_t L^{\infty}} \lesssim 1 + \overline{\bbE}\big[\norme{\bfv}^2_{\mcE^\beta_T}\big] +
   | \xi_{> n} |_{C_T C^{\alpha - 2}} | u_2^{\sharp} |_{t^{\alpha / 2}
   C^{\alpha / 2}_t L^{\infty}} \]
so
\begin{equation}
  | u |_{t^{\alpha / 2} C^{\alpha / 2}_t L^{\infty}} \lesssim | u_1^{\sharp}
  |_{t^{\alpha / 2} C^{\alpha / 2}_t L^{\infty}} + | u_2^{\sharp} |_{t^{\alpha
  / 2} C^{\alpha / 2}_t L^{\infty}} + 1 \lesssim 1 + | u_2^{\sharp}
  |_{t^{\alpha / 2} C^{\alpha / 2}_t L^{\infty}} +\overline{\bbE}\big[\norme{\bfv}^2_{\mcE^\beta_T}\big]
\end{equation}
With $| u |_{\alpha / 2, \alpha} \lesssim 1 + |
u^{\sharp} |_{\alpha / 2, \alpha}$ from the paracontrolled structure of $u$,
it follows that
\begin{align*} \| u \|_{\alpha / 2, \alpha} = | u |_{\alpha / 2, \alpha} + | u
   |_{t^{\alpha / 2} C^{\alpha / 2}_t L^{\infty}} &\lesssim 1 + | u^{\sharp}
   |_{\alpha / 2, \alpha} + | u_2^{\sharp} |_{t^{\alpha / 2} C^{\alpha / 2}_t
   L^{\infty}}  + \overline{\bbE}\big[\norme{\bfv}^2_{\mcE^\beta_T}\big]
   \\
   &\lesssim 1 + \| u_2^{\sharp} \|_{\alpha / 2, \alpha} + \overline{\bbE}\big[\norme{\bfv}^2_{\mcE^\beta_T}\big]. \end{align*}
\end{Dem}

\subsection{Bounds on $u_2^{\sharp}$\boldmath{.} \hspace{0.15cm}}

We will prove that $\| u_2^{\sharp} \|_{\alpha / 2, 2 - \epsilon_1} \lesssim 1$ and $| u_2^{\sharp} |_{L^{\infty}_{t, x}} \lesssim 1 + \| u_2^{\sharp} \|^r_{\alpha / 2, 2 - {\color[HTML]{000000}\epsilon_1}}$ for some appropriate choice of $0 < \epsilon_1 < 1$ and exponent $0 < r < 1$.

\ssk

\begin{prop} \label{prop_AAAAAA4}
Assume $4 / 5 < \alpha < 1$. If one chooses the constant $0 < \epsilon_1 < 3 \alpha - 2$ from Proposition 2 such that $\epsilon_1 > 2 (1 - \alpha_{\circ})$, for all $1 - \alpha < \epsilon_2
\leq \alpha$ there is a vector-valued function $g_u \in C^{1 + \epsilon_2}$ that depends on $u$ such that one has the uniform in time estimate
\begin{equation}\label{eq_refinedparalinearisation} 
\big\| f_o (\bfu\vert\bfv )^{\sharp} \odot \xi_{> n} - (g_u \odot \xi_{> n}) \nabla u^{\sharp}\big\|_{\infty} \lesssim \| \xi_{> n} \|_{\alpha - 2} \left\{ 1 + | u
   |_{\alpha} + \left| {u_1^{\sharp}}  \right|_{2 \alpha} + \left|
   {u_2^{\sharp}}  \right|_{2 - \epsilon_1}  \right\}\Big(1+\overline\bbE[\norme{\bfv}_{\mathcal{E}_T^\beta}^2]\Big) \end{equation}
\textit{pointwise in time, and the function $g_u$ has a decomposition $g_u = g_u^1 + g_u^2$ where}
\begin{equation}
  | g_u^1 |_{1 + \alpha} \lesssim 1 + | u_1^{\sharp} |_{\alpha} + \overline{\bbE}[\norme{\bfv}_{\mcE^\beta_T}] , \qquad |
  g_u^2 |_{1 + \epsilon_2} \lesssim_{\kappa} | u_2^{\sharp} |_{1 +
  \kappa}^{\epsilon_2}
\end{equation}
\textit{for any} $0 < \kappa \leq 2 \alpha - 1$.
\end{prop}

\ssk

\begin{Dem}
Applying Lemma 4.1 from \cite{ShenZhuZhu} to $F$, and setting the decomposition $v^\# = v_1^\# + v_2^\#$ with $v_2^\#=0$ (so as to paralinearise with respect to $\bfv$), we obtain
\begin{align*}
& \norme{F(\bfu_t \vert\bfv_t)^\#\odot \xi_{>n}  -  ( h_1(\bfu_t\vert \bfv_t)\odot\xi_{>n})\nabla u^\# + ( h_2(\bfu_t\vert\bfv_t)\odot\xi_{>n}) \nabla v^\# }_{L^\infty} 
\\
&\hspace{4.5cm}\lesssim \norme{\xi_{>n}}_{\alpha-2} \big(1 + \norme{u}_\alpha + \norme{u_1^\#}_{2\alpha}   + \norme{u_2^\#}_{2-\eps} + \norme{v}_\alpha+\norme{v^\#}_{2\alpha}\big)
\end{align*}
where both $h_1$ and $h_2$ decomposes as $h_i=h_i^1+h_i^2$ with $| h_i^1 |_{1 + \alpha} \lesssim 1 + | u_1^{\sharp} |_{\alpha} + \norme{v^\#}_{\alpha} $ and $|
  h_i^2 |_{1 + \epsilon_2} \lesssim_{\kappa} | u_2^{\sharp} |_{1 +   \kappa}^{\epsilon_2}$. Integrating this decomposition on $\overline{\Omega}$ gives indeed the proposition with $g_u^j=\overline\bbE[h_1^j(\bfu\vert\bfv)]$, as the term $\overline{\bbE}[h_2(\bfu_t\vert\bfv_t)\odot\xi_{>n}) \nabla v^\#]$ can be estimated by the right hand side of \eqref{eq_refinedparalinearisation}.
\end{Dem}

\ssk

The parameters $\kappa$ and $\epsilon_2$ will be chosen later. We see in particular, using Proposition \ref{prop_AAAAAA2}, that
\begin{equation}\label{eq_estimeegu1bis}
  | g_u^1 |_{1 + \alpha} \lesssim 1 + \overline{\bbE}[\norme{\bfv}_{\mcE^\beta_T}].
\end{equation} 
Write here as a shorthand notation
\[ 
\scrL u^{\sharp}_2 = (g_u \odot \xi_{> n}) \nabla u_2^{\sharp} + (\star) 
\]
where 
\[
\Delta_u \defeq f_o(\bfu,\bfv)^{\sharp} \odot \xi_{> n} - (g_u \odot \xi_{> n}) \nabla u_2^{\sharp}
\]
and 
\begin{align*} 
(\star) &= (\delta_z f_o f_o) (u\vert v) \{ (\xi\odot X)_{> n} - (\xi\odot X)_{> n_1} + \zeta_n \}   \\
&\quad+ \delta_zf_o (u\vert v) \big({\sf C} (f_o (u\vert v), X_{> n}, \xi_{> n}) + u^{\sharp} \odot \xi_{> n}\big) 
   \\
&\quad+   \overline{\bbE}\Big[{\sf C}\big(\partial_2 F(u,v)v', \overline{X},\xi\big)\Big]
   + {\sf C} (\delta_zf_o (u\vert v), u, \xi_{> n}) + \Delta_u + (g_u\odot \xi_{> n}) \nabla u_1^{\sharp} .
\end{align*}
The refined maximum principle (Theorem 4.2 in \cite{ShenZhuZhu}) applied to the parabolic dynamics of $u_2^{\sharp}$ tells us that
\begin{equation} \label{eq_staru2dash}
  | u_2^{\sharp} |_{L^{\infty}_{t, x}} \lesssim | (\star) |_{\alpha / 2, 0} .
\end{equation}

\medskip

\begin{lem} \label{lem_AAAAAA5} Assume $0 < \epsilon_1 < 2 - 2 \alpha$. There
is an explicit exponent $0 < r < 1$ such that one has for all $0 < \mu$
sufficently small
\begin{equation} \label{eq_staru2dashu2dash}
  | (\star) |_{\alpha / 2, 0} \lesssim \big(1+\overline{\bbE}\big[\norme{\bfv}_{\mcE^\beta_T}\big]^2\big) \Big( 1 + \| u_2^{\sharp} \|^r_{\alpha / 2, 2 - \epsilon_1} + \mu \| u_2^{\sharp} \|_{L^{\infty}} \Big).
\end{equation}
\end{lem}

\medskip

Together \ref{eq_staru2dash} and \ref{eq_staru2dashu2dash} give
\[ | u_2^{\sharp} |_{L^{\infty}_{t, x}} \lesssim \big(1+\overline{\bbE}\big[\norme{\bfv}_{\mcE^\beta_T}\big]^2\big) \big( 1 + \| u_2^{\sharp} \|^r_{\alpha / 2, 2 - \epsilon_1} \big) . \]
   
\medskip   
   
\begin{Dem} We estimate each term separately.

\ssk

-- Since $F \in C^1_b$ the first term is bounded by a constant multiple of $2^{(2 - 2 \alpha) n} + 2^{(2 - 2 \alpha) n_1}$.

\ssk

-- For the second term one has from Proposition \ref{prop_AAAAAA3}
\begin{align*}
  | \delta_z f_o (u\vert v) ({\sf C} (f_o (u\vert v), X_{> n}, \xi_{> n}) + u^{\sharp} \odot \xi_{>
  n}) |_{\alpha / 2, 0} & \lesssim   \big(| u |_{\alpha / 2, \alpha} + |
  u^{\sharp} |_{\alpha / 2, \alpha}  + \overline{\bbE}[\norme{\bfv}_{\mcE^\beta_T}]   \big) 2^{- (\alpha_{\circ} - \alpha) n}
  \\
  & \lesssim   \big(1 + | u_2^{\sharp} |_{\alpha / 2, \alpha} + \overline{\bbE}[\norme{\bfv}_{\mcE^\beta_T}] \big) 2^{-
  (\alpha_{\circ} - \alpha) n}
\end{align*}

\ssk 

-- For the third term we write as above
\[ 
{\sf C} (\delta_z f_o(u\vert v), u, \xi_{> n}) = {\sf C} ( \delta_z f (u\vert v), f (u\vert v) \varolessthan X_{> n}, \xi_{> n}) + {\sf C} (\delta_z f_o (u\vert v), u^{\sharp}, \xi_{> n}) 
\]
and use that $F \in C^2_b$, put the time weight on $| X_{> n} |_{\alpha}$ and see with Proposition \ref{prop_AAAAAA3} that $| {\sf C} (\delta_zf_o(u\vert v), u, \xi_{> n}) |_{\alpha /
2, 0}$ is bounded above by
\begin{equation*} \begin{split}
&\Big(1 +\overline{\bbE}[\norme{\bfv}_{\mcE^\beta_T}] +  \| u^{\sharp}_2 \|_{\alpha / 2, \alpha} + | u_2^{\sharp} |_{\alpha / 2, 2 \alpha}\Big) \, 2^{- (\alpha_{\circ} - \alpha) n}   \\
&\lesssim \Big(1 +\overline{\bbE}[\norme{\bfv}_{\mcE^\beta_T}]+ \| u^{\sharp}_2 \|_{\alpha / 2, \alpha} + | u_2^{\sharp} |_{\alpha / 2, 2 - \epsilon_1}\Big) \, 2^{- (\alpha_{\circ} - \alpha) n} 
\end{split} \end{equation*} 
since $0 < \epsilon_1 < 2 - 2 \alpha$. Likewise we have
\[
\Big|\overline{\bbE}\Big[{\sf C}\big(\partial_2 F(u,v)v', \overline{X},\xi\big)\Big] \Big|_{\alpha/2,0} \lesssim (1 +\overline{\bbE}[\norme{\bfv}_{\mcE^\beta_T}] + \| u^{\sharp}_2 \|_{\alpha / 2, \alpha} + | u_2^{\sharp} |_{\alpha / 2, 2 - \epsilon_1}) 2^{- (\alpha_{\circ} - \alpha) n}.
\]

\ssk

-- We already know from Proposition \ref{prop_AAAAAA4} that $\| \Delta_u \|_{\alpha / 2, 0}$ is bounded above by
\begin{align*}
 &2^{- (\alpha_{\circ} - \alpha) n} \left( 1 + | u |_{\alpha / 2, \alpha} +
   | {u_1^{\sharp}} |_{\alpha / 2, 2 \alpha} + |
   {u_2^{\sharp}} |_{\alpha / 2, 2 - \epsilon_1} +\overline{\bbE}\big[\norme{\bfv}_{\mcE^\beta_T}\big] \right) 
   \\ 
&\lesssim 
   2^{-(\alpha_{\circ} - \alpha) n} \Big(1 + \| u^{\sharp}_2 \|_{\alpha / 2, \alpha} +
   | u_2^{\sharp} |_{\alpha / 2, 2 - \epsilon_1} +\overline{\bbE}\big[\norme{\bfv}_{\mcE^\beta_T}\big]\Big)
\end{align*}
from Proposition \ref{prop_AAAAAA2} and Proposition \ref{prop_AAAAAA3}.

\ssk

-- Recall $(1 + \epsilon_2) + (\alpha - 2) > 0$ . Last we have from Proposition \ref{prop_AAAAAA4} and \ref{eq_estimeegu1bis}
 \begin{align*}
     | (g_u \odot \xi_{> n}) \nabla u_1^{\sharp} |_{\alpha / 2, 0} & \leq 
     \max_{0 \leq t \leq T} t^{\alpha / 2} \| g_u \odot \xi_{> n}
     \|_{L^{\infty}} \| {u_1^{\sharp}}  (t) \|_{C^1}
     \\
     & \lesssim  2^{- n (\alpha_{\circ} - \alpha)} \big(1 + | u_2^{\sharp}
     |^{\epsilon_2}_{\alpha / 2, 1 + \kappa} + \overline{\bbE}\big[\norme{\bfv}_{\mcE^\beta_T}\big]  \big) \max_{0 \leq t \leq T} t^{\alpha
     / 2 (1 - \epsilon_2)} | {u_1^{\sharp}}  (t) |_{1 + \kappa}
     \\
     & \lesssim  2^{- n (\alpha_{\circ} - \alpha)} \big( 1 + | u_2^{\sharp}
     |^{\epsilon_2}_{\alpha / 2, 1 + \kappa} +  \overline{\bbE}\big[\norme{\bfv}_{\mcE^\beta_T}\big]  \big)  \| {u_1^{\sharp}} 
     \|_{C_T C^{\alpha}}^{(2 \alpha - 1 - \kappa) / \alpha} 
     \\
     &\hspace{6.5cm} \times \max_{0 \leq
     t \leq T} t^{{\alpha / 2 (1 - \epsilon_2)}} |
     {u_1^{\sharp}}  (t) |^{(1 + \kappa - \alpha) / \alpha}_{2 \alpha}
   \end{align*}
We used the elementary interpolation result
$ | {u_1^{\sharp}}  (t) |_{1 + \kappa} \lesssim |
   {u_1^{\sharp}}  (t) |^{(2 \alpha - 1 - \kappa) / \alpha}_{\alpha}
   | {u_1^{\sharp}}  (t) |^{(1 + \kappa - \alpha) / \alpha}_{ 2
   \alpha} $
   
in the third inequality. Choose $\kappa$ to have $(1 + \kappa
- \alpha) / \alpha = \alpha (1 - \epsilon_2) / 2$. Use the fact that $\|
u_1^{\sharp} \|_{C_T C^{\alpha}} \lesssim 1$, from \ref{prop_AAAAAA2}, and 
\[
| u_1^{\sharp}|_{\alpha / 2, 2 \alpha} \lesssim 2^{- n (\alpha_{\circ} - \alpha)} \| u_2^{\sharp} \|_{\alpha / 2, \alpha} + 1 + \overline{\bbE}[\norme{\bfv}_{\mcE^\beta_T}]
\] 
from Proposition \ref{prop_AAAAAA3}, to obtain
\begin{align*}
     &| (g_u \odot \xi_{> n}) \nabla u_1^{\sharp} |_{\alpha / 2, 0}
     \\
     &\qquad\quad \lesssim
      2^{- (\alpha_{\circ} - \alpha) n} (1 + | u_2^{\sharp}
     |^{\epsilon_2}_{\alpha / 2, 1 + \kappa})  \big(2^{- (\alpha_{\circ} - \alpha) n} \| u_2^{\sharp} \|_{\alpha / 2, \alpha} + 1\big)^{(1 - \epsilon_2)}\big(1+\overline{\bbE}\big[\norme{\bfv}_{\mcE^\beta_T}\big]^2\big)
     \\
     & \qquad\quad\lesssim 2^{- (\alpha_{\circ} - \alpha) n} (1 + | u_2^{\sharp}
     |^{\epsilon_2}_{\alpha / 2, 1 + \kappa})  \big(2^{- (\alpha_{\circ} - \alpha)
     (1 - \epsilon_2) n} \| u_2^{\sharp} \|^{{(1 - \epsilon_2)}}_{\alpha / 2, \alpha} + 1\big)\big(1+\overline{\bbE}\big[\norme{\bfv}_{\mcE^\beta_T}\big]^2\big)
   \end{align*}
   
Altogether this gives for $| (\star) |_{\alpha / 2, 0}$ the upper bound
\begin{align*}
     (1+ \overline{\bbE}[\norme{\bfv}_{\mcE^\beta_T}]^2)\Big( &1  +  2^{(2 - 2 \alpha) n} + 2^{(2 - 2 \alpha) n_1}
     \\
     & +  2^{- (\alpha_{\circ} - \alpha) n} (1 + \| u^{\sharp}_2 \|_{\alpha
     / 2, \alpha} + | u_2^{\sharp} |_{\alpha / 2, 2 - \epsilon_1})
     \\
     & +  2^{- (\alpha_{\circ} - \alpha) n} (1 + |
     u_2^{\sharp} |^{\epsilon_2}_{\alpha / 2, 1 + \kappa})  (2^{-
     (\alpha_{\circ} - \alpha) (1 - \epsilon_2) n} \| u_2^{\sharp} \|^{(1 -
     \epsilon_2)}_{\alpha / 2, \alpha} + 1) \Big) .
   \end{align*} 
Set
\[ \ominus (\mu) \defeq c (\mu) \| u_2^{\sharp} \|_{\alpha / 2, 2 -
   \epsilon_1} + \mu \| u_2^{\sharp} \|_{L^{\infty}} . \]
We have
\[ \begin{array}{lll}
     | u_2^{\sharp} |_{\alpha / 2, \alpha} \lesssim | u_2^{\sharp} |_{\alpha /
     2, 2 - \epsilon_1}^{\alpha / (2 - \epsilon_1)} \| u_2^{\sharp}
     \|_{L^{\infty}}^{1 - \alpha / (2 - \epsilon_1)} & \lesssim & \|
     u_2^{\sharp} \|_{\alpha / 2, 2 - \epsilon_1}^{\alpha / (2 - \epsilon_1)}
     \| u_2^{\sharp} \|_{L^{\infty}}^{1 - \alpha / (2 - \epsilon_1)} \lesssim
     \ominus (\mu)
   \end{array} \]
from the interpolation bounds and Young inequality, for $\mu > 0$ small and an associated constant $c (\mu) > 0$. Similarly we have
\[ \begin{array}{lll}
     | u_2^{\sharp} |_{\alpha / 2, 1 + \kappa} \lesssim
     | u_2^{\sharp} |_{\alpha / 2, 2 - \epsilon_1}^{(1 + \kappa) / (2 -
     \epsilon_1)} \| u_2^{\sharp} \|_{L^{\infty}}^{1 - (1 + \kappa) / (2 -
     \epsilon_1)} & \lesssim & \| u_2^{\sharp} \|_{\alpha / 2, 2 -
     \epsilon_1}^{(1 + \kappa) / (2 - \epsilon_1)} \| u_2^{\sharp}
     \|_{L^{\infty}}^{1 - (1 + \kappa) / (2 - \epsilon_1)} \lesssim \ominus
     (\mu) .
   \end{array} \]
We can also use  $0 < \alpha \leq \beta < 2, 0 \leq \rho \leq \frac{\beta \gamma}{\alpha}$ and $0 < \mu < 1$ the useful interpolation bound
\[
  \| h \|_{\gamma, \alpha} \lesssim \| h \|_{\rho, \beta}^{\alpha / \beta}  \|
  h \|_{\mathbb{L}^{\infty}}^{1 - \alpha / \beta} \lesssim c (\mu) \| h
  \|_{\rho, \beta} + \mu \| h \|_{\mathbb{L}^{\infty}} .
\]
It gives $\| u_2^{\sharp} \|_{\alpha / 2, \alpha} \lesssim \ominus (\mu)$. We also have $| u_2^{\sharp} |_{\alpha / 2, 2 - \epsilon_1} \leq \| u_2^{\sharp} \|_{\alpha / 2, 2 - \epsilon_1}$ and $| u_2^{\sharp} |_{\alpha / 2, 1 + \kappa} \leq \| u_2^{\sharp} \|_{\alpha / 2, 1 + \kappa} \lesssim \ominus (\mu)$. We thus get for $| (\star) |_{\alpha / 2, 0}$ the upper bound
 \begin{align*}
     \big(1+\overline{\bbE}[\norme{\bfv}_{\mcE^\beta_T}]^2\big) \Big(&1  +  2^{(2 - 2 \alpha) n} + 2^{(2 - 2 \alpha) n_1} + 2^{-
     (\alpha_{\circ} - \alpha) n} \ominus (\mu)
     \\
     & +  2^{- (\alpha_{\circ} - \alpha) n} (1 + \ominus (\mu)^{\epsilon_2})
     (2^{- (\alpha_{\circ} - \alpha) (1 - \epsilon_2) n} \ominus (\mu)^{1 -
     \epsilon_2} + 1\Big)
   \end{align*} 
Now we have for all $0 < a < 1$ and $x > 0$ the elementary inequality $1 + x^{\epsilon_2} + a^{1 - \epsilon_2} x^{1 - \epsilon_2} + a^{1 - \epsilon_2} x \lesssim a^{- \epsilon_2} + a^{1 - \epsilon_2} x$. By choosing $2 - 2 \alpha_{\circ} < \epsilon_2 = \epsilon_1 < 2 - 2 \alpha$, it gives in the end the estimate 
\[ 
| (\star) |_{\alpha / 2, 0} \lesssim \begin{array}{lll}
     1 & + & 2^{(2 - 2 \alpha) n} + 2^{(2 - 2 \alpha) n_1} + 2^{-
     (\alpha_{\circ} - \alpha) n} \ominus (\mu) + 2^{- (\alpha_{\circ} -
     \alpha) (2 - \epsilon_1) n} \ominus (\mu) .
   \end{array} \]
The parameter $n$ was first chosen after \ref{eq_choixaftern} as a function of $\widehat{\xi}^+$. We choose $n$ even larger to have 
\[ 
2^{- (\alpha_{\circ} - \alpha) n} \| u_2^{\sharp} \|_{\alpha / 2, 2 -
   \epsilon_1} \simeq \| u_2^{\sharp} \|_{\alpha / 2, 2 - \epsilon_1}^r + C 
\]
for a large constant $C > 0$ and a constant $r$ that we chose below. We chose $n_1$ in \ref{eq_conditionn1}. We can choose an even larger $n_1$ that ensures that $2^{- n_1 (3 \alpha - 2 - \epsilon_1) } \| u_2^{\sharp} \|_{1 - \alpha + \epsilon_1 / 2, \alpha} \simeq 1$. Then we have 
\[ 
2^{(2 - 2 \alpha) n_1} \simeq \max (\| u_2^{\sharp} \|_{1 - \alpha +
   \epsilon_1 / 2, 2 - \epsilon_1}, C)^{(2 - 2 \alpha) / (3 \alpha - 2 -
   \epsilon_1)} 
\]
with $(2 - 2 \alpha) / (3 \alpha - 2 - \epsilon_1) < 1$ for $0 < \epsilon_1 < 5 \alpha - 4$. This finally gives for $| (\star) |_{\alpha / 2, 0}$ the upper bound
\begin{equation}\label{eq_AAA23}
(1+\overline{\bbE}[ \norme{\bfv}_{\mcE^\beta_T} ]^2)\Big(1 + \| u_2^{\sharp} \|_{\alpha / 2, 2 -
  \epsilon_1}^{2 (1 - r) (1 - \alpha) / (\alpha_{\circ} - \alpha)} + \|
  u_2^{\sharp} \|_{1 - \alpha + \epsilon_1 / 2, \alpha}^{(2 - 2 \alpha) / (3
  \alpha - 2 - \epsilon_1)} + \| u_2^{\sharp} \|_{\alpha / 2, 2 -
  \epsilon_1}^r + \mu \| u_2^{\sharp} \|_{L^{\infty}}\Big) .
\end{equation}
We $\| u_2^{\sharp} \|_{1 - \alpha + \epsilon_1 / 2, \alpha} \lesssim \|
u_2^{\sharp} \|_{\alpha / 2, 2 - \epsilon_1}^{\alpha / (2 - \epsilon_1)} \|
u_2^{\sharp} \|^{1 - \alpha / (2 - \epsilon_1)}_{L^{\infty}}$ from (22) since
$\alpha / 2 \leq (1 - \alpha + \epsilon_1 / 2) (2 - \epsilon_1) / \alpha$ for
an appropriate choice of $\epsilon_1$ small. We then have from Young
inequality
\begin{equation}\label{eq_AAA24}
  \| u_2^{\sharp} \|_{1 - \alpha + \epsilon_1 / 2, \alpha}^{(2 - 2 \alpha) /
  (3 \alpha - 2 - \epsilon_1)} \lesssim c (\mu) \| u_2^{\sharp} \|_{\alpha /
  2, 2 - \epsilon_1}^s + \mu \| u_2^{\sharp} \|_{L^{\infty}}
\end{equation}
with $s = \frac{1}{2 - \epsilon_1}  \frac{2 \alpha (1 - \alpha)}{3 \alpha - 2
- \epsilon_1 - (2 - 2 \alpha) (1 - \alpha / 2 - \epsilon_1)} < 1$. Taking $1 -
\frac{\alpha_{\circ} - \alpha}{2 (1 - \alpha)} < r < 1$, all the exponents of
$\| u_2^{\sharp} \|_{\alpha / 2, \alpha}$ in \ref{eq_AAA23} and \ref{eq_AAA24} are strictly
smaller than $1$. 
\end{Dem}

Schauder estimate and Lemma \ref{lem_AAAAAA5} tell us that
\begin{equation} \label{eq_AAA25}
\begin{split}
  \| u_2^{\sharp} \|_{\alpha / 2, 2 - \epsilon_1} &\lesssim | (g_u \odot \xi_{>
  n}) \nabla u_2^{\sharp} + (\star) |_{\alpha / 2, 0}   \\
  &\lesssim | (g_u \odot
  \xi_{> n}) \nabla u_2^{\sharp} |_{\alpha / 2, 0} + (1 + \| u_2^{\sharp}
  \|^r_{\alpha / 2, 2 - \epsilon_1}) \big(1+\overline{\bbE}\big[\norme{\bfv}_{\mcE^\beta_T}\big]^2\big) .
\end{split} \end{equation}

\medskip

\begin{lem} \label{lem_AAAAAA4}One has
\begin{equation}
  | (g_u \odot \xi_{> n}) \nabla u_2^{\sharp} |_{\alpha / 2, 0} \lesssim \big(\mu
  \| u_2^{\sharp} \|_{\alpha / 2, 2 - \epsilon_1} + \| u_2^{\sharp}
  \|^{r'}_{\alpha / 2, 2 - \epsilon_1}\big) \big(1+\overline{\bbE}[\norme{\bfv}_{\mcE^\beta_T}]\big)
\end{equation}
for a constant $\mu > 0$ that can be chosen arbitrarily small and
for some other exponent $0 < r' < 1$.
\end{lem}
\medskip

This inequality allows to conclude from \ref{eq_AAA25} that $\| u_2^{\sharp} \|_{\alpha
/ 2, 2 - \epsilon_1} \lesssim 1$. As said above, the bounds from Propositions \ref{prop_AAAAAA2} and \ref{prop_AAAAAA3} then
show that the paracontrolled distribution $(u, u^{\sharp})$ remains bounded.

\medskip

\begin{Dem}
To see that the lemma holds we write at any fixed time $0 \leq t \leq T$
\[ 
t^{\alpha / 2} | (g_u \odot \xi_{> n}) |_0  | u_2^{\sharp} |_{1 + \kappa} \lesssim t^{\alpha / 2} (1 + | u_2^{\sharp} |_{1 + \kappa}^{\epsilon_2}  + \overline{\bbE}[\norme{\bfv}_{\mcE^\beta_T}]  ) | u_2^{\sharp} |_{1 + \kappa} 
\]
using first the Proposition \ref{prop_AAAAAA4} and estimate \ref{eq_estimeegu1bis}. The point is that this bound involves the relatively weak $| \cdot |_{1 + \kappa}$ norm so we can use interpolation
to bound it in terms of the $L^{\infty}$ and $| \cdot |_{2 \alpha}$ norms
\[ 
| u_2^{\sharp} |_{1 + \kappa} + | u_2^{\sharp} |_{1 + \kappa}^{1 +
   \epsilon_1} \lesssim | u_2^{\sharp} |_{L^{\infty}}^{\ell_0^-}  |
   u_2^{\sharp} |_{2 - \epsilon_1}^{\ell_0^+} + | u_2^{\sharp}
   |_{L^{\infty}}^{(1 + \epsilon_1) \ell_0^-}  | u_2^{\sharp} |_{2 -
   \epsilon_1}^{(1 + \epsilon_1) \ell_0^+} 
\]
with $\ell_0^- = \frac{1 - \epsilon_1 - \kappa}{2 - \epsilon_1}$ and $\ell_0^+ = \frac{1 + \kappa}{2 - \epsilon_1}$. Taking the supremum over $0 \leq t \leq
T$ we thus have
\[ 
| (g_u \odot \xi_{> n}) \nabla u_2^{\sharp} |_{\gamma, 0} \lesssim \Big( |
   u_2^{\sharp} |_{L^{\infty}}^{\ell_0^-}  | u_2^{\sharp} |_{2 -
   \epsilon_1}^{\ell_0^+} + | u_2^{\sharp} |_{L^{\infty}}^{(1 + \epsilon_1)
   \ell_0^-}  | u_2^{\sharp} |_{2 - \epsilon_1}^{(1 + \epsilon_1) \ell_0^+}\Big) \big(1+\overline{\bbE}[\norme{\bfv}_{\mcE^\beta_T}]\big) .
\]
For any $0 < \mu < 1$ we then get from Young inequality that
\[ 
| (g_u \odot \xi_{> n}) \nabla u_2^{\sharp} |_{\alpha / 2, 0} \lesssim \big(\mu
   \| u_2^{\sharp} \|_{\alpha / 2, 2 - \epsilon_1} + c (\mu) (| u_2^{\sharp}
   |^{q \ell_0^-}_{L^{\infty}_{t, x}} + | u_2^{\sharp} |^{q (1 + \epsilon_1)
   \ell_0^-}_{L^{\infty}_{t, x}}) \big) \big(1+\overline{\bbE}[\norme{\bfv}_{\mcE^\beta_T}]\big) 
\]
with $q = \frac{2 - \epsilon_1}{2 - \epsilon_1 - (1 + \epsilon_1) (1 + \kappa)}$, for some positive constant $c (\mu)$. Recalling the bounds $| u_2^{\sharp} |_{L^{\infty}_{t, x}} \lesssim \big(1 + \| u_2^{\sharp} \|^r_{\alpha / 2, 2 - \epsilon_1}\big) \big(1+\overline{\bbE}[\norme{\bfv}_{\mcE^\beta_T}]\big)$, we see that what matters is that $q(1 + \epsilon_1) \ell_0^- r < 1$. This strict inequality holds true when $\kappa = 0$ and $\epsilon_1 > 2 (1 - \alpha_{\circ})$ is small enough. It suffices then to assume from scratch that $\alpha_{\circ} < 1$ is close enough to $1$ to conclude that an appropriate choice of $\kappa > 0$ gives $q (1 + \epsilon_1) \ell_0^- r < 1$.
\end{Dem}

\ssk

Combining all the previous lemmas, we are finally ready to prove Theorem \ref{thm_nonexplosionA}. Indeed we obtain the bound for $\norme{u_2^\#}_{\alpha/2,2\alpha}$ from Lemma \ref{lem_AAAAAA4} and estimate \ref{eq_AAA25} by taking $\mu$ small enough. Then, we get the bound for $\norme{u}_{\alpha/2,\alpha}$ from Proposition \ref{prop_AAAAAA3}. 

}


\bigskip

\noindent \textcolor{gray}{$\bullet$} {\sf I. Bailleul} -- Univ Brest, CNRS, LMBA - UMR 6205, F-29238 Brest, France   \\
\noindent {\it E-mail}: ismael.bailleul@univ-brest.fr   

\bigskip

\noindent \textcolor{gray}{$\bullet$} {\sf N. Moench} -- Univ. Rennes, CNRS, IRMAR - UMR 6625, F-35000 Rennes, France   \\
\noindent {\it E-mail}: nicolas.moench@univ-rennes1.fr   

\bigskip

\noindent \textbf{\textsf{Acknowledgement --}} I.B. acknowledges some support from the ANR through the ANR grant ANR-22-CE40-0017.

\end{document}